\documentclass[11pt,reqno,twoside,makeidx]{amsart}
\usepackage[margin=2.7cm]{geometry}

\usepackage{caption}
\usepackage{graphicx}
\usepackage{dsfont}
\usepackage{amsthm}
\usepackage{amsmath}
\usepackage{amsfonts}
\usepackage{amscd}
\usepackage{fancyhdr}
\usepackage{enumerate}
\usepackage{enumitem}
\usepackage{url}
\usepackage{lipsum}
\usepackage{hyperref} 

\usepackage{tikz}
\usepackage{tikz-cd}
    \usetikzlibrary{cd} 
    \usetikzlibrary{arrows}
    \usetikzlibrary{matrix}
    \usepackage{pictexwd,dcpic}
    \usepackage{sidecap}
    \graphicspath{./imagine/}

\usepackage{caption}
\usepackage{subcaption}

\usepackage{cancel}
\usepackage{import}

\usepackage{amssymb}

\usepackage{mathtools}
\usepackage{layout}

\usepackage{array}
\newcolumntype{C}[1]{>{\centering\let\newline\\\arraybackslash\hspace{0pt}}m{#1}}

\usepackage{romannum} 
\AtBeginDocument{\pagenumbering{arabic}} 

\let\originalleft\left
\let\originalright\right
\renewcommand{\left}{\mathopen{}\mathclose\bgroup\originalleft}
\renewcommand{\right}{\aftergroup\egroup\originalright}

\newcommand{\fund}[1]{ \pi_1  \left( #1 \right)}
\newcommand{\compl}[1]{S^3 \setminus #1}
\newcommand{\cyclic}[1]{\nicefrac{\mathbb{Z}}{#1 \mathbb{Z}}}

\newcommand\restr[2]{{
  \left.\kern-\nulldelimiterspace 
  #1 
  \vphantom{|} 
  \right|_{#2} 
  }}
  
 \newcommand{\ZSU}{\mathcal{Z}\left(SU(2)\right)}
 \newcommand{\Ui}{U(1)}
 \newcommand{\normalsubgroup}[1]{\langle\! \langle #1 \rangle \! \rangle}
\newcommand{\R}{\mathcal{R}_{\Ui}}
\newcommand{\alg}[2]{\Delta \left( #1 ,#2\right)}

\tikzset{cross/.style={cross out, draw=black, minimum size=2*(#1-\pgflinewidth), inner sep=0pt, outer sep=0pt},
cross/.default={1pt}}

\DeclareMathOperator{\ima}{Im}
\DeclareMathOperator{\Tr}{Tr}
\DeclareMathOperator\coker{coker}

\usepackage{mathrsfs}
\newcommand{\calS}{\mathscr{S}}

\usepackage{multicol}
\usepackage{lscape}
\usepackage{cjhebrew}

\usepackage{bbold}

\usepackage{graphicx,nicefrac}
 \usepackage{framed}

\usepackage{algorithm}
\usepackage{algpseudocode}
\algtext*{EndFor}
\algtext*{EndIf}
\usepackage{tikz}
\usepackage{subcaption}
\makeatletter
\newtheorem*{rep@teo}{\rep@title}
\newcommand{\newreptheorem}[2]{%
\newenvironment{rep#1}[1]{%
 \def\rep@title{#2 \ref{##1}}%
 \begin{rep@teo}}%
 {\end{rep@teo}}}
 
 \newreptheorem{teo}{Theorem}
 \newreptheorem{cor}{Corollary}
 \newreptheorem{prop}{Proposition}
 \newreptheorem{defn}{Definition}
 \newreptheorem{lemma}{Lemma}
 \newreptheorem{condition}{Condition}

\theoremstyle{definition}
\newtheorem{defn}{Definition}[subsection]
\newtheorem{exmp}[defn]{Example}
\newtheorem{rmk}[defn]{Remark}

\newtheorem*{PN*}{Please Note}
\newtheorem{claim}{Claim}

\theoremstyle{plain}

\newtheorem{teo}[defn]{Theorem}
\newtheorem*{teo*}{Theorem}
\newtheorem{cor}[defn]{Corollary}
\newtheorem{lemma}[defn]{Lemma}
\newtheorem{prop}[defn]{Proposition}

\newtheorem*{prob*}{Problem}
\newtheorem{fact}[defn]{Fact}
\newtheorem*{notation*}{Notation}
\newtheorem{conj}[defn]{Conjecture}

\newtheorem{condition}{Condition}

\numberwithin{equation}{section}
\makeatother

\counterwithin*{defn}{section}
\renewcommand*{\thedefn}{
  \ifnum\value{subsection}<1 
    \thesection
  \else
    \thesubsection
  \fi
  .\arabic{defn}
}

\usepackage[
backend=biber,
style=alphabetic,
sorting=ynt,
natbib=true,
sorting=anyvt
]{biblatex}

\bibliography{biblio.bib}

\def\appendix{\par\c@section\z@ \c@subsection\z@
   \let\sectionname\appendixname
   \def\thesection{\@Alph\c@section}}

\usepackage{overpic}

\title[SU(2)-abelian graph manifolds with a single JSJ torus]{SU(2)-abelian graph manifolds with a single JSJ torus}
\author{Giacomo Bascape }

\begin{document}

\begin{abstract}
    A 3-manifold is called \emph{SU(2)}-abelian if every \emph{SU(2)}-representation of its fundamental group has abelian image. We classify, in terms of the Seifert coefficients, \emph{SU(2)}-abelian 3-manifolds among the family of graph manifolds obtained by gluing two Seifert spaces both fibred over a disk and with two singular fibers. Finally, we prove that these \emph{SU(2)}-abelian manifolds are Heegaard Floer homology L-spaces.
\end{abstract}
\maketitle
\section{Introduction}
This paper aims to classify which graph manifolds obtained by gluing together the boundaries of two Seifert fibred manifolds each with two singular fibers and disk base space are $SU(2)$-abelian.

The theory of $SU(2)$-representations provides a fruitful technique for studying the fundamental group of closed $3$-manifolds. In particular, $SU(2)$-abelian rational homology spheres are interesting from the instanton Floer homology perspective. For instance, in \cite{SivekBalwinSteinFilling} Baldwin and Sivek proved that an $SU(2)$-abelian rational homology sphere $Y$ is an instanton Floer L-space if and only if every $SU(2)$-representations of $\fund{Y}$ is non degenerate in a Morse-Bott sense.

In \cite{Menangerie} Sivek and Zentner provided the classification of closed $SU(2)$-abelian Seifert fibred spaces.
The first examples of $SU(2)$-abelian toroidal graph manifolds are contained in \cite{Motegi_haken_manifolds}, where it is proven that a specific gluing of torus knot complements is $SU(2)$-abelian. Here we provide more examples of $SU(2)$-abelian toroidal graph manifolds.

For $p_i \ge 2$, let $M_1=\mathbb{D}^2(\nicefrac{p_1}{q_1},\nicefrac{p_2}{q_2})$ and $M_2=\mathbb{D}^2(\nicefrac{p_3}{q_3},\nicefrac{p_4}{q_4})$ be Seifert fibred manifolds with torus boundary and singular fibers of orders $(p_1,p_2)$ and $(p_3,p_4)$ respectively (the notation is made explicit in Section \ref{section: notation}).
Denote by $M$ a graph manifold with a single JSJ torus $\Sigma \subset M$ such that $M_1$ and $M_2$ are the closures of the components of $M \setminus \Sigma$. 
We will write
$
M= M_1 \cup_{\Sigma} M_2.
$
In Theorem \ref{MAIN THEOREM} we give necessary and sufficient conditions for the manifold $M$ to be $SU(2)$-abelian.
Let $g_1=\gcd(p_1,p_2)$ and $g_2=\gcd(p_3,p_4)$, $o_1$ and $o_2$ be the orders of the rational longitudes $\lambda_{M_1}$ and $\lambda_{M_2}$, and lastly let $h_1$ and $h_2$ be the regular fibers of $M_1$ and $M_2$ in $\partial M_1 = \partial M_2 \subset M$. We denote by $\Delta(\gamma_1,\gamma_2)$ the geometric intersection number between the curves $\gamma_1$ and $\gamma_2$ on $\Sigma$.
\begin{teo}\label{MAIN THEOREM}
    Let $M= M_1 \cup_{\Sigma} M_2$ be a graph manifold as above and let us suppose that $g_1 \le g_2$. The manifold $M$ is $SU(2)$-abelian if and only if $\Delta(h_1,h_2)=1$ and each of the conditions $A$, $B$, and $C$ below are satisfied.
\end{teo}
\begin{condition}
One of the following holds:
\begin{itemize}
    \item $g_1=1$;
    \item $g_1=2$ and $o_1= 1$;
    \item $g_2=2$ and $o_2 = 1$;
    \item $g_1=2$, $o_2\equiv_2 1$, and $\Delta(\lambda_{M_2},h_1) \equiv_2 0$.
\end{itemize}
\end{condition}
\begin{condition}
One of the following holds:
    \begin{itemize}
        \item $\Delta(\lambda_{M_2},h_1)=0$, $p_1=2$, and $o_2\equiv_2 1$;
        \item $\Delta(\lambda_{M_2},h_1)=1$, and $o_2\le 2$;
        \item $\Delta(\lambda_{M_2},h_1)=4$, $p_1=2$, $p_2=4$, and $o_2=1$.
    \end{itemize}
\end{condition}
\begin{condition}
One of the following holds:
    \begin{itemize}
        \item $\Delta(\lambda_{M_1},h_2)=1$ and $o_1\le 2$;
        \item $\Delta(\lambda_{M_1},h_2)=2$, $p_3=p_4=4$, and $o_1=1$;
        \item $\Delta(\lambda_{M_1},h_2)=3$, $p_3=p_4=3$, $o_1=1$, and $\Delta(\lambda_{M_1},\lambda_{M_2})\equiv_2 0$;
        \item $\Delta(\lambda_{M_1},h_2)=4$, $p_3=2$, $p_4=4$, and $o_1=1$.
    \end{itemize}
\end{condition}

Theorem \ref{MAIN THEOREM} leads to a classification of $SU(2)$-abelian graph manifolds of the form $M= M_1 \cup_{\Sigma} M_2$ in terms of Seifert coefficients. Theorem \ref{teo: classification} states that there are only seven classes of $SU(2)$-abelian manifolds of the form $M= M_1 \cup_{\Sigma} M_2$; see Table $\ref{table: classification}$.
\begin{table}[b]
\begin{tabular}{p{0.5cm}p{1.8cm}p{2.7cm}p{4.5cm}C{2.7cm}}
$\#$ & $M_1$ & $M_2$ & Additional Requirements & $(\Delta_1,\Delta_2)$ \\[4pt]
\hline
\\ [-5pt]
1)&$\mathbb{D}^2\left( \frac{2}{1},\frac{p_2}{q_2}\right)$ & $\mathbb{D}^2\left( \frac{p_3}{q_3},\frac{p_3}{p_3-q_3}\right)$ &$p_3\equiv_2 1$ & $(0,1)$ \\[3pt]
&& & $2q_2 + p_2 \equiv_{2p_2} \pm o_1g_1$ \\[2pt]

2)&$\mathbb{D}^2\left( \frac{2}{1},\frac{p_2}{q_2}\right)$ & $\mathbb{D}^2\Big( \frac{3}{1},\frac{3}{2}\Big)$ &$o_1=1$ & $(0,3)$ \\[3pt]
 & & & $2q_2 + p_2 \equiv_{2p_2} \pm 3g_1$ \\[4pt]

3)& $\mathbb{D}^2\left( \frac{2}{1},\frac{4}{q_2}\right)$ & $\mathbb{D}^2\left( \frac{p_3}{q_3},\frac{p_4}{q_4}\right)$ &$o_2=1 $ & $(4,1)$ \\[3pt]
& & & $p_3q_4 + p_4q_3 \equiv_{p_3p_4} \pm 4g_2$\\[4pt]

4)& $\mathbb{D}^2\left( \frac{2}{1},\frac{4}{q_2}\right)$ & $\mathbb{D}^2\left( \frac{3}{q_3},\frac{3}{q_3}\right)$ & \text{\emph{none}}
& $(4,3)$ \\[8pt]

5)& $\mathbb{D}^2\left( \frac{p_1}{q_1},\frac{p_2}{q_2}\right)$ & $\mathbb{D}^2\left( \frac{3}{q_3},\frac{3}{q_3}\right)$ &$p_1 p_2 \equiv_2 1 $ & $(1,3)$ \\[3pt]
& & & $p_1q_2 + p_2q_1 \equiv_{p_1p_2} \pm 3$ \\[3pt]
& & & $g_1 =1$ \\[4pt]

6)& $\mathbb{D}^2\left( \frac{p_1}{q_1},\frac{p_2}{q_2}\right)$ & $\mathbb{D}^2\left( \frac{4}{q_3},\frac{4}{q_3}\right)$ &$o_1=1 $ & $(1,2)$ \\[3pt]
& & & $g_1 \le 2$ \\[3pt]
& & & $p_1q_2 +p_2q_1 \equiv_{p_1p_2} \pm 2g_1$\\[4pt]
 
 7)&$\mathbb{D}^2\left( \frac{p_1}{q_1},\frac{p_2}{q_2}\right)$ & $\mathbb{D}^2\left( \frac{p_3}{q_3},\frac{p_4}{q_4}\right)$ &Condition A & $(1,1)$ \\[3pt]
 & & & $o_1 \le 2$, $o_2 \le 2$ \\[3pt]
 & & & $(p_1 q_2 +p_2 q_1)\equiv_{p_1p_2} \pm o_1g_1$ \\[3pt]
 & & & $(p_3 q_4 +p_4 q_3)\equiv_{p_3p_4} \pm o_2g_2$\\[4pt]
\hline
\end{tabular}
\caption{Classes of $SU(2)$-abelian graph manifolds of the form $M_1 \cup_{\Sigma} M_2.$ We write $\Delta_1$ for $\Delta(\lambda_{M_2},h_1)$ and $\Delta_2$ for $\Delta(\lambda_{M_1},h_2)$.}
\label{table: classification}
\end{table}
\begin{teo}
    Let $M= M_1 \cup_{\Sigma} M_2$ be as before, and further suppose that $M_1$ and $M_2$ are presented in such a way that $0<q_i<p_i$ and $g_1 \le g_2$. The manifold $M$ is $SU(2)$-abelian if and only if $\Delta(h_1,h_2)=1$ and it is contained in one of the seven classes listed in Table \ref{table: classification}.
    \label{teo: classification}
\end{teo}

If $M_1$ and $M_2$ are as in Table \ref{table: classification}, Proposition \ref{prop: l'incollamento esiste sempre se soddisma le CC} provides a diffeomorphism $\varphi \colon \partial M_1 \to \partial M_2$ such that the graph manifold $M_1\cup_{\Sigma}M_2=M_1\cup_{\varphi}M_2$ is $SU(2)$-abelian. Thus, each class in Table \ref{table: classification} is nonempty.
Such a gluing is in general not unique, as shown in Example \ref{exmp: non unique diffeomorphism}.

Let $Y$ be the gluing of torus knot exteriors defined in \cite{Motegi_haken_manifolds}. In the latter Motegi proved that $Y$ is $SU(2)$-abelian.
\begin{cor}\label{cor: torus knots complements are unique}
    For $i \in \{1,2 \}$, let $E(T_i)$ be the exterior of a open tubular neighborhood of the torus knot $T_i \subset S^3$. We denote by $\lambda_i$ and $\mu_i$ the null-homologous longitude and the meridian of $T_i$. The manifold $E(T_1)\cup_{\Sigma} E(T_2)$ is $SU(2)$-abelian if and only if $\Delta(\lambda_1,\mu_2)=0$ and $\Delta(\lambda_2,\mu_1)=0$.
\end{cor}
Corollary \ref{cor: torus knots complements are unique} implies that if $Y$ is a graph manifold whose JSJ pieces are two torus knot exteriors, then $Y$ is one of the gluing of torus knot exteriors as in \cite{Motegi_haken_manifolds}. Since $Y$ has two JSJ pieces which are two torus knot exteriors, $Y$ belongs to class $7$ of Table \ref{table: classification}.

Let $Y$ be a rational homology sphere such that $\fund{Y}$ is \emph{cyclically finite} (see \cite{BoyerNicasAndrewVarietiesOfGroupsRepresentations}). If $Y$ is not an instanton L-space, then \cite[Theorem 4.6]{SivekBalwinSteinFilling} states that $\fund{Y}$ admits an irreducible $SU(2)$-representation. Moreover, it has been conjectured in \cite[Conjecture 7.24]{kronheimer2010knots} that the instanton Floer homology group of a rational homology $3$-sphere $Y$ is isomorphic, as a $\mathbb{C}$-vector space, to the $\mathbb{C}$-tensored Heegaard Floer homology group $\widehat{HF}(Y;\mathbb{C})$. Hence, if $Y$ is not a Heegaard Floer L-space, then $\fund{Y}$ is expected to admit an irreducible $SU(2)$-representation. This motivates the following conjecture.
\begin{conj}[\cite{zhang2019remarks}]
    If a rational homology sphere is $SU(2)$-abelian, then it is an L-space.
    \label{conj: SU(2)-abelian implies L-space}
\end{conj}

By means of an explicit calculation of the L-space interval of $M_1$ and $M_2$, and by applying the ``gluing theorem" for graph L-space manifolds (i.e. \cite[Proposition 1.5]{LspaceIntervalGraphManifold}) we prove the following:
\begin{teo}\label{teo: M SU(2)-abelian then Lspace}
    Let $M=M_1 \cup_{\Sigma} M_2$ be a $3$-manifold as in Theorem \ref{teo: classification}. If $M$ is an $SU(2)$-abelian rational homology sphere, then $M$ is an L-space.
\end{teo}
This gives further evidence to the conjecture that if a rational homology sphere is $SU(2)$-abelian, then it is a Heegaard Floer L-space (see Conjecture \ref{conj: SU(2)-abelian implies L-space}). 

We denote by $\Ui$ the subgroup of $SU(2)$ of diagonal matrices.
Let $Y$ be a $3$-manifold, and let $\Sigma \subset Y$ be an embedded torus. We define the \emph{$\Ui$-representation space of $\Sigma$ relative $Y$}, denoted by $T(Y,\Sigma)$, as the set of representations $\fund{\Sigma} \to \Ui$ that can be extended to a representation $\fund{Y}\to SU(2)$.
The strategy to prove Theorem \ref{MAIN THEOREM} is to analyse the spaces $T(M_1,\partial M_1)$ and $T(M_2,\partial M_2)$ and how they are related.

Let $Y_1$ and $Y_2$ be two $3$-manifolds with torus boundary and let $\varphi \colon \partial Y_1 \to \partial Y_2$ be a diffeomorphism. Let $Y=Y_1 \cup_{\varphi} Y_2$ be the $3$-manifold obtained by gluing $Y_1$ and $Y_2$ along $\varphi$ and let $\Sigma$ be the torus embedded in $Y$ corresponding to $\partial Y_1=\partial Y_2$. The spaces $T(Y_1,\Sigma)$ and $T(Y_2,\Sigma)$ both live in the same representation variety $\text{Hom}(\fund{\Sigma},\Ui)$. Therefore, we can consider their intersection in $\text{Hom}(\fund{\Sigma},\Ui)$.

\begin{repprop}{prop: Y SU(2)-abelian iff every point in the intersection contains SU(2)-abelian reprs}
    Let $Y=Y_1 \cup_{\varphi} Y_2$ be the $3$-manifold defined above. The manifold $Y$ is $SU(2)$-abelian if and only if every $SU(2)$-representation of $\fund{Y}$ that restricts to a representation in $T(Y_1,\Sigma) \cap T(Y_2,\Sigma)$ is $SU(2)$-abelian.
\end{repprop}

Let $Y=Y_1 \cup_{\varphi} Y_2$ be the manifold defined as before. We also determine an obstruction to the existence of an irreducible $SU(2)$-representation of $\fund{Y}$.

\begin{repcor}{cor: Y SU(2) allora i due pezzi sono SU(2)-abelian}
    If the manifold $Y = Y_1 \cup_{\varphi} Y_2$ is $SU(2)$-abelian, then the manifolds
    \[
       Y_1(\lambda_{Y_2}) \quad \text{and} \quad Y_2(\lambda_{Y_1})
    \]
    are $SU(2)$-abelian. Here $\lambda_{Y_1}$ and $\lambda_{Y_2}$ are the rational longitudes of $Y_1$ and $Y_2$.
\end{repcor}

An analogue of Corollary \ref{cor: Y SU(2) allora i due pezzi sono SU(2)-abelian} in the Heegaard Floer world is obtained as a consequence of \cite[Theorem 1.14]{hanselman2023bordered}.

\subsection*{Organization} In Section \ref{section: notation} we give a summary of the notation we shall use. In Section \ref{section: structure-less statements} we prove Proposition \ref{prop: Y SU(2)-abelian iff every point in the intersection contains SU(2)-abelian reprs} and Corollary \ref{cor: Y SU(2) allora i due pezzi sono SU(2)-abelian}. In Section \ref{section: the topology set up} we introduce the problem for the graph manifold $M=M_1 \cup_{\Sigma} M_2$ and study the $\Ui$-representation space of $\partial M_1$ relative $M_1$. In Section \ref{sec: P1 cap P2}, \ref{sec: A1 e H2}, and \ref{sec: H1 cap H2} we study the intersections of $T(M_1,\Sigma)$ with $T(M_2,\Sigma)$. In Section \ref{section: the main theorem}, we use the previous results for proving Theorem \ref{MAIN THEOREM} and the corresponding classification in Theorem \ref{teo: classification}. In Section \ref{sec: L-spaces} we prove Theorem \ref{teo: M SU(2)-abelian then Lspace}.

\subsection*{Acknowledgements}
I am deeply grateful to my research advisors, Duncan McCoy and Steven Boyer, for their constant patience, for guiding me along the path of my PhD and for the useful discussions during which many key ideas were discussed. I would also like to thank Antonio Alfieri for introducing me to gauge theory. Lastly, I would like to acknowledge my colleague and friend Patricia Sorya for constantly supporting me.

\section{Notation}
\label{section: notation}
\begin{defn}
A group $G$ is said to be \emph{$SU(2)$-abelian} (resp. \emph{$SU(2)$-central}) if every representation $G \to SU(2)$ has abelian image (resp. has image contained in the center $\ZSU=\{\pm 1\}$). A $3$-manifold is said to be \emph{$SU(2)$-abelian} if its fundamental group is $SU(2)$-abelian. 
\label{defn: su(2)-abelian 3 man}
\end{defn}

A representation $G \to SU(2)$ is said to be \emph{abelian} (resp. \emph{central}) if its image is abelian (resp. contained in the center $\ZSU$).
A non-abelian representation is also called \emph{irreducible}.

Let $z$ be a non-central element of $SU(2)$. We denote by $\Lambda_z \subset SU(2)$ the centralizer subgroup of the element $z$. The subgroups $\Lambda_z$ and $\Ui$ are known to be conjugate in $SU(2)$, details can be found in \cite[Lecture 13]{saveliev}. 
The following is well known.
\begin{fact}
Let $x,y$ be two non-central elements of $SU(2)$. The elements $x$ and $y$ commute if and only if the two centralizer subgroups $\Lambda_x$ and $\Lambda_y$ coincide.
\label{fact: two elements commute iff same centralizer}
\end{fact}
    
\begin{defn}\label{defn: angle that supports an interval}
    Let $I=(a_1,a_2) \subset [-2,2]$. Let $\theta_1,\theta_2 \in [0,\pi]$ such that $\Re (2e^{i \theta_j})= 2\cos(\theta_j) =a_j$. We say that the interval $I$ is \emph{supported} by the angle $\alpha(I)\coloneq |\theta_2 - \theta_1|$. See Figure \ref{Pic: Angle that supports}.
\end{defn}
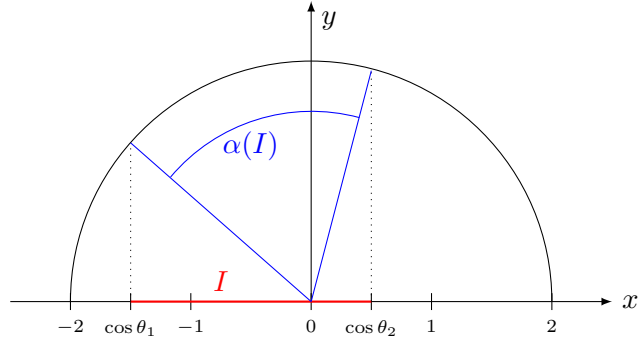
\begin{figure}[h]
 \centering
    \begin{tikzpicture}[>=latex,scale=1.6]
            \draw[->] (-2.5,0) -- (2.5,0) node[right] {$x$};
            \foreach \x /\n in {-2/$-2$,-1/$-1$,-1.5/$\cos\theta_1$ ,0/$0$,0.5/$\cos\theta_2$, 1/$1$,2/$2$} \draw[shift={(\x,0)}] (0pt,2pt) -- (0pt,-2pt) node[below] {\tiny \n};
            \draw[->] (0,0) -- (0,2.5) node[below right] {$y$};
            \draw (2,0) arc (0:180:2);
            \draw[red, thick] (0.5,0)--(-1.5,0);
            \draw [dotted] (-1.5,0)--(-1.5,2*0.66);
            \draw [dotted] (0.5,0)--(0.5,2*0.96);
            \draw[blue] (0,0)--(-1.5,2*0.66);
            \draw[blue] (0,0)--(0.5,2*0.96);
            \draw[blue] (0.4,1.53) arc (75:140:1.53);
            \node [above,red] at (-0.75,0) {$I$};
            \node [below,blue] at (-0.5,1.5) {$\alpha(I)$};
            \end{tikzpicture}
    \caption{The angle (in blue) that supports the interval $I$ (in red).}
    \label{Pic: Angle that supports}
\end{figure}

\begin{defn}
    Let $M$ be a closed and orientable $3$-manifold.
    A \emph{JSJ decomposition} on $M$ is a minimal collection of disjointly embedded incompressible tori $\mathcal{T} \subset M$ such that each connected component of $M \setminus \mathcal{T}$ is either atoroidal or Seifert fibered.
\end{defn} 
In \cite{Jaco1979SeifertFS}, it is proven that every closed $3$-manifold admits a JSJ decomposition. Moreover, they proved that a JSJ decomposition is unique up to isotopy.

\begin{defn}
    Let $\{\mu,\lambda\}$ be a basis of $H_1\left( \Sigma;\mathbb{Z} \right)$, where $\Sigma$ is a $2$-torus.
    Let $\gamma_1$ and $\gamma_2$ be two simple closed curves on a $2$-torus $\Sigma$ such that
    $[\gamma_1]= p_1 \mu + q_1 \lambda$ and $[\gamma_2]=p_2 \mu + q_2 \lambda$ in $H_1\left( \Sigma;\mathbb{Z} \right)$.
    The \emph{geometric intersection number} between $\gamma_1$ and $\gamma_2$ is $\Delta(\gamma_1,\gamma_2)\coloneqq |p_1q_2-q_2p_1|$.
\end{defn}

\subsection{Convention for Seifert fibred spaces}

We give \cite{martelli} as a reference. Let us suppose that $p_i \ge 2$.
Let $\{\mathbb{D}^2_i\}_{i \in \{1,\cdots, n\}}$ be a system of $n$ open disjoint disks in $S^2$. We define $\hat{Y}$ as
\[
    \hat{Y}= S^1 \times \left( S^2 \setminus  \coprod_{i=1}^n   \mathbb{D}^2_i \right)=S^1 \times S^2 \setminus  \left(\coprod_{i=1}^n S^1 \times  \mathbb{D}^2_i \right).
\]
We define $m_i$ and $l_i$ on the $i^{th}$-torus boundary component of $\hat{Y}$ as
\[
    m_i = \{\ast\} \times \partial \mathbb{D}^2_i \quad \text{and} \quad l_i=S^1 \times \{\ast\}.
\]
The closed manifold $Y=S^2(\nicefrac{p_1}{q_1},\cdots \nicefrac{p_n}{q_n})$ is the result of performing a Dehn filling of $\hat{Y}$ along the curve $p_i m_i + q_i l_i$ for every $i$. This construction gives the following presentation for the fundamental group of $Y$:
\begin{equation*}
  \fund{Y}=\fund{S^2\left(\nicefrac{p_1}{q_1},\cdots \nicefrac{p_n}{q_n}\right)}= \left\langle x_1,\cdots, x_n,h \,\middle| \, x_i^{p_i}h^{q_i}, [x_i,h], x_1x_2\cdots x_n\right\rangle.
  \label{eq: la uso una volta per il bordo}
\end{equation*}
Here $[x,y]$ denotes the commutator $xyx^{-1}y^{-1}$.
Let $h$ be a regular fiber of the Seifert fibred space $Y=S^2(\nicefrac{p_1}{q_1},\cdots \nicefrac{p_n}{q_n})$, we define the space $\mathbb{D}^2(\nicefrac{p_1}{q_1},\cdots \nicefrac{p_n}{q_n})$ as
\[
    \mathbb{D}^2(\nicefrac{p_1}{q_1},\cdots \nicefrac{p_n}{q_n}) = Y \setminus \nu(h),
\]
where $\nu(h)$ is a small fibred open neighborhood of the fiber $h$. If $M=\mathbb{D}^2(\nicefrac{p_1}{q_1},\cdots \nicefrac{p_n}{q_n})$, then its fundamental group is 
\begin{equation}
    \fund{M}=\fund{\mathbb{D}^2(\nicefrac{p_1}{q_1},\cdots \nicefrac{p_n}{q_n})} =\left\langle x_1,\cdots, x_n,h \,\middle| \, x_i^{p_i}h^{q_i}, [x_i,h]\right\rangle.
    \label{eq: definition of the seifert invariants}
\end{equation}

\begin{defn}\label{defn: meridian of the fibration}
Let $M=\mathbb{D}^2(\nicefrac{p_1}{q_1},\cdots \nicefrac{p_n}{q_n})$. Let us suppose that $\fund{M}$ is presented as in \eqref{eq: definition of the seifert invariants}. We define the \emph{fibration meridian} of this presentation as the, unique up to isotopy, simple closed curve $\mu \subset \partial M$, such that
\begin{align*}
    [\mu] = x_1x_2\cdots x_{n-1}x_n \in \fund{\mathbb{D}^2\left(\nicefrac{p_1}{q_1},\cdots ,\nicefrac{p_n}{q_n}\right)}.
\end{align*}
\end{defn}

The classification of Seifert fibrations \cite[Proposition 10.3.11]{martelli} states that there exists a fiber-preserving diffeomorphism
\begin{align}\label{eq: description of M1 in terms of fraction and disk. It says that q1 and q2 can be taken}
\mathbb{D}^2 \left( \frac{p_1}{q_1}, \cdots , \frac{p_n}{q_n}\right) \longrightarrow \mathbb{D}^2 \left( \frac{p_1}{ q_1},  \cdots , \frac{p_i}{kp_i+q_i}, \cdots ,\frac{p_n}{q_n}\right)
\end{align}
for every $k \in \mathbb{Z}$ and $i \in \{1,\cdots, n\}$.
In particular, we can suppose $q_i$ to be odd for every $i$.

Let $M=\mathbb{D}^2(\nicefrac{p_1}{q_1},\cdots \nicefrac{p_n}{q_n})$ with $n \ge 2$ and $p_i \ge 2$. The Seifert fibration on $M$ is unique up to isotopy unless $M \cong \mathbb{D}^2(\nicefrac{2}{1},\nicefrac{2}{1})$; see \cite[Proposition 10.4.16]{martelli}. The latter admits exactly two isotopy classes of Seifert fibration. One has base orbifold the M\"{o}bius band and the other has base orbifold the disk with two cone points. We will always consider the fibration over the disk for such a manifold.
\section{Some general results}
\label{section: structure-less statements}
In this section we do not require the studied $3$-manifolds to be compact, oriented, or fibred. 

\begin{defn}
    Let $\Sigma$ be a $2$-torus. We define $\R \left(\Sigma\right)$ as $\mbox{Hom}(\pi_1(\Sigma), \Ui$).
\end{defn}

Let $Y$ be a manifold, we denote by $\mathcal{R}(Y)$ the set $\text{Hom}(\fund{Y},SU(2))$. There exists a natural inclusion $\R (\Sigma) \subset \mathcal{R}(\Sigma)=\mbox{Hom}\left(\fund{\Sigma}, SU(2)\right)$.

\begin{defn}
    Let $Y$ be a $3$-manifold and $\Sigma$ an embedded torus via the map $\iota \colon \Sigma \to Y$. We define the \emph{$\Ui$-representation space of $\Sigma$ relative to $Y$}
    as the set of representations $\fund{\Sigma} \to \Ui$ that extend to a representation $\fund{Y} \to SU(2)$. We denote this by $T(Y,\Sigma) \subseteq \R(\Sigma)$.
\end{defn}
We can thus see $T(Y,\partial Y)$ as the $SU(2)$-version of the A-polynomial of $Y$, a formal definition of this polynomial can be found in \cite{PlaneCurvesAssociatedToCharacterVAriety}.
The set $T(Y,\partial Y)$ was previously defined in \cite[Definition 2.10]{TheAPolyHolonomyPerturbation} in terms of the flat connections of a trivialized rank $2$ unitary bundle over $Y$. 

Let $Y_1$ and $Y_2$ be a pair of $3$-manifolds with torus boundary and $\varphi \colon \partial Y_1 \to \partial Y_2$ a diffeomorphism, we denote by $Y= Y_1 \cup_{\varphi} Y_2$ be the manifold obtained by gluing $Y_1$ and $Y_2$ along the map $\varphi$.
Let $\Sigma$ be the torus $\partial Y_1 = \partial Y_2$ in $Y$ and $\iota \colon \Sigma \to Y$ the inclusion.
Then, the spaces $T(Y_1,\Sigma)$ and $T(Y_2,\Sigma)$ are both contained in $\R (\Sigma)$. Therefore, the intersection $T(Y_1,\Sigma) \cap T(Y_2,\Sigma)$ is well defined in $\R(\Sigma)$. If $\eta \colon \fund{\Sigma} \to \Ui$ is in $T(Y_1,\Sigma)\cap T(Y_2,\Sigma)$, then there exist $\rho_1 \in \mathcal{R}(Y_1)$ and $\rho_2 \in \mathcal{R}(Y_2,)$ such that
$
\restr{\rho_1}{\iota_\ast \fund{\Sigma}} \equiv \eta \equiv \restr{\rho_2}{\iota_\ast\fund{\Sigma}}.
$
Consequently, there exists a representation $\rho \colon \fund{Y} \to SU(2)$ such that $\restr{\rho}{\fund{Y_1}}\equiv \rho_1$, $ \restr{\rho}{\fund{Y_2}}\equiv\rho_2$, and $\restr{\rho}{\iota_\ast\fund{\Sigma}} \equiv \eta.$ 

\begin{prop}
    Let $Y=Y_1 \cup_{\varphi} Y_2$ be the $3$-manifold defined above. The manifold $Y$ is $SU(2)$-abelian if and only if every $SU(2)$-representation of $\fund{Y}$ that restricts to a representation in $T(Y_1,\Sigma) \cap T(Y_2,\Sigma)$ is $SU(2)$-abelian.
    \label{prop: Y SU(2)-abelian iff every point in the intersection contains SU(2)-abelian reprs}
    \begin{proof}
        If $Y$ is $SU(2)$-abelian, the conclusion is trivial.
        If $Y$ is not $SU(2)$-abelian, there exists a representation $\rho$ whose image is not abelian. Up to conjugation, we can suppose that $\rho(\iota_\ast\fund{\Sigma}) \subseteq \Ui$. The restriction $\restr{\rho}{\iota_\ast \fund{\Sigma}}$ extends to both $\fund{Y_1}$ and $\fund{Y_2}$. Hence, $\restr{\rho}{\iota_\ast\fund{\Sigma}} \in T(Y_1,\Sigma) \cap T(Y_2,\Sigma)$.
    \end{proof}
\end{prop}

\begin{defn} \label{defn: A1 H1 e P1}
    Let $Y$ be a $3$-manifold with torus boundary, let $\iota \colon \partial Y \to Y$ be the natural inclusion. We define the sets $A_Y, H_Y$, and $P_Y$ as:
    \begin{alignat*}{2}
        A_Y& \coloneqq \left\{ \eta \in \R \left(\partial Y \right) \,\middle|\, \exists \rho \in \mathcal{R}(Y) \text{ such that } \restr{\rho}{\iota_\ast \fund{\partial Y}}\equiv \eta\text{ and $\rho$ is abelian} \right\}, \\
        H_Y& \coloneqq \left\{ \eta \in \R \left(\partial Y \right) \,\middle|\, \exists \rho \in \mathcal{R}(Y) \text{ such that } \restr{\rho}{\iota_\ast \fund{\partial Y}}\equiv \eta\text{ and $\rho$ is irreducible} \right\}, \\
        P_Y& \coloneqq \Big\{ \eta \in \R \left(\partial Y \right) \,\Big|
            \begin{aligned}[t]
            &\exists \rho \in \mathcal{R}(Y) \text{ such that } \restr{\rho}{\iota_\ast \fund{\partial Y}}\equiv \eta,\\
            &\text{$\eta$ is central, and $\rho$ is abelian and non-central} \Big\}.
    \end{aligned}
    \end{alignat*}
\end{defn}
Equivalently, a representation $\fund{\partial Y} \to \Ui$ is in $A_Y$ (resp. $H_Y$) if and only if it extends to an abelian (resp. an irreducible) representation $\fund{Y} \to SU(2)$. Similarly, a representation $\fund{\partial Y} \to \ZSU$ is in $P_Y$ if and only if it extends to an abelian representation $\fund{Y} \to SU(2)$ whose image is not in $\ZSU$. 
Notice that $P_Y \subset A_Y$ and $A_Y \cup H_Y=T(Y,\partial Y) \subseteq \R(\partial Y)$. 

Let $Y=Y_1 \cup_{\varphi} Y_2$ and $\rho \in \mathcal{R}(Y)$. We write $\rho_1$ for $\restr{\rho}{\fund{Y_1}}$ and $\rho_2$ for $\restr{\rho}{\fund{Y_2}}$.

\begin{prop}\label{prop: if rho1 and rho2 abelian and rho non abelian implies P1 and P2}
    Let $Y=Y_1 \cup_{\varphi} Y_2$ and $\Sigma \subset Y$ the torus corresponding to $\partial Y_1=\partial Y_2$. If $\rho \in \mathcal{R}(Y)$ is an irreducible representation such that $\rho_1$ and $\rho_2$ are both abelian, then $\rho(\iota_\ast\fund{\Sigma}) \subseteq \ZSU$.
    \begin{proof}
        Since $\rho$ is irreducible, neither $\ima \rho_1$ nor $\ima \rho_2$ is central. Let us suppose, by contradiction, that there exists $g \in \iota_\ast\fund{\Sigma}$ such that $\rho_1(g)$ is not in the center. Since $g = \varphi_\ast(g)$ in $\fund{Y}$, we have that
        $
        \rho_1(g) = \rho_2(\varphi_\ast(g)).
        $
        In particular $\ima{\rho_1} \subseteq \Lambda_{\rho_1(g)} $ and $\ima{\rho_2} \subseteq \Lambda_{\rho_2(\varphi_\ast(g))}$, where $\Lambda_z \subset SU(2)$ is the centralizer subgroup of the non-central element $z \in SU(2)$.
        Fact \ref{fact: two elements commute iff same centralizer} implies that $\Lambda_{\rho_1(g)}=\Lambda_{\rho_2(\varphi_\ast(g))}$. This implies that $\ima \rho \subseteq \Lambda_{\rho_1(g)}=\Lambda_{\rho_2(\varphi_\ast(g))}$, which contradicts the irreducibility of the representation $\rho$.
    \end{proof}
\end{prop}

\begin{prop}
    Let $Y=Y_1 \cup_{\varphi} Y_2$ be as above. If there exists $\rho \in \mathcal{R}(Y)$ such that $\rho(\iota_\ast \fund{\Sigma}) \subseteq \ZSU$ and neither $\rho_1$ nor $\rho_2$ is central, then the manifold $Y$ is not $SU(2)$-abelian.
    \label{prop: if the image of the torus is central the either p1 or p2 has to be central}
    \begin{proof}
        If $\rho$ is irreducible we get our conclusion.
        Let us suppose that $\rho$ is reducible and therefore both $\rho_1$ and $\rho_2$ are reducible as well. Without loss of generality, we assume that the images $\ima \rho_1$ and $\ima \rho_2$ are both in $\Ui$. Let $z \in SU(2) \setminus \Ui$. A non-central element in $z \Ui z^{-1}$ does not commute with a non-central element of $\Ui$. Let $ \gamma \colon \fund{Y} \to SU(2)$ be defined as
        \[
            \gamma(x) = 
            \begin{dcases}
                \rho_1(x) & \text{if } x \in \fund{Y_1} \\ 
                z\rho_2(x)z^{-1} & \text{if } x \in \fund{Y_2}
            \end{dcases}.
        \]
        Since neither $\rho_1$ nor $\rho_2$ is central, the representation $\gamma$ has non-abelian image. This implies that $\gamma$ is irreducible.
    \end{proof}
\end{prop}

Let $Y_1$ and $Y_2$ be two $3$-manifolds with torus boundary and let $Y=Y_1 \cup_{\varphi} Y_2$. Let $\Sigma \subset Y$ be the torus corresponding to $\partial Y_1=\partial Y_2$, let $\iota \colon \Sigma \to Y$ be inclusion. In order to make the notation a little lighter, the sets 
$A_{Y_i}$, $H_{Y_i}$, and $P_{Y_i}$ of Definition \ref{defn: A1 H1 e P1} will be denoted by $A_i$, $H_i$, and $P_i$, where $i \in \{1,2\}$.
Proposition \ref{prop: if rho1 and rho2 abelian and rho non abelian implies P1 and P2} and Proposition \ref{prop: if the image of the torus is central the either p1 or p2 has to be central} imply that there exists an irreducible representation $\rho \in \mathcal{R}(Y)$ such that $\restr{\rho}{\iota_\ast \fund{\Sigma}} \in A_1 \cap A_2$ if and only if $\restr{\rho}{\iota_\ast \fund{\Sigma}} \in P_1 \cap P_2$.

\begin{teo}\label{teo: M SU(2)-abeliano se e solo se i pezzi sono empty}
    Let $Y_1$ and $Y_2$ be two $3$-manifolds with torus boundary. The manifold $Y=Y_1 \cup_{\varphi} Y_2$ is $SU(2)$-abelian if and only if $H_1\cap H_2$, $H_1 \cap A_2$, $A_1 \cap H_2$, and $P_1 \cap P_2$ are empty.
    \begin{proof}
        If $Y$ is $SU(2)$-abelian, then, for every $\rho \in \mathcal{R}(Y)$, the restrictions $\rho_1=\restr{\rho}{\iota_\ast\fund{Y_1}}$ and $\rho_2=\restr{\rho}{\iota_\ast\fund{Y_2}}$ are abelian. 
        This implies that $H_1\cap H_2 = \emptyset$, $A_1\cap H_2 = \emptyset$, and $H_1\cap A_2 = \emptyset$.
        
        Let $\eta \in P_1 \cap P_2$, then there exists a representation $\rho \in \mathcal{R}(Y)$ such that $\rho_1$ and $\rho_2$ are non-central. Proposition \ref{prop: if the image of the torus is central the either p1 or p2 has to be central} implies that $Y$ is not $SU(2)$-abelian, which is a contradiction. Thus, $P_1 \cap P_2 = \emptyset$.
        
        Conversely, let $\rho \in \mathcal{R}(Y)$ be an irreducible representation. In particular, $Y$ is not $SU(2)$-abelian. Up to conjugation, we can suppose that $\rho(\iota_\ast \fund{\Sigma}) \subseteq \Ui$. 
        If $\restr{\rho}{\iota_\ast\fund{\Sigma}}$ is either in $H_1\cap H_2 $, $A_1\cap H_2 $, or in $H_1\cap A_2$, then we get our conclusion. If $\restr{\rho}{\iota_\ast \fund{\Sigma}}$ is neither in $H_1\cap H_2 $, $A_1\cap H_2 $, nor in $H_1\cap A_2$, then $\rho_1$ and $\rho_2$ are both abelian. Proposition \ref{prop: if rho1 and rho2 abelian and rho non abelian implies P1 and P2} implies that $\rho(\iota_\ast \fund{\Sigma}) \subseteq \ZSU$. Since $\rho$ is irreducible and the restrictions $\rho_1$ and $\rho_2$ are $SU(2)$-abelian, neither $\rho_1$ nor $\rho_2$ is $SU(2)$-central. This implies that $\restr{\rho}{\iota_\ast \fund{\Sigma}}$ is in $P_1 \cap P_2$. Therefore, $P_1 \cap P_2$ is nonempty.
    \end{proof}
\end{teo}

Let $Y$ be a compact orientable $3$-manifold with torus boundary, we identify the group $H_1(\partial Y;\mathbb{Z}) \cong \mathbb{Z}^2$ with the group $\fund{\partial Y} \cong \mathbb{Z}^2$ in the natural way. With an abuse of notation, we consider the group 
\[
\ker \left(\iota_\ast: H_1(\partial Y ; \mathbb{Z}) \to H_1( Y ; \mathbb{Z})\right) \le H_1(\partial Y ; \mathbb{Z}),
\]
where $\iota \colon \partial Y \to Y$ is the inclusion,
as a subgroup of $\fund{\partial Y}$.
\begin{prop}\label{prop: a representation from the boundary to the whole manifold extends iff is trivial on the kernel of iota}
    Let $Y$ be a $3$-manifold with torus boundary and $\iota \colon \partial Y \to Y$ the natural inclusion. A representation $\eta \colon \fund{\partial Y} \to SU(2)$ extends to an abelian representation $\rho:\fund{Y} \to SU(2)$ if and only if $\restr{\eta}{\ker{\iota_\ast}} \equiv 1$.
    \begin{proof}
        The abelianization map $\fund{Y} \twoheadrightarrow H_1(Y;\mathbb{Z})$ is denoted by $\mathcal{A}$. Let us suppose that the representation $\eta \colon\fund{\partial Y} \to SU(2)$ extends to an abelian representation $\rho:\fund{Y}\to SU(2)$. That is, $\eta = \rho\circ \iota_\ast$. Since $\rho$ is abelian, there exists a representation $\widetilde{\rho} \colon H_1(Y;\mathbb{Z}) \to SU(2)$, such that $\rho=\widetilde{\rho}\circ \mathcal{A}$. This implies that all the triangles in diagram \eqref{eq: diagram theorem delle cose abeliane} commute. 
       \begin{equation}
            \begin{tikzcd}
             H_1(Y;\mathbb{Z}) \arrow[rrrd, bend left, "\widetilde{\rho}"]\\
             & \fund{Y} \arrow[ul, twoheadrightarrow, "\mathcal{A}"'] \arrow[rr, "\rho"]
             & & SU(2). \\
            H_1(\partial Y;\mathbb{Z}) \cong \fund{\partial Y} \arrow[uu, "\iota_\ast"] \arrow[ur, "\iota_\ast"'] \arrow[rrru, bend right, "\eta"']
            \end{tikzcd}
           \label{eq: diagram theorem delle cose abeliane}
       \end{equation}
       Thus, $\eta \equiv \widetilde{\rho} \circ \iota_\ast$. Therefore,
        $
            \eta(\ker \iota_\ast)= \widetilde{\rho} \circ \iota_\ast \left( \ker \iota_\ast\right)=1.
        $
        This concludes one direction.

       Conversely, let $\eta \colon \fund{\partial Y} \to SU(2)$ be a representation that is trivial on $\ker \iota_\ast$. Up to conjugating, we can suppose that $\eta$ has image in $\Ui \subset SU(2)$. Therefore, we define the (abelian) representation
       $\gamma \colon \ima \iota_\ast \to \Ui$ as $\eta(x)=\gamma(x)$. Since group $\Ui$ is divisible, \cite[Proposition 3.19]{Lam_1999} implies that there exists a representation $\widetilde{\rho} \colon H_1(Y;\mathbb{Z})\to \Ui$ that extends $\gamma$. The representation $\rho$ is given by pre-composing $\widetilde{\rho}$ with the abelianization homomorphism $\mathcal{A}$.
    \end{proof}
\end{prop}

Let $Y$ be a $3$-manifold with torus boundary and let $\iota\colon \partial Y \to Y$ be the natural inclusion. According to the ``half lives and half dies" Theorem \cite[Corollary 9.1.5]{martelli}, we have that
\[
     \dim \left(\ker \left( H_1(\partial Y; \mathbb{Q}) \overset{\iota_\ast}{\longrightarrow} H_1(Y;\mathbb{Q})\right)\right)=1.
\]
We define the \emph{rational longitude} as the unique, up to isotopy, simple and closed curve $\lambda_{Y} \subset \partial Y$ such that its class in homology generates $\ker\iota_\ast \subseteq H_1(\partial Y;\mathbb{Q})$.
Equivalently, the rational longitude of $Y$ is a simple and closed curve $\lambda_{Y} \subset \partial Y$ such that its class is a torsion element of $H_1(Y; \mathbb{Z})$, which is unique up to isotopy.
We want to describe
$\ker \iota_\ast \subseteq H_1(\partial Y; \mathbb{Z})$ in terms of the rational longitude $\lambda_Y$ and obtain a more operational formulation of Proposition \ref{prop: a representation from the boundary to the whole manifold extends iff is trivial on the kernel of iota}.

In this paper an abuse of notation is in use: for a given simple closed curve $\gamma$ in the torus $\Sigma$, when we refer to its homotopy class $[\gamma] \in \fund{\Sigma}$ we omit the brackets. Consequently, $\gamma$ indicates both a curve in $\Sigma$ and its homotopy class $\gamma \in \fund{\Sigma}$.

\begin{cor}\label{cor: an abelian representation extends if and only if the rational longitude is zero}
    Let $Y$ be a $3$-manifold with torus boundary and let $\lambda_Y$ be its rational longitude. Let $n$ be the order of $\lambda_Y$ in $H_1(Y;\mathbb{Z})$. A representation $\eta \colon \fund{\partial Y}\to SU(2)$ extends to an abelian representation $\fund{Y}\to SU(2)$ if and only if $\eta(\lambda_{Y})^{n}=1$.
    \begin{proof}
        The subgroup $\ker \iota_\ast \le H_1( \partial Y;\mathbb{Z})$ is generated by the element $n \cdot \lambda_{Y} \in H_1(\partial Y;\mathbb{Z})$. Hence, $\eta(\ker \iota_\ast)=1$ if and only if $\eta(\lambda_{Y})^n=1$. The conlcusion holds by Proposition \ref{prop: a representation from the boundary to the whole manifold extends iff is trivial on the kernel of iota} 
    \end{proof}
\end{cor}

In what follows $\lambda_{Y_1}$ and $\lambda_{Y_2}$ are the rational longitudes of $Y_1$ and $Y_2$. Furthermore, $o_1$ and $o_2$ are the orders of the corresponding rational longitudes in $H_1(Y_1;\mathbb{Z})$ and $H_1(Y_2;\mathbb{Z})$. 

\begin{prop}
    Let $Y_1$ and $Y_2$ be two $3$-manifolds with torus boundary and let $\varphi \colon \partial Y_1 \to \partial Y_2$ be a diffeomorphism. If the manifold $Y = Y_1 \cup_{\varphi} Y_2$ is $SU(2)$-abelian, then the groups
    \[
        \frac{\fund{Y_1}}{\normalsubgroup{{\lambda_{Y_2}}^{o_2}}} \quad \text{and} \quad  \frac{\fund{Y_2}}{\normalsubgroup{{\lambda_{Y_1}}^{o_1}}}
    \]
    are $SU(2)$-abelian. Where $\normalsubgroup{\lambda}$ denotes the smallest normal subgroup containing $\lambda$.
    \label{prop: Y SU(2)-abelian then il pezzo orbifold e SU(2)-abelian}
    \begin{proof}
        Let us suppose that the group $\fund{Y_1}/\normalsubgroup{{\lambda_{Y_2}}^{o_2}}$ admits an irreducible $SU(2)$-representation. Hence, there exists an irreducible representation $\rho_1\colon \fund{Y_1} \to SU(2)$ such that $\rho_1({\lambda_{Y_2}}^{o_2})=1$. Let $\Sigma$ be the embedded torus in $Y$ corresponding to $\partial Y_1 =\partial Y_2$ and let $\iota \colon \Sigma \to Y$ be the natural inclusion. Let $\eta$ be the restriction $\restr{\rho_1}{\iota_\ast \fund{\Sigma}}$. Since $\eta(\lambda_{Y_2})^{o_2}=\rho_1(\lambda_{Y_2})^{o_2}=1$, Corollary \ref{cor: an abelian representation extends if and only if the rational longitude is zero} implies that the representation $\eta$ extends to an abelian representation $\rho_2 \colon \fund{Y_2}\to SU(2)$. Thus, there exists a representation $\rho \colon \fund{Y} \to SU(2)$ such that $\restr{\rho}{\fund{Y_1}}\equiv \rho_1$ and $\restr{\rho}{\fund{Y_2}}\equiv \rho_2$. The representation $\rho$ is irreducible. The conclusion for $\fund{Y_1}/\normalsubgroup{{\lambda_{Y_2}}^{o_2}}$ holds similarly. 
    \end{proof}
\end{prop}

\begin{cor}
    If $Y$ is $SU(2)$-abelian, then the manifolds $Y_1(\lambda_{Y_2})$ and $Y_2(\lambda_{Y_1})$ are $SU(2)$-abelian as well.
    \label{cor: Y SU(2) allora i due pezzi sono SU(2)-abelian}
    \begin{proof}
        The group $\frac{\fund{Y_1}}{\normalsubgroup{{\lambda_{Y_2}}^{o_2}}}$ surjects onto $\fund{Y_1(\lambda_{Y_2})}$.
        The conclusion holds by Proposition \ref{prop: Y SU(2)-abelian then il pezzo orbifold e SU(2)-abelian}.
    \end{proof}
\end{cor}
\section{The graph manifold \texorpdfstring{$M_1 \cup_{\Sigma} M_2$}{M\_1  U M\_2}}
\label{section: the topology set up}
From now on, we will only work with compact $3$-manifolds. Let $M_1$ and $M_2$ be two Seifert manifolds with torus boundary, each with base orbifold a disk with two singular fibers. Let $M$ be a graph manifold with a single JSJ torus $\Sigma \subset M$ such that $M_1$ and $M_2$ are the closures of the components of $M \setminus \Sigma$. We will write
\[
M= M_1 \cup_{\Sigma} M_2.
\]
The orientation of $M$ gives an orientation to the manifolds $M_1$ and $M_2$. Let $\varphi \colon \partial M_1 \to \partial M_2$ be a diffeomorphism such that
\begin{equation}
M= M_1 \cup_{\Sigma} M_2= M_1 \cup_{\varphi} M_2.
\label{eq: the manifold M}
\end{equation}
Since $M$ is an oriented manifold, the diffeomorphism $\varphi$ is an orientation reversing diffeomorphism. We will determine which of these manifolds are $SU(2)$-abelian.

Let us suppose that the manifolds $M_1$ and $M_2$ are parameterized as
\begin{align}
M_1 = \mathbb{D}^2\left(\frac{p_1}{q_1},\frac{p_2}{q_2}\right), \quad M_2 = \mathbb{D}^2\left(\frac{p_3}{q_3},\frac{p_4}{q_4}\right);
\label{the two seifert pieces}
\end{align}
with $\gcd(p_i,q_i)=1$, $p_2 \ge p_1 \ge 2$, and $p_4 \ge p_3 \ge 2$ (see Section \ref{section: notation} for the notation).
Let $g_1\coloneqq \gcd(p_1,p_2)$ and $g_2 \coloneqq \gcd(p_3,p_4)$. Without loss of generality, we assume $g_1\le g_2$.
The presentations \eqref{the two seifert pieces} give the following presentations for their fundamental groups:
\begin{align}
    \begin{split}
        \fund{M_1}=\langle a_{1},b_{1},h_1 \,|\, a_1^{p_{1}}h_1^{q_{1}}, b_1^{p_2}h_1^{q_{2}},[h_1,a_1],[h_1,b_1]\rangle \quad \text{and} \\[4pt]
        \fund{M_2}=\langle a_{2},b_{2},h_2 \,|\, a_1^{p_{3}}h_2^{q_{3}}, b_2^{p_4}h_2^{q_{4}},[h_2,a_2],[h_2,b_2]\rangle;
    \end{split}
    \label{eq: presentation of the fundamental group of M1 and M2}
\end{align}
where $[x,y]\coloneqq xyx^{-1}y^{-1}$. 

Given these presentations, the curves $\mu_1 \subset \partial M_1$ and $\mu_2 \subset \partial M_2$ will represent the fibration meridians of $M_1$ and $M_2$ as in Definition \ref{defn: meridian of the fibration}, respectively. In particular, when $\mu_i$ is considered as an element of $\fund{M_i}$, then $\mu_i = a_ib_i$, for $i \in \{1,2\}$. The groups $\fund{\partial M_1}$ and $\fund{\partial M_2}$ admit the following presentations:
\begin{align*}
    \fund{\partial M_1}=\langle \mu_1,h_1 \,|\, [h_1,\mu_1]\rangle,\quad \text{and} \quad
    \fund{\partial M_2}=\langle \mu_2,h_2 \,|\, [h_2,\mu_2]\rangle.
\end{align*}

We set the convention that the vectors $\left(\begin{smallmatrix}1\\0\end{smallmatrix}\right),\left(\begin{smallmatrix}0\\1\end{smallmatrix}\right)\in \mathbb{Z}^2 \cong \fund{\partial M_1}$ correspond to $\mu_1$ and $h_1$ respectively. Similarly, the vectors $\left(\begin{smallmatrix}1\\0\end{smallmatrix}\right),\left(\begin{smallmatrix}0\\1\end{smallmatrix}\right)\in \mathbb{Z}^2 \cong \fund{\partial M_2}$ correspond to $\mu_2$ and $h_1$.
We will use the matrix
\[
    \begin{bmatrix}
     \alpha & \beta \\ \gamma & \delta
    \end{bmatrix}  \quad \text{with} \quad \alpha\delta - \beta\gamma = -1
\]
to represent the map $\varphi_* \colon \fund{\partial M_1}\to \fund{\partial M_2}$ with respect to these ordered bases.
In particular, we have that $\Delta(h_1,h_2)=|\beta|$.
This gives an explicit presentation for the fundamental group of $M=M_1 \cup_{\varphi}M_2$:
\begin{equation}
    \fund{M} = \frac{\fund{M_1}\ast \fund{M_2}}{\normalsubgroup{a_1b_1=(a_2b_2)^\alpha h_1^\gamma, h_1 = (a_2b_2)^\beta h_2 ^\delta }}
    \label{fund(M) descrizione totale}
\end{equation}
If the diffeomorphism $\varphi :\partial M_1 \to \partial M_2$ is such that $\Delta(h_1,h_2)=\beta=0$, then the matrix $\varphi_*$ is lower triangular.
Thus, we have that $\varphi(h_1)=h_2$. This means
the fibration of $M_1$ can be extended over $M_2$. Therefore, the manifold $M=M_1 \cup_{\varphi} M_2$ is a closed Seifert fibred manifold. In particular, $M$ has four singular fibers and it is fibred over $S^2$. Therefore, $M$ is not $SU(2)$-abelian by \cite[Theorem 1.2]{Menangerie}. Thus, we assume that $|\beta| \ge 1$.

\begin{lemma}\label{prop: one between h1 and h2 has to go in the center}
If $\fund{M}$ admits an irreducible $SU(2)$-representation $\rho$, then either $\rho(h_1) \in \ZSU$ or $\rho(h_2)\in \ZSU$.
\begin{proof}
We recall that $\Sigma=\partial M_1 =\partial M_2$. Since the torus $\Sigma$ contains the regular fibers $h_1$ and $h_2$, their homotopy classes commute in $\fund{M}$.

Let us suppose by contradiction that $\rho(h_1) \notin \ZSU$ and $\rho(h_2) \notin \ZSU$. According to Fact \ref{fact: two elements commute iff same centralizer}, the centralizers $\Lambda_{\rho(h_1)}$ and $\Lambda_{\rho(h_2)}$ coincide. Since $[a_1,h_1]=[b_1,h_1]=1$ and $[a_2,h_2]=[b_2,h_2]=1$, we conclude that $\rho(\fund{M_1})\subset\Lambda_{\rho(h_1)}$ and $\rho(\fund{M_2})\subset\Lambda_{\rho(h_2)}$. Then, the elements $\rho(a_1)$, $\rho(b_1)$, $\rho(h_1)$, $\rho(a_2)$, $\rho(b_2)$, and $\rho(h_2)$ are all contained in
$\Lambda_{\rho(h_1)}=\Lambda_{\rho(h_2)}$.
This implies that $\rho$ is abelian, which is a contradiction.
\end{proof}
\end{lemma}

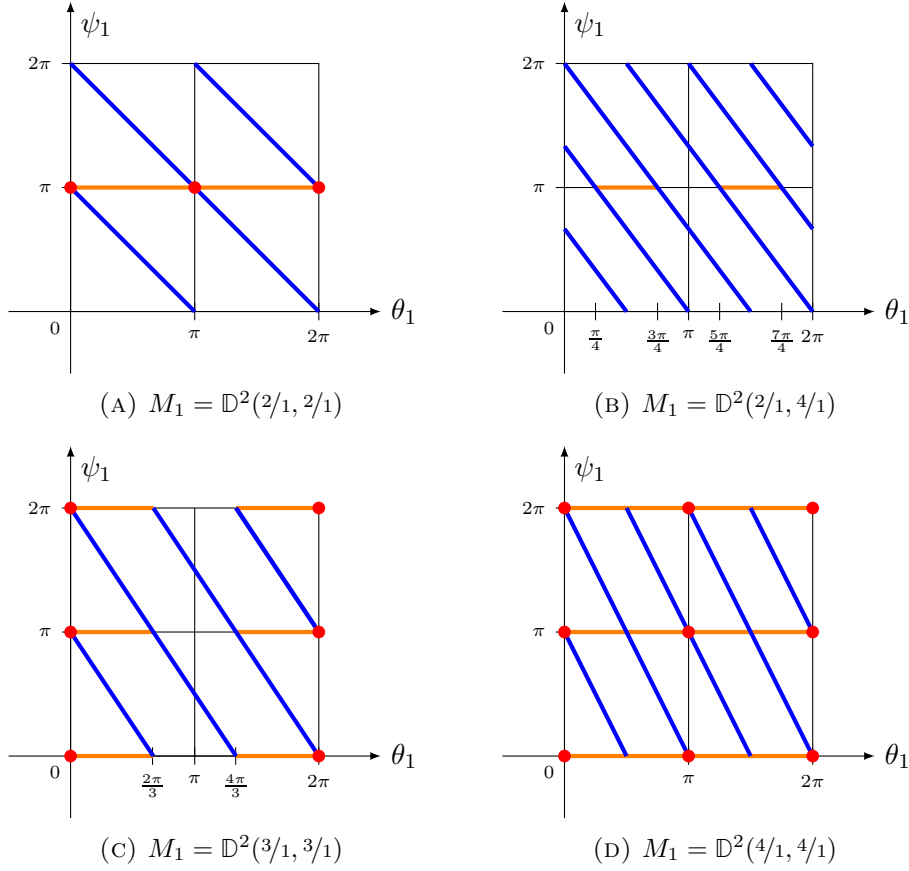
\begin{figure}[t]
    \centering
    \begin{subfigure}[b]{0.35\textwidth}
        \begin{tikzpicture}[>=latex,scale=1.65]
            \draw[->] (-.5,0) -- (2.5,0) node[right] {$\theta_1$};
            \foreach \x /\n in {1/$\pi$,2/$2\pi$} \draw[shift={(\x,0)}] (0pt,2pt) -- (0pt,-2pt) node[below] {\tiny \n};
            \draw[->] (0,-.5) -- (0,2.5) node[below right] {$\psi_1$};
            \foreach \y /\n in {1/$\pi$,2/$2\pi$}
            \draw[shift={(0,\y)}] (2pt,0pt) -- (-2pt,0pt) node[left] {\tiny \n};
            \node[below left] at (0,0) {\tiny $0$};
             \draw (0,0) rectangle (2,2);
             \draw (1,0) -- (1,2); 
             \draw (0,1) -- (2,1);
             \draw [orange,  ultra thick] (0,1) -- (2,1);
             \draw [blue,  ultra thick] (0,2) -- (2,0);
             \draw [blue,  ultra thick] (1,2) -- (2,1);
             \draw [blue,  ultra thick] (0,1) -- (1,0);
             \node at (0,1) [circle,fill,inner sep=1.7pt, red]{};
             \node at (1,1) [circle,fill,inner sep=1.7pt, red]{};
             \node at (2,1) [circle,fill,inner sep=1.7pt, red]{};
        \end{tikzpicture}
        \caption{$M_1=\mathbb{D}^2(\nicefrac{2}{1},\nicefrac{2}{1})$}
    \end{subfigure}
    \qquad 
    \begin{subfigure}[b]{0.35\textwidth}
        \begin{tikzpicture}[>=latex,scale=1.65]
        \draw[->] (-.5,0) -- (2.5,0) node[right] {$\theta_1$};
        \foreach \x /\n in {1/$\pi$,2/$2\pi$, 0.25/$\frac{\pi}{4}$, 0.75/$\frac{3\pi}{4}$, 1.25/$\frac{5\pi}{4}$, 1.75/$\frac{7\pi}{4}$} \draw[shift={(\x,0)}] (0pt,2pt) -- (0pt,-2pt) node[below] {\tiny \n};
        \draw[->] (0,-.5) -- (0,2.5) node[below right] {$\psi_1$};
        \foreach \y /\n in {1/$\pi$,2/$2\pi$}
        \draw[shift={(0,\y)}] (2pt,0pt) -- (-2pt,0pt) node[left] {\tiny \n};
        \node[below left] at (0,0) {\tiny $0$};
         \draw (0,0) rectangle (2,2);
         \draw (1,0) -- (1,2); 
         \draw (0,1) -- (2,1);
         \draw [orange, ultra thick] (0.25,1) -- (0.75,1);
         \draw [orange, ultra thick] (1.25,1) -- (1.75,1);
         \draw [blue,  ultra thick] (0,2) -- (1.5,0);
         \draw [blue,  ultra thick] (1.5,2) -- (2,4/3);
         \draw [blue,  ultra thick] (0,4/3) -- (1,0);
         \draw [blue,  ultra thick] (1,2) -- (2,2/3);
         \draw [blue,  ultra thick] (0,2/3) -- (0.5,0);
         \draw [blue,  ultra thick] (0.5,2) -- (2,0);
    \end{tikzpicture}
        \caption{$M_1=\mathbb{D}^2(\nicefrac{2}{1},\nicefrac{4}{1})$}
    \end{subfigure}
    \\[8pt]
    \begin{subfigure}[b]{0.35\textwidth}
        \begin{tikzpicture}[>=latex,scale=1.65]
        \draw[->] (-.5,0) -- (2.5,0) node[right] {$\theta_1$};
        \foreach \x /\n in {1/$\pi$,2/$2\pi$, 0.66/$\frac{2\pi}{3}$, 1.33/$\frac{4\pi}{3}$} \draw[shift={(\x,0)}] (0pt,2pt) -- (0pt,-2pt) node[below] {\tiny \n};
        \draw[->] (0,-.5) -- (0,2.5) node[below right] {$\psi_1$};
        \foreach \y /\n in {1/$\pi$,2/$2\pi$}
        \draw[shift={(0,\y)}] (2pt,0pt) -- (-2pt,0pt) node[left] {\tiny \n};
        \node[below left] at (0,0) {\tiny $0$};
         \draw (0,0) rectangle (2,2);
         \draw (1,0) -- (1,2); 
         \draw (0,1) -- (2,1);
         \draw [orange, ultra thick] (0,0) -- (2/3,0);
         \draw [orange, ultra thick] (4/3,0) -- (2,0);
         \draw [orange, ultra thick] (0,1) -- (2/3,1);
         \draw [orange, ultra thick] (4/3,1) -- (2,1);
         \draw [orange, ultra thick] (0,2) -- (2/3,2);
         \draw [orange, ultra thick] (4/3,2) -- (2,2);
         \draw [blue,  ultra thick] (0,2) -- (4/3,0);
         \draw [blue,  ultra thick] (4/3,2) -- (2,1);
         \draw [blue,  ultra thick] (4/3,2) -- (2,1);
         \draw [blue,  ultra thick] (0,1) -- (2/3,0);
         \draw [blue,  ultra thick] (2/3,2) -- (2,0);
        
         \node at (0,0) [circle,fill,inner sep=1.7pt, red]{};
         \node at (0,1) [circle,fill,inner sep=1.7pt, red]{};
         \node at (0,2) [circle,fill,inner sep=1.7pt, red]{};
         \node at (2,0) [circle,fill,inner sep=1.7pt, red]{};
         \node at (2,1) [circle,fill,inner sep=1.7pt, red]{};
         \node at (2,2) [circle,fill,inner sep=1.7pt, red]{};
    \end{tikzpicture}
        \caption{$M_1=\mathbb{D}^2(\nicefrac{3}{1},\nicefrac{3}{1})$}
    \end{subfigure}
    \qquad 
    \begin{subfigure}[b]{0.35\textwidth}
        \begin{tikzpicture}[>=latex,scale=1.65]
        \draw[->] (-.5,0) -- (2.5,0) node[right] {$\theta_1$};
        \foreach \x /\n in {1/$\pi$,2/$2\pi$} \draw[shift={(\x,0)}] (0pt,2pt) -- (0pt,-2pt) node[below] {\tiny \n};
        \draw[->] (0,-.5) -- (0,2.5) node[below right] {$\psi_1$};
        \foreach \y /\n in {1/$\pi$,2/$2\pi$}
        \draw[shift={(0,\y)}] (2pt,0pt) -- (-2pt,0pt) node[left] {\tiny \n};
        \node[below left] at (0,0) {\tiny $0$};
         \draw (0,0) rectangle (2,2);
         \draw (1,0) -- (1,2); 
         \draw (0,1) -- (2,1);
         \draw [orange,  ultra thick] (0,0) -- (2,0);
         \draw [orange,  ultra thick] (0,1) -- (2,1);
         \draw [orange,  ultra thick] (0,2) -- (2,2);
         \draw [blue,  ultra thick] (0,2) -- (1,0);
         \draw [blue,  ultra thick] (1,2) -- (2,0);
         \draw [blue,  ultra thick] (0,1) -- (0.5,0);
         \draw [blue,  ultra thick] (0.5,2) -- (1.5,0);
         \draw [blue,  ultra thick] (1.5,2) -- (2,1);
         \node at (0,0) [circle,fill,inner sep=1.7pt, red]{};
         \node at (1,0) [circle,fill,inner sep=1.7pt, red]{};
         \node at (0,1) [circle,fill,inner sep=1.7pt, red]{};
         \node at (1,1) [circle,fill,inner sep=1.7pt, red]{};
         \node at (2,1) [circle,fill,inner sep=1.7pt, red]{};
         \node at (1,2) [circle,fill,inner sep=1.7pt, red]{};
         \node at (2,2) [circle,fill,inner sep=1.7pt, red]{};
         \node at (0,2) [circle,fill,inner sep=1.7pt, red]{};
         \node at (2,0) [circle,fill,inner sep=1.7pt, red]{};
    \end{tikzpicture}
        \caption{$M_1=\mathbb{D}^2(\nicefrac{4}{1},\nicefrac{4}{1})$}
    \end{subfigure} 
        \caption{Four examples of $T(M_1,\partial M_1)\subset \R(\partial M_1)$. The set $H_1$, $A_1$, and $P_1$ are respectively in orange, blue, and red.}
        \label{Some examples of toroidal representation spaces}
\end{figure}

With the notation above, the space $\R(\partial M_1)$ is homeomorphic to the torus $[0,2\pi]^2/_{\sim}$ in the following way:
the point $(\theta_1, \psi_1) \in [0,2\pi]^2/_{\sim}$ is associated to the unique representation
\begin{equation}\label{eq: coordinate for the torus mu1 and h1}
\fund{\partial M_1}\to \Ui, \quad \text{with }\quad \mu_1 \mapsto \begin{bmatrix}
    e^{i\theta_1} & 0 \\ 0 & e^{-i\theta_1}
\end{bmatrix} \quad \text{and} \quad h_1 \mapsto \begin{bmatrix}
    e^{i\phi_1} & 0 \\ 0 & e^{-i\phi_1}
\end{bmatrix}.
\end{equation}
From now on, the sets $H_{M_i}$, $A_{M_i}$, and $P_{M_i}$, for $i \in \{1,2\}$ of Definition \ref{defn: A1 H1 e P1} are denoted by $H_i$, $A_i$, and $P_i$.
Figure \ref{Some examples of toroidal representation spaces} shows the spaces $T(M_1,\partial M_1)$ for $M_1$ homeomorphic to $\mathbb{D}^2(\nicefrac{2}{1},\nicefrac{2}{1})$, $\mathbb{D}^2(\nicefrac{2}{1},\nicefrac{4}{1})$, $\mathbb{D}^2(\nicefrac{3}{1},\nicefrac{3}{1})$, or $\mathbb{D}^2(\nicefrac{4}{1},\nicefrac{4}{1})$. In Figure \ref{Some examples of toroidal representation spaces} the tori $\R(\partial M_1)$ are parameterized with coordinates $(\theta_1,\psi_1)$ as above.
Since $P_1$ consists of central representations, it is contained in $\{0,\pi\}^2 \subset \R(\partial M_1)$. According to Fact \ref{fact: two elements commute iff same centralizer}, if the representation $\rho \in \mathcal{R}(Y)$ is such that $\restr{\rho}{\fund{\partial M_1}} \in H_1$, then $\rho(h_1) \in \ZSU$. This implies that $H_1$ is contained in the ``horizontal" lines $\{\psi_1 \equiv_{\pi} 0 \} \subset \R(\partial M_1)$.

Here is an example of how Theorem \ref{teo: M SU(2)-abeliano se e solo se i pezzi sono empty} can be applied to determinate the $SU(2)$-abelian status of a $3$-manifold.
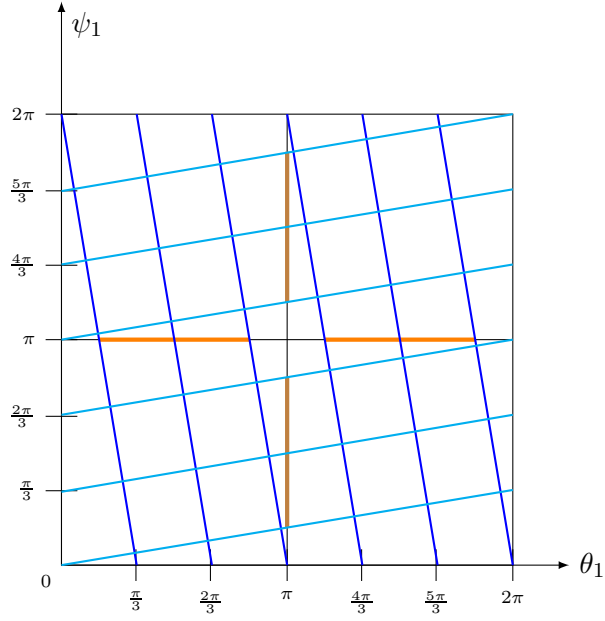
\begin{figure}[ht]
        \begin{tikzpicture}[>=latex,scale=3]
        \draw[->] (0,0) -- (2.25,0) node[right] {$\theta_1$};
        \foreach \x /\n in {1/$\pi$,2/$2\pi$,0.66/$\frac{2\pi}{3}$,0.33/$\frac{\pi}{3}$,1.33/$\frac{4\pi}{3}$,1.66/$\frac{5\pi}{3}$} \draw[shift={(\x,0)}] (0pt,2pt) -- (0pt,-2pt) node[below] {\tiny \n};
        \draw[->] (0,0) -- (0,2.5) node[below right] {$\psi_1$};
        \foreach \y /\n in {1/$\pi$,2/$2\pi$,0.66/$\frac{2\pi}{3}$,0.33/$\frac{\pi}{3}$,1.33/$\frac{4\pi}{3}$,1.66/$\frac{5\pi}{3}$}
        \draw[shift={(0,\y)}] (2pt,0pt) -- (-2pt,0pt) node[left] {\tiny \n};
        \node[below left] at (0,0) {\tiny $0$};
         \draw (0,0) rectangle (2,2);
         \draw (1,0) -- (1,2); 
         \draw (0,1) -- (2,1);
         \draw [orange,  ultra thick] (1/6,1) -- (5/6,1);
         \draw [orange,  ultra thick] (7/6,1) -- (11/6,1);
        \draw [brown,   ultra thick] (1,1/6) -- (1,5/6);
        \draw [brown,   ultra thick] (1,7/6) -- (1,11/6);
         
         \draw [blue,    thick] (0,2) -- (2/6,0);
         \draw [blue,    thick] (2/6,2) -- (4/6,0);
         \draw [blue,    thick] (4/6,2) -- (1,0);
         \draw [blue,    thick] (1,2) -- (8/6,0);
         \draw [blue,    thick] (8/6,2) -- (10/6,0);
         \draw [blue,    thick] (10/6,2) -- (2,0);

         \draw [cyan,    thick] (2,2) -- (0,10/6-1/125);
         \draw [cyan,    thick] (2,10/6) -- (0,8/6);
         \draw [cyan,    thick] (2,8/6) -- (0,1);
         \draw [cyan,    thick] (2,1) -- (0,4/6);
         \draw [cyan,    thick] (2,4/6) -- (0,2/6-1/125);
         \draw [cyan,    thick] (2,2/6) -- (0,0);

    \end{tikzpicture}
        \caption{The intersection $T(M_1,\partial M_1)\cap T(M_2,\partial M_2)$ as in Example \ref{exmp: two trefoil}.}
        \label{picture: T(M1,partial) cap T(M2,partial) of the manifold T(2,3,2,3)}
    \end{figure}
\begin{exmp}
Let $M_1=M_2= S^3 \setminus \nu(T_{2,3})$, where $T_{2,3}$ is the trefoil in $S^3$. It is known that $M_1=M_2=\mathbb{D}^2(\nicefrac{2}{1},\nicefrac{3}{1})$. Let $M$ be the manifold obtained by gluing $M_1$ and $M_2$ along a diffeomorphism $\varphi \colon \partial M_1 \to \partial M_2$ that sends a knot meridian of $M_1$ into a regular fiber of $M_2$ and a regular fiber of $M_1$ into a knot meridian of $M_2$. With respect to our bases,
\[
    \varphi_\ast= \begin{bmatrix}
        0 & 1 \\ 1 & 0
    \end{bmatrix}.
\]
Let $\Sigma$ be the torus in $M$ corresponding to $\partial M_1 = \partial M_2$. Figure \ref{picture: T(M1,partial) cap T(M2,partial) of the manifold T(2,3,2,3)} exhibits the sets $T(M_1,\Sigma) $ and $ T(M_2,\Sigma)$ as a subsets of $\R(\Sigma)$. This latter comes with coordinates $(\theta_1,\psi_1)$ as in \eqref{eq: coordinate for the torus mu1 and h1}. In particular, $H_1$ and $H_2$ are in orange and brown, $A_1$ and $A_2$ are in blue and light blue.
For explicit descriptions of the sets $A_i$ and $H_i$, we refer the reader to the sections \ref{sec: the set A1} and \ref{sec: H1} below.
Let $\rho_1 \colon \fund{M_1}\to SU(2)$ be a representation such that $\rho_1(\fund{\Sigma}) \subset \ZSU$. This implies that $\rho_1(\mu) = \pm 1$, where $\mu \in \fund{M_1}$ is the homotopy class of the knot meridian. Since every element in the knot group $\fund{M_1}= \fund{\compl{T_{2,3}}}$ is conjugate to $\mu$, the image of $\rho_1$ is contained in the center $\ZSU$. This implies that $P_1 = \emptyset$.
Therefore, $P_1 \cap P_2 = \emptyset$. Figure \ref{picture: T(M1,partial) cap T(M2,partial) of the manifold T(2,3,2,3)} implies that 
\[
    T(M_1,\Sigma) \cap T(M_2,\Sigma) = (A_1 \cap A_2) \subset \R(\Sigma).
\]
In particular, we obtain that $H_1 \cap H_2$, $A_1 \cap H_2$, $H_1 \cap A_2$, and $P_1 \cap P_2$ are empty. Theorem \ref{teo: M SU(2)-abeliano se e solo se i pezzi sono empty} implies that $M$ is $SU(2)$-abelian.
\label{exmp: two trefoil}
\end{exmp}
\section{Abelian representations and the set \texorpdfstring{$A_1$}{A\_1}.}
\label{sec: the set A1}
Corollary \ref{cor: an abelian representation extends if and only if the rational longitude is zero} implies that the set $A_1 \subset T(M_1,\partial M_1)$ can be written as
\begin{equation}
    A_1 = \left\{ \eta \in \R(\partial M_1) \,\middle|\, \eta(\lambda_{M_1})^{o_1}=1\right\}.
    \label{eq: description of A1}
\end{equation}

\begin{lemma}
Let $M_1=\mathbb{D}^2 (\nicefrac{p_1}{q_1},\nicefrac{p_2}{q_2})$ and let us suppose that $\fund{M_1}$ is presented as in \eqref{eq: presentation of the fundamental group of M1 and M2}.
Let $\mu_1 \subset \partial M_1$ be the fibration meridian of the chosen presentation and $h_1 \subset \partial M_1$ a regular fiber.
Let $g_1=\gcd(p_1,p_2)$ and $o_1$ be the order of the rational longitude of $M_1$ in $H_1(M_1;\mathbb{Z})$. Then
\[
 \lambda_{M_1} \coloneqq \left(\frac{p_1p_2}{o_1g_1} \right) \mu_1+ \left( \frac{p_1q_2 + p_2q_1}{o_1g_1}\right) h_1 \subset \partial M_1
\]
is the rational longitude of $M_1$. Furthermore, $o_1 = \gcd \left( \frac{p_1p_2}{g_1}, \frac{p_1q_2+p_2q_1}{g_1} \right)$ and $o_1$ divides $g_1$.
\label{lemma: description of the rational longitude}
\begin{proof}

The $\mathbb{Z}$-module $H_1( M_1 ; \mathbb{Z})$ is then generated by the set $\{a_1,b_1,h_1\}$ and related by the equations $p_1a_1 + q_1h_1=0$ and $p_2b_1 + q_2h_1=0$. This presentation implies that $H_1( M_1 ; \mathbb{Z}) = \mathbb{Z} \times \cyclic{g_1}$. Since the element $\lambda_{M_1}$ is a torsion element of $H_1( M_1 ; \mathbb{Z})$, its order $o_1$ divides $g_1$.

We remind the reader that $\mu_1=a_1b_1 \in \fund{\partial M_1}$. Thus, $[\mu_1]=a_1+b_1 \in H_1(M_1;\mathbb{Z})$.
The following is a simple closed curve on $\partial M_1$:
\[
\lambda_{M_1} \coloneqq \left(\frac{p_1p_2}{g_1k} \right) \mu_1+ \left( \frac{p_1q_2 + p_2q_1}{g_1k}\right) h_1, \quad \text{where} \quad k=\gcd\left(\frac{p_1p_2}{g_1},\frac{p_1q_2 + p_2q_1}{g_1}\right).
\]
Since $\lambda_{M_1}$ is a torsion element of $H_1(M_1;\mathbb{Z})$, the curve $\lambda_{M_1}\subset \partial M_1$ as above is the rational longitude of $M_1$.

We prove now the expression for the order $o_1$. Since $o_1\cdot \lambda_{M_1}$ is trivial in $H_1( M_1 ; \mathbb{Z})$, there exist two integers $n$ and $m$ such that the following identities hold in $H_1( M_1 ; \mathbb{Z})$:
\begin{equation}
    o_1\left(\frac{p_1p_2}{g_1k}\right) (a_1+b_1) + o_1\left(\frac{p_1q_2+p_2q_1}{g_1k}\right) h_1 = o_1 \cdot \lambda_{M_1}= n(p_1a_1+q_1h_1) + m(p_2 b_1 + q_2 h_1) .
    \label{eq: la uso una volta, order of the }
\end{equation}
Thus, we obtain that
$o_1\left(\frac{p_1p_2}{g_1k}\right) = np_1 = mp_2$. Therefore, there exists a $j \in \mathbb{Z}_{\neq 0}$ such that
$
    j \,\text{lcm}(p_1,p_2)= o_1\left(\frac{p_1p_2}{g_1k}\right).
$
Since $\text{lcm}(p_1,p_2)=\frac{p_1p_2}{g_1}$, we obtain that
$
j \left(\frac{p_1p_2}{g_1}\right)= o_1 \left(\frac{p_1p_2}{g_1k}\right)=np_1=mp_2.
$
This latter implies that $n= \frac{jp_2}{g_1}$ and $m= \frac{jp_1}{g_1}$. Hence, the \eqref{eq: la uso una volta, order of the } becomes
\[
o_1 \cdot \lambda_{M_1} = \frac{jp_2}{g_1}(p_1a_1+q_1h_1)+  \frac{jp_1}{g_1}(p_2a_2+q_2h_1)= j \left(\frac{p_2}{g_1}(p_1a_1+q_1h_1)+ \frac{p_1}{g_1}(p_2a_2+q_2h_1) \right).
\]
The presentation of the first homology implies that the quantity
$
    \frac{p_2}{g_1}(p_1a_1+q_1h_1)+ \frac{p_1}{g_1}(p_2a_2+q_2h_1)
$
is trivial in $H_1(M_1;\mathbb{Z})$. The minimality of $o_1$ implies $j=\pm 1$. Hence,
$
    \frac{p_1p_2}{g_1}= o_1\left(\frac{p_1p_2}{g_1k}\right).
$
This implies that $o_1=k=\gcd\left(\frac{p_1p_2}{g_1},\frac{p_1q_2 + p_2q_1}{g_1}\right)$.
\end{proof}
\end{lemma}
\section{Central representations and the set \texorpdfstring{$P_1$}{P\_1}}
The diffeomorphism in \eqref{eq: description of M1 in terms of fraction and disk. It says that q1 and q2 can be taken} implies that $q_1$ and $q_2$ can be chosen to be odd. In this section we consider $\fund{M_1}$ to be presented as in \eqref{eq: presentation of the fundamental group of M1 and M2} with $q_1$ and $q_2$ odd.

We consider $\fund{\partial M_1}$ generated by the basis $\{\mu_1,h_1\}$. We remind that $\mu_1=a_1b_1$ when it is considered as an element of $\fund{M_1}$.
If the representation $\rho \in \mathcal{R}(M_1)$ is such that the restriction $\restr{\rho}{\fund{\partial M_1}}$ is central, then $\rho(a_1b_1)=\pm 1$. The latter implies that
$\rho(a_1)$ commutes with $\rho(b_1)$, and hence that $\rho$ is an abelian representation.
We remind the reader that $g_1\coloneqq\gcd(p_1,p_2)$.

\begin{lemma}
    Let $\eta \colon \fund{\partial M_1} \to SU(2)$ be a representation such that $\eta(h_1)=-1$ and let $g_1 \equiv_2 0$.
    If $\eta$ extends to a representation $\fund{M_1} \to SU(2)$, then every such extension is non-central.
    \begin{proof}
        Suppose that the representation $\eta$ extends to a central representation $\rho \colon \fund{M_1} \to SU(2)$. This implies that $\rho(a_1)$ and $\rho(b_1)$ are both in $\ZSU$ and $\rho(h_1)=\eta(h_1)=-1$. Since $p_1$ is even and $q_1$ is odd, we obtain that
        \[
        1=\rho(a_1)^{p_1}\rho(h_1)^{q_1}=(\pm 1)^{p_1}(-1)^{q_1}=-1.
        \]
        This is a contradiction.
    \end{proof}
    \label{lemma: g1 pari e h1=-1 allora non-central}
\end{lemma}

\begin{lemma}
    If $g_1 \ge 3$, then
    the trivial representation $\eta \colon \fund{\partial M_1} \to 1$ extends to a non-central representation $\fund{M_1} \to SU(2)$.
    \begin{proof}
        Let $\mathcal{A}$ be the abelianization homomorphism and $G$ the abelian group
        $
            G= \mathcal{A}\left( \frac{\fund{M_1}}{\normalsubgroup{a_1b_1,h_1}}\right).
        $
        The presentation \eqref{eq: presentation of the fundamental group of M1 and M2} of $\fund{M_1}$ implies that $G$ is isomorphic to $\cyclic{g_1}$.
        Let $q: \fund{M_1} \twoheadrightarrow G$ be the quotient map.
        Since $g_1 \ge 3$, the group $G$ admits a non-central $SU(2)$-representation. Let $\gamma \colon G \to SU(2)$ be a non-central representation, then
        $\gamma \circ q \colon \fund{M_1}\to SU(2)$ is a non-central $SU(2)$-representation of $\fund{M_1}$. Furthermore,
        $\restr{(\gamma \circ q)}{\fund{\partial M_1}}$ is the trivial representation $\eta$. This implies that the representation $\eta$ admits a non-central extension.
    \end{proof}
    \label{lemma: g1>=3 the trivial representation is not central}
\end{lemma}

\begin{lemma}\label{lemma: g1 >=3, allora non-central}
    Let $\eta \colon \fund{\partial M_1} \to \ZSU$ be a representation and let $g_1 \ge 3$.
    If $\eta$ extends to a representation $\fund{M_1} \to SU(2)$, then it admits a non-central extension.
    \begin{proof}
        We split the proof in the cases $g_1 \equiv_20$ and $g_1 \equiv_21$.

        Let us suppose that $g_1 \equiv_2 0$. This implies that $g_1 \ge 4$. If $\eta(h_1)=-1$, then Lemma \ref{lemma: g1 pari e h1=-1 allora non-central} gives the conclusion. If $\eta$ is the trivial representation, the conclusion holds by Lemma \ref{lemma: g1>=3 the trivial representation is not central}. Therefore, we prove the remaining case: we suppose that $\eta(a_1b_1)=-1$ and $\eta(h_1)=1$. Since $p_1 \equiv_2 p_2 \equiv_2 0$, we define a representation $\rho \colon \fund{M_1} \to SU(2)$ as
        \[
            \rho(a_1)= 
            \begin{bmatrix}
                e^{\frac{2\pi i}{g_1}} & 0 \\ 0 & e^{-\frac{2\pi i}{g_1}}
            \end{bmatrix}
            , 
            \quad
            \rho(b_1)=-
            \begin{bmatrix}
            e^{-\frac{2\pi i}{g_1}} & 0 \\ 0 & e^{\frac{2\pi i}{g_1}}
            \end{bmatrix}
            , \quad \text{and} \quad \rho(h_1)=1.
        \]
        Since $g_1 \ge 4$, such a representation $\rho$ is non-central and it restricts to $\eta$.

        Let us suppose that $g_1 \equiv_2 1$ and $g_1 \ge 3$.
        Since $g_1$, $q_1$, and $q_2$ are odd, we obtain that
        $
            \frac{p_1p_2}{g_1} \equiv_2 p_1p_2$ and $\frac{p_1 q_2+ p_2q_1}{g_1}\equiv_2 p_1 + p_2. 
        $
        Let $\eta \colon \fund{\partial M_1} \to \ZSU$ be a representation that extends to $\fund{M_1}$. As we said before, every such an extension is abelian. By Corollary \ref{cor: an abelian representation extends if and only if the rational longitude is zero}, we obtain that $\eta(\lambda_{M_1})^{o_1}=1$.
        According to Lemma \ref{lemma: description of the rational longitude}, we obtain that
        \begin{equation}
        1=\eta(\lambda_{M_1})^{o_1}= \eta(a_1b_1)^{\frac{p_1p_2}{g_1}}\eta(h_1)^{\frac{p_1q_2 +p_2q_1}{g_1}}=\eta(a_1b_1)^{p_1p_2}\eta(h_1)^{p_1 +p_2}.
        \label{eq: la uso due volte g1>=3}
        \end{equation}
        If $p_1p_2 \equiv_2 1$, then $p_1$ and $p_2$ are both odd and $p_1+p_2 \equiv_2 0$. Equation \eqref{eq: la uso due volte g1>=3} implies that $\eta(a_1b_1)=1$. If $\eta(h_1)=1$, then the conclusion holds by Lemma \ref{lemma: g1>=3 the trivial representation is not central}. Without loss of generality, we can suppose that $\eta(h_1)=-1$. The representation $\rho \colon \fund{M_1}\to SU(2)$ with
        \[
            \rho(a_1)=
            \begin{bmatrix}
                e^{i\frac{\pi}{g_1}} & 0 \\ 0 & e^{-i\frac{\pi}{g_1}}
            \end{bmatrix}
            , \quad \rho(b_1)=
            \begin{bmatrix}
            e^{-i\frac{\pi}{g_1}} & 0 \\ 0 & e^{i\frac{\pi}{g_1}}
            \end{bmatrix}
            \quad
            \text{and} \quad \rho(h_1)=-1,
        \]
        restricts to $\eta$. Since $g_1 \ge 3$, the representation $\rho$ is non-central.
        
        If $p_1p_2 \equiv_2 0$, then $p_1 + p_2 \equiv_2 1$. Equation \eqref{eq: la uso due volte g1>=3} implies that $\eta(h_1)=1$. Again, if $\eta(a_1b_1)=1$, then the conclusion is implied by Lemma \ref{lemma: g1>=3 the trivial representation is not central}. Without loss of generality, we can suppose $\eta(a_1b_1)=-1$.
        Since $o_1$ divides $g_1$ by Lemma \ref{lemma: description of the rational longitude}, the quantity $o_1$ is odd.
        Since $g_1 \equiv_2 1$ and $p_1p_2 \equiv_20$, we can also assume that $p_1 \equiv_20$ and $p_2 \equiv_2 1$. Let $g_1=2n+1$ with $n \in \mathbb{N}$. The representation $\eta$ extends to the representation $\rho \colon \fund{M_1} \to SU(2)$ with
        \[
            \rho(a_1)=\begin{bmatrix}
                e^{i\frac{\pi}{g_1}} & 0 \\ 0 & e^{-i\frac{\pi}{g_1}}
            \end{bmatrix}, \quad \rho(b_1)=\begin{bmatrix}
            e^{i\frac{2 \pi n}{g_1}} & 0 \\ 0 & e^{-i\frac{2 \pi n}{g_1}}
            \end{bmatrix}, \quad \text{and} \quad \rho(h_1)=1.
        \]
        If $g_1 \ge 3$, then $\rho$ has a non-central image.
    \end{proof}
\end{lemma}

\begin{lemma}
    Let $\eta \colon \fund{\partial M_1} \to \ZSU$ be a representation, then $\eta$ extends to a representation $\rho \colon \fund{M_1}\to SU(2)$ if and only if $\eta(\lambda_{M_1})^{o_1}=1$. If this happens, the following hold:
    \begin{itemize}
        \item If $g_1=1$, then $\eta$ extends only to central representations;
        \item If $g_1=2$, then $\eta$ extends to a non-central representation if and only if $o_1=2$ and $\eta(h_1)=-1$;
        \item If $g_1 \ge 3$, then $\eta$ extends to a non-central representation.
    \end{itemize}
    \label{lemma: description of the repersentation P1}
    \begin{proof}
        As we stated before, if $\eta \colon \fund{\partial M_1}\to \ZSU$ extends to the representation $\rho \colon \fund{M_1}\to SU(2)$, then $\rho$ has abelian image. Thus, the first part of the statement is a consequence of Corollary \ref{cor: an abelian representation extends if and only if the rational longitude is zero}. Moreover, if $g_1 \ge 3$, then the conclusion holds by Lemma \ref{lemma: g1 >=3, allora non-central}.

        Let $\varepsilon_1,\varepsilon_2 \in \{\pm 1\} = \ZSU$ and let $\tau_1,\tau_2 \in \{1,2\}$ be their orders in $SU(2)$. 
        Let $\rho \colon \fund{M_1} \to SU(2)$ be such that $\rho(h_1)=\varepsilon_1$ and $\rho(a_1b_1)=\varepsilon_2$. Let $\eta$ be the restriction of $\rho$ to $\fund{\partial M_1}$. We define the group $F_{\varepsilon_1,\varepsilon_2}$ as the the abelianization of 
        $
            \frac{\fund{M_1}}{ \normalsubgroup{h_1^{\tau_1},(a_1b_1)^{\tau_2}}},
        $
        for the corresponding $\tau_1$ and $\tau_2$. Consequently, the following diagram commutes:
        \begin{equation}
            \begin{tikzcd}
            && SU(2) && \\
            \fund{\partial M_1} \arrow[rr, "\iota_\ast"] \arrow[rru, bend left=10, "\eta" ] & & 
            \fund{M_1} \arrow[rr, twoheadrightarrow, "\mathcal{F}"] \arrow[u, "\rho"] & &
            \mathcal{A}\left(\frac{\fund{M_1}}{ \normalsubgroup{h_1^{\tau_1},(a_1b_1)^{\tau_2}}}\right) \eqqcolon F_{\varepsilon_1,\varepsilon_2}.
            \arrow[llu, bend right=10,"\bar{\rho}"']        
            \end{tikzcd}
            \label{eq: factorization group diagram}
        \end{equation}
        Here $\mathcal{F}$ is the quotient map and $\iota_\ast$ is the map induced by the inclusion $\iota: \partial M_1 \to M_1$.

        Let us suppose that $g_1=1$. If $\eta$ extends to a representation $\rho$. Diagram \ref{eq: factorization group diagram} implies that $\rho$ factors through the group 
        \[
            F_{\varepsilon_1,\varepsilon_2}= \coker 
            \begin{bmatrix}
                p_1 & 0  & 0 & \tau_2 \\
                0 & p_2 & 0 & \tau_2 \\
                q_1 & q_2 & \tau_1 & 0
            \end{bmatrix}.
        \]
        Then Smith's algorithm implies that $F_{\varepsilon_1,\varepsilon_2} \in \{\{1\}, \cyclic{2}\}$. Thus, the representation $\rho$ factors through the center $\ZSU$. This implies the conclusion for the case $g_1=1$.
        
        Let us suppose that $g_1=2$. We are going to prove that $\eta$ extends to a non-central representation if and only if $\eta(h_1)=-1$ and $o_1=2$.
        Let us suppose that $o_1=2$.
        Let $\eta\colon \fund{\partial M_1}\to \ZSU$ be a representation such that $\eta(h_1)=-1$. Since $\eta$ is a central representation, we obtain that $\eta(\lambda_{M_1})^{o_1}=1$.
        According to Corollary \ref{cor: an abelian representation extends if and only if the rational longitude is zero}, the representation $\eta$ extends to $\fund{M_1}$. Lemma \ref{lemma: g1 pari e h1=-1 allora non-central} implies that $\eta$ extends to a non-central representation. This concludes one direction of the case $g_1=2$.
        
        Conversely, suppose that $\eta\colon \fund{\partial M_1}\to \ZSU $ extends to a non-central representation $\fund{M_1}\to SU(2)$. We need to show that this implies that $o_1=2$ and $\eta(h_1)=-1$.
        For this purpose, we first prove that if the representation $\rho \colon \fund{\partial M_1}\to \ZSU$ extends to $\fund{M_1}$ and $\rho(h_1)=1$, then every extension of $\rho$ is central.
        
        Let $\eta \colon \fund{\partial M_1}\to \ZSU$ be a representation that extends to $\fund{M_1}$ such that $\eta(h_1)=1$. Let $\rho \colon \fund{M_1} \to SU(2)$ be an extension of $\eta$. According to diagram \ref{eq: factorization group diagram}, the representation $\rho$ factors through the group
        \[
            F_{1,\varepsilon_2}= \coker 
            \begin{bmatrix}
                p_1 & 0  & 0 & \tau_2 \\
                0 & p_2 & 0 & \tau_2 \\
                q_1 & q_2 & 1 & 0
            \end{bmatrix}= \coker 
            \begin{bmatrix}
                p_1 & 0  & 0 & \tau_2 \\
                0 & p_2 & 0 & \tau_2 \\
                0 & 0 & 1 & 0
            \end{bmatrix}.
        \]
        The Smith's algorithm implies that $F_{1,\varepsilon_2} \in \{ \cyclic{2},\cyclic{2}\times \cyclic{2}\}$. Hence, $\rho$ is central. This implies that if the representation $\eta$ extends to a non-central representation, then $\eta(h_1)=-1$.

        Let $\eta \colon \fund{\partial M_1}\to \ZSU$ be a representation that extends to an $SU(2)$-representation of $\fund{M_1}$ whose image is non-central. We proved before that $\eta(h_1)=-1$. We are going to prove that $o_1=2$.
        Since $g_1=2$, the integer $\frac{p_1p_2}{2}$ is even.
        Since $\eta$ extends to $\fund{M_1}$ by hypothesis, Corollary \ref{cor: an abelian representation extends if and only if the rational longitude is zero} and Lemma \ref{lemma: description of the rational longitude} imply that
        \begin{align}
        1=\eta(\lambda_{M_1})^{o_1} =\left( \cancel{\eta(a_1b_1)^{\frac{p_1p_2}{2}}}\eta(h_1)^{\frac{p_1q_2+p_2q_1}{2}}\right)^{o_1} = (-1)^{o_1 \left(\frac{p_1q_2+p_2q_1}{2}\right)}.
        \label{eq: uso una volta per contraddizione}
        \end{align}
        In particular, equation \eqref{eq: uso una volta per contraddizione} implies that $o_1\left(\frac{p_1q_2+p_2q_1}{2}\right)$ is even.
        Lemma \ref{lemma: description of the rational longitude} implies that $o_1$ is even if and only if the integer $\frac{p_1q_2+p_2q_1}{2}$ is even.
        Hence, if $o_1$ is odd, then $o_1 \left(\frac{p_1q_2+p_2q_1}{2}\right)$ is also odd, which contradicts \eqref{eq: uso una volta per contraddizione}.
        This implies that $o_1$ has to be even. Since, by Lemma \ref{lemma: description of the rational longitude}, $o_1$ divides $g_1=2$, we obtain that $o_1=2$.
        \end{proof}
\end{lemma}

\section{The irreducible representations and the set \texorpdfstring{$H_1$}{H\_1}}
\label{sec: H1}
Let $\rho \colon \fund{M_1}\to SU(2)$ be such that $\restr{\rho}{\fund{\partial M_1}} \coloneq \eta \in H_1$,
Fact \ref{fact: two elements commute iff same centralizer} implies that $\rho(h_1)=\eta(h_1) \in \ZSU$. Hence, we divide $H_1$ into the sets $H_{1,0}$ and $ H_{1,\pi}$ where
\[
    H_{1,0}= \left\{\eta \in H_1 \, \middle| \, \eta(h_1)=1\right\}
    \quad \text{and} \quad
    H_{1,\pi}=\left\{\eta \in H_1 \, \middle| \, \eta(h_1)=-1\right\}.
\]
Clearly $H_1 = H_{1,0} \cup H_{1,\pi}$.
We use the coordinates $(\theta_1,\psi_1)$ for the space $\R(\partial M_1)$ with respect to the basis $\{\mu_1,h_1\}$ as in \eqref{eq: coordinate for the torus mu1 and h1}.
With this parameterization of $\R(\partial M_1)$, we have that $H_{1,0} \subset \{\psi_1 = 0 \} \subset \R (\partial M_1)$ and $H_{1,\pi} \subset \{\psi_1 = \pi \} \subset \R (\partial M_1)$.

\begin{lemma}\label{lemma: H1 doesn't contain -2 and 2}
    Let $\fund{M_1}$ be presented as in \eqref{eq: presentation of the fundamental group of M1 and M2}. If the representation $\eta \colon \fund{\partial M_1} \to SU(2)$ is such that $\Tr \eta(a_1b_1) = \pm 2$, then $\eta$ is not in $H_1$.
    \begin{proof}
        Let us suppose that $\rho: \fund{M_1} \to SU(2)$ is an extension of $\eta$. Since $\Tr\eta(a_1b_1)=\Tr \rho(a_1b_1)=\pm 2$, then $\rho(a_1b_1)=\pm 1$. Thus, $[\rho(a_1),\rho(b_1)]=1$. This implies that $\rho$ is an abelian representation and hence $\restr{\rho}{\fund{\partial M_1}} =\eta \notin H_1$.
    \end{proof}
\end{lemma}
Lemma \ref{lemma: H1 doesn't contain -2 and 2} implies that the set $H_1$ has no intersection with the lines
\[
    \left\{\theta_1 \equiv_{2\pi} 0\right\} \cup \left\{\theta_1 \equiv_{2\pi} \pi\right\} \subset \R(\partial M_1).
\]

\begin{prop}
    Given $a,b,c \in (-2,2)$, there exist two matrices $A,B\in SU(2)$ with $\Tr A=a$, $\Tr B=b$, $\Tr AB=c$, and $AB \neq BA$ if and only if
    \[
        \frac{1}{2}\left(a b-\sqrt{(4-a^2)(4-b^2)}\right) < c < \frac{1}{2}\left(a b+\sqrt{(4-a^2)(4-b^2)}\right).
    \]
    \label{prop: esistenza di AB basato sulle traccie}
    \begin{proof}
        In this proof we think of $S^1$ as a subset of $\mathbb{C}$ in the usual way. Let $z$ be a complex number, we denote by $\Re(z)$ and $\Im(z)$ its real and imaginary part.
        Let us suppose that there exist two matrices $A,B\in SU(2)$ such that $\Tr A=a$, $\Tr B=b$, $\Tr AB=c$, and $AB \neq BA$. Up to conjugation, we can suppose that
        \[
            A= \begin{bmatrix}
                u & 0 \\ 0 & \overline{u}
            \end{bmatrix} \quad \text{and} \quad 
            B= \begin{bmatrix}
                \alpha & -\overline{\beta} \\ \beta & \overline{\alpha}
            \end{bmatrix},
        \]
        where $u \in S^1$ and $\alpha,\beta \in \mathbb{C}$. Since $A$ and $B$ are assumed not to commute, $B$ is not diagonal. This implies that $\beta \neq 0$ and $|\alpha|^2 < 1$. Thus, there exist $v \in S^1$ and $t \in (0,1)$ such that
        $
            \alpha= (1-t) v + t\overline{v}.
        $
        By hypothesis, the complex numbers $\alpha$, $u$, and $v$ are such that
        \[
         2\Re( u) = \Tr A=a, \quad 2\Re(\alpha)= 2\Re(v)= \Tr B= b, \quad \text{and} \quad  \Tr{AB}=c.
        \]
        Let us consider the following:
        \begin{equation}
            \begin{aligned}
            c= &\Tr{AB} = u\alpha+ \overline{u}\overline{\alpha} = 2 \Re (u \alpha) = 2 \left((1-t)\Re(uv)+t\Re(u\overline{v})\right) \\
            &= 2 \Re(u)\Re(v) +2(2t-1)\Im(u)\Im(v)= \frac{ab}{2}+2(2t-1)\sqrt{1-\frac{a^2}{4}}\sqrt{1-\frac{b^2}{4}}  \\
            & = \frac{1}{2}\left( ab+(2t-1)\sqrt{(4-a^2)(4-b^2)}\right).
            \end{aligned}
            \label{eq: la uso una volta Tr(AB)}
        \end{equation}
        Since $t$ is neither $0$ nor $1$, we conclude that
        \[
        \frac{1}{2}\left(a b-\sqrt{(4-a^2)(4-b^2)}\right) < c < \frac{1}{2}\left(a b+\sqrt{(4-a^2)(4-b^2)}\right).
        \]
        
        Conversely, let $u,v \in S^1$ with
        \[
        u = \frac{a}{2} + i \sqrt{1-\frac{a^2}{4}} \quad \text{and} \quad v = \frac{b}{2} + i \sqrt{1-\frac{b^2}{4}}.
        \]
        For $t \in (0,1)$ we set $\alpha(t)= (1-t)v+ t\overline{v}$ and define $A,B(t) \in SU(2)$ by
        \[
        A= \begin{bmatrix}
            u & 0 \\ 0 & \overline{u}
        \end{bmatrix}
        \quad \text{and} \quad B(t)= 
        \begin{bmatrix}
            \alpha(t) & \sqrt{1-|\alpha(t)|^2} \\ -\sqrt{1-|\alpha(t)|^2} & \overline{\alpha}(t)
        \end{bmatrix}.
        \]
        It is easy to see that $\Tr{A}=a$ and $\Tr{B(t)}=b$. Since $t$ is neither $0$ nor $1$, the matrix $B$ is not diagonal and hence $AB(t) \neq B(t)A$ for every $t \in (0,1)$. The computation in \eqref{eq: la uso una volta Tr(AB)} implies that the trace of the multiplication $AB(t)$ is
        \begin{equation*}
            \begin{aligned}
            \Tr{AB(t)}= u\alpha(t)+ \overline{u}\overline{\alpha}(t)= \frac{1}{2}\left( ab+(2t-1)\sqrt{(4-a^2)(4-b^2)}\right).
            \end{aligned}
        \end{equation*}
        Therefore, there is some $t \in (0,1)$ for which $\Tr{AB(t)}=c$.
    \end{proof}
\end{prop}

\begin{defn}\label{defn: interval I(a,b)}
For $a,b\in [-2,2]\subset \mathbb{R}$ we define $I(a,b)$ as the open interval
\[
    \left( \frac{1}{2}\left(a b-\sqrt{(4-a^2)(4-b^2)}\right), \frac{1}{2}\left(ab+\sqrt{(4-a^2)(4-b^2)}\right)\right) \subseteq [-2,2].
\]
\end{defn}
\begin{defn}\label{defn: interval J0(p1,p2)}
    Let $2 \le p_1  \le p_2$ be two natural numbers, we define the $J_0(p_1,p_2)$ as
    \[
        J_0(p_1,p_2) = \bigcup_{k_1,k_2 \in \mathbb{Z}} I \left(2 \cos\left( \frac{2\pi k_1}{p_1}\right), 2 \cos\left( \frac{2\pi k_2}{p_2}\right) \right) \subseteq [-2,2].
    \]
\end{defn}
\begin{defn} \label{defn: interval Jpi(p1,p2)}
    Let $2 \le p_1  \le p_2$ be two natural numbers, we define the $J_\pi(p_1,p_2)$ as
    \[
        J_\pi(p_1,p_2) = \bigcup_{\substack{k_1,k_2 \in \mathbb{Z} \\ k_i \text{ odd} }} I \left(2 \cos\left( \frac{\pi k_1}{p_1}\right), 2 \cos\left( \frac{\pi k_2}{p_2}\right) \right) \subseteq [-2,2].
    \]
\end{defn}

In order to make the notation lighter, the interval $I\left(2\cos \left(\frac{2\pi k_1}{p_1}\right),2\cos\left(\frac{2\pi k_2}{p_2}\right)\right)$ as in Definition \ref{defn: interval J0(p1,p2)} will be denoted as $I\left(\nicefrac{k_1}{p_1},\nicefrac{k_2}{p_2}\right)$.
Thus, $J_\pi(p_1,p_2)$ is the union of $I\left(\nicefrac{k_1}{2p_1},\nicefrac{k_2}{2p_2}\right)$ with $k_1$ and $k_2$ odd.
The length of $I\left(\nicefrac{k_1}{p_1},\nicefrac{k_2}{p_2}\right)$ can be computed from Definition \ref{defn: interval I(a,b)} and it equals
\begin{align}  
    m \left( I \left(\frac{ k_1}{p_1},\frac{ k_2}{p_2}\right) \right) = 4 \left| \sin \left(\frac{2\pi k_1}{p_1}\right) \sin\left(\frac{2\pi k_2}{p_2}\right) \right|.
    \label{eq: length of the interval J(a,b)}
\end{align}
Similarly, the length of $I\left(\nicefrac{k_1}{2p_1},\nicefrac{k_2}{2p_2}\right)$ is computed from \eqref{eq: length of the interval J(a,b)}.  
\begin{lemma}
    Let $2\le p_1 \le p_2$. The set $J_0(p_1,p_2)$ is empty if and only if $p_1=2$. Furthermore, the set $J_\pi(p_1,p_2)$ is not empty.
    \begin{proof}
       The identity \eqref{eq: length of the interval J(a,b)} implies that $I\left(\nicefrac{k_1}{p_1},\nicefrac{k_2}{p_2}\right)$ has positive measure if and only if $\sin \left(\nicefrac{2\pi k_1}{p_1}\right)$ and $\sin \left(\nicefrac{2\pi k_2}{p_2}\right)$ are both non-zero. This implies that $J_0(2,p_2)$ is empty. Conversely, if $3 \le p_1 \le p_2$, then there exists $(k_1,k_2) \in \mathbb{Z}^2$ so that $\sin \left(\nicefrac{2\pi k_1}{p_1}\right)\neq0$ and $\sin \left(\nicefrac{2\pi k_2}{p_2}\right)\neq0$. This implies that $I\left(\nicefrac{k_1}{p_1},\nicefrac{k_2}{p_2}\right)$ is nonempty, and therefore that $J_0(p_1,p_2)$ is nonempty as well.
       
       The set $J_\pi(p_1,p_2)$ contains the interval $I\left(\nicefrac{1}{2p_1},\nicefrac{1}{2p_2}\right)$ that is not empty by the \eqref{eq: length of the interval J(a,b)}. Hence, the set $J_\pi(p_1,p_2)$ is not empty.
    \end{proof}
    \label{lemma: J0(2,p2) è vuoto}
\end{lemma}

\begin{rmk}
    Let $n \ge 1$ be a natural number. Let $A \in SU(2)$. The matrix $A$ is such that $A^n=1$ if and only if $\Tr A= 2\cos\left(\frac{2\pi k }{n}\right)$ for a $k \in \mathbb{Z}$. Moreover, $A^n=-1$ if and only if $\Tr A= 2\cos\left(\frac{\pi k}{n}\right)$ for an odd $k \in \mathbb{Z}$.
    \label{rmk: A n =-+1 se e solo se traccia}
\end{rmk}

The next lemma and Proposition \ref{prop: esistenza di AB basato sulle traccie} are, the two key results for constructing the set $H_1 \subset T(M_1,\partial M_1)$ and thus to understand $H_1 \cap H_2$, as we shall do in Section \ref{sec: H1 cap H2}.

\begin{lemma}
   Let $\fund{M_1}$ be presented as in \eqref{eq: presentation of the fundamental group of M1 and M2} with $q_1$ and $q_2$ odd. Then the maps
   \begin{alignat*}{2}
       f_0 \colon  H_{1,0}  &\to J_0(p_1,p_2) \quad \quad \text{and} \quad \quad f_{\pi} \colon H&_{1,\pi} \to  J_\pi(p_1,p_2) \\
       \eta &\mapsto \Tr\eta(a_1b_1) &\eta \mapsto \Tr\eta(a_1b_1) 
   \end{alignat*}
   are well-defined and surjective.
\begin{proof}
    Let us suppose that $p_1=2$. The set $J_0(2,p_2)$ is empty by Lemma \ref{lemma: J0(2,p2) è vuoto}. Since $p_1=2$, every representation $\rho \in \mathcal{R}(M_1)$ with $\rho(h_1)=1$ is such that $\rho(a_1)\in \ZSU$. Hence, the representation $\rho$ has abelian image. This implies that $H_{1,0}$ is empty.
    
    Let us suppose that $3 \le p_1 \le p_2$. If $\eta \in H_{1,0}$, then there exists an irreducible representation $\rho \in \mathcal{R}(M_1)$ such that $\restr{\rho}{\fund{\partial M_1}} \equiv \eta$ and $\rho(h_1)=1$. In particular, we have that
    \[
        \rho(a_1)^{p_1}\rho(h_1)^{q_1}=\rho(a_1)^{p_1}=1  \quad \text{and} \quad \rho(b_1)^{p_2}\rho(h_1)^{q_2}=\rho(b_1)^{p_2}=1.
    \]
    Since $\rho$ is irreducible, we obtain that $\rho(a_1)\rho(b_1)\neq \rho(b_1)\rho(a_1)$.
    According to Remark \ref{rmk: A n =-+1 se e solo se traccia}, there exist two integers $k_1$ and $k_2$ such that $\Tr \rho(a_1)=2\cos\left(\frac{2\pi k_1}{p_1} \right)$ and $\Tr \rho(b_1)=2\cos\left(\frac{2\pi k_2}{p_2}\right)$. Proposition \ref{prop: esistenza di AB basato sulle traccie} implies that
    \[
    \Tr \eta(a_1b_1)=\Tr \rho(a_1b_1) \in I \left( \frac{k_1}{p_1},\frac{k_2}{p_2}\right) \subseteq J_0(p_1,p_2).
    \]
    This implies that the map $f_0$ is well defined. 
    
    Let $z \in J_0(p_1,p_2)$. Thus, there exist two integers $k_1$ and $k_2$ such that $z \in I \left( \nicefrac{k_1}{p_1},\nicefrac{k_2}{p_2}\right) \subseteq J_0(p_1,p_2)$. According to Proposition \ref{prop: esistenza di AB basato sulle traccie}, there exist two matrices $A$ and $B$ of $SU(2)$ such that
    \[
        \Tr A = 2\cos \left(\frac{2\pi k_1}{p_1}\right), \quad \Tr B = 2\cos \left(\frac{2\pi k_2}{p_2}\right), \quad \Tr AB =z, \quad \text{and} \quad AB \neq BA.
    \]
    Remark \ref{rmk: A n =-+1 se e solo se traccia} implies that $A^{p_1}=1$ and $B^{p_2}=1$. Let $\rho \in \mathcal{R}(M_1)$ be the representation defined as
    \[
    \rho(a_1)=A, \quad \rho(b_1)=B, \quad \text{and} \quad \rho(h_1)=1.
    \]
    Since $A$ and $B$ do not commute, $\rho$ is irreducible. Up to conjugation, we can suppose that $\restr{\rho}{\fund{\partial M_1}} \in \R(\partial M_1)$. In particular this implies that $\restr{\rho}{\fund{\partial M_1}} \in H_{1,0}$. Thus, $f_0(\restr{\rho}{\fund{\partial M_1}})= \Tr \rho(a_1b_1) = \Tr AB=z$. This implies that the map $f_0$ is surjective.

    Let $\eta \in H_{1,\pi}$. This implies that there exists an irreducible representation $\rho \in \mathcal{R}(M_1)$ such that $\restr{\rho}{\fund{\partial M_1}}\equiv \eta$ and $\rho(h_1)=-1$.
    Since $q_1$ and $q_2$ are odd, we obtain that
    \[
        \rho(a_1)^{p_1}=-1 \quad \text{and} \quad \rho(b_1)^{p_2}=-1.
    \]
    As a result of Remark \ref{rmk: A n =-+1 se e solo se traccia}, there exist two odd integers $k_1$ and $k_2$ such that $\Tr \rho(a_1)=\frac{\pi k_1}{p_1}$ and $\Tr \rho(b_1)=\frac{\pi k_2}{p_2}$. Proposition \ref{prop: esistenza di AB basato sulle traccie} implies that
    \[
    \Tr \eta(a_1b_1)=\Tr \rho(a_1b_1) \in I \left( \frac{ k_1}{2p_1},\frac{ k_2}{2p_2}\right) \subseteq J_\pi(p_1,p_2).
    \]
    This implies that $f_\pi$ is well defined. 
    
    Let $z \in J_\pi(p_1,p_2)$. As a result of Remark \ref{rmk: A n =-+1 se e solo se traccia}, there exist two odd integers $k_1$ and $k_2$ such that $z \in I \left( \nicefrac{ k_1}{2p_1},\nicefrac{ k_2}{2p_2}\right) \subseteq J_\pi(p_1,p_2)$. According to Proposition \ref{prop: esistenza di AB basato sulle traccie}, there exist two matrices $A$ and $B$ of $SU(2)$ such that
    \[
        \Tr A = 2\cos \left( \frac{\pi k_1}{p_1}\right), \quad \Tr B = 2\cos \left(\frac{\pi k_2}{p_2} \right), \quad \Tr AB =z, \quad \text{and} \quad AB \neq BA.
    \]
    Let $\rho \colon \fund{M_1}\to SU(2)$ be the representation determined by the following:
    \[
        \rho(a_1)=A, \quad \rho(b_1)=B, \quad \text{and} \quad \rho(h_1)=-1.
    \]
    Since $A$ and $B$ do not commute, $\rho$ is irreducible. Up to conjugation, we can suppose that $\restr{\rho}{\fund{\partial M_1}} \in \R(\partial M_1)$. In particular, this implies that $\restr{\rho}{\fund{\partial M_1}} \in H_{1,\pi}$. Moreover, $f_\pi(\restr{\rho}{\fund{\partial M_1}})= \Tr \rho(a_1b_1)=\Tr AB=z$. This implies that $f_\pi$ is surjective.
\end{proof}
\label{lemma: surjective maps}
\end{lemma}

Let us define
$
    S: \mathbb{N}_{\ge 2} \longrightarrow \mathbb{N}
$
as the function that maps a natural number $n\ge 2$ into the smallest number $1\le k \le n$ which maximizes the quantity $\left|\sin \left( \frac{2\pi k}{n}\right)\right|$. It is easy to see that
\begin{equation}
S(n)= \frac{n+x_n}{4} \quad \text{with} \quad x_n = \begin{dcases}
    0 & \text{if } n \equiv_4 0, \\
    -1 & \text{if }n \equiv_4 1, \\
    -2 & \text{if }n\equiv_4 2, \\
    1 & \text{if } n \equiv_4 3.
\end{dcases}
\label{definition of S and x_p}
\end{equation}
Let us assume that $3 \le p_1 \le p_2$. The interval $J_0(p_1,p_2)$ is, by definition, the union of the sub-intervals $I(\nicefrac{k_1}{p_1},\nicefrac{k_2}{p_2})$. According to the \eqref{eq: length of the interval J(a,b)}, the interval $I \left( \nicefrac{S(p_1)}{p_1},\nicefrac{S(p_2)}{p_2}\right)$ is the one of maximum length among these.

\begin{lemma}
Let $3 \le p_1$ and $\alpha \in [0,\pi]$ be the angle that supports the interval $I \left( \frac{S(p_1)}{p_1},\frac{S(p_2)}{p_2}\right)$ as in Definition \ref{defn: angle that supports an interval}, then $\alpha \ge \frac{2\pi}{3}$.
\label{lemma: the interval I(S,S) has angle more than 2pi/3}
\begin{proof}
    According to Lemma \ref{lemma: J0(2,p2) è vuoto}, the interval $I \left( \nicefrac{S(p_1)}{p_1},\nicefrac{S(p_2)}{p_2}\right)$ is nonempty.
    Let $x_{p_{i}}$ be defined as in \eqref{definition of S and x_p} for $i \in \{1,2\}$.
    Definition \ref{defn: interval I(a,b)} implies that the end points the interval $I \left( \frac{S(p_1)}{p_1},\frac{S(p_2)}{p_2}\right)$ are
\begin{multline*}
2\cos \left(\frac{2\pi S(p_2)}{p_2}+\frac{2\pi S(p_1)}{p_1}\right)= 2\cos\left( \pi + \frac{\pi}{2}\left(\frac{x_{p_2}}{p_2}+\frac{x_{p_1}}{p_1}\right)\right)  \quad \text{and} \\ 2\cos \left(\frac{2\pi S(p_2)}{p_2}-\frac{2\pi S(p_1)}{p_1}\right)= 2\cos\left( \frac{\pi}{2}\left(\frac{x_{p_2}}{p_2}-\frac{x_{p_1}}{p_1}\right)\right)
\label{eq: end points of the interval}
\end{multline*}
    Let $(c_1,c_2)$ be an interval in $[-2,2]$. The angle that supports the interval $(c_1,c_2)$ is equal to $|\arccos (\nicefrac{c_1}{2})-\arccos (\nicefrac{c_2}{2})|$. If $c_1=2\cos(\pi \pm \gamma_1)$ and $c_2=2\cos(\pm \gamma_2)$ with $\gamma_1,\gamma_2 \in [0,\nicefrac{\pi}{2}]$, then the angle that supports the interval $(c_1,c_2)$ is $\pi - \gamma_1-\gamma_2$.
    Hence, we obtain
    \begin{align*}
        \alpha\left(I \left( \frac{S(p_1)}{p_1},\frac{S(p_2)}{p_2} \right) \right) &= \pi - \frac{\pi}{2} \left|\frac{x_{p_2}}{p_2}+\frac{x_{p_1}}{p_1} \right| - \frac{\pi}{2} \left| \frac{x_{p_2}}{p_2}-\frac{x_{p_1}}{p_1}\right|  \\[5pt]
        & = \pi - \pi \frac{\left|x_{p_1}p_2+x_{p_2}p_1 \right|+\left|x_{p_1}p_2-x_{p_2}p_1 \right|}{2p_1p_2}\\[5pt]
        & \ge \min \left\{\pi - \pi \frac{|x_{p_1}|}{p_1},\pi - \pi \frac{|x_{p_2}|}{p_2} \right\} \ge \frac{2 \pi}{3}.
    \end{align*}
\end{proof}
\end{lemma}

\begin{lemma}
Let $3 \le p_1$, then $J_0(p_1,p_2)$ is connected.
\label{lemma: J(p_1,p_2) is connected}
\begin{proof}
The intervals $J_0(3,3)=I(\nicefrac{1}{3},\nicefrac{1}{3})$, $J_0(3,4)=I(\nicefrac{1}{3},\nicefrac{1}{4})$, and $J_0(4,4)=I(\nicefrac{1}{4},\nicefrac{1}{4})$ are connected. Hence, let us suppose that $p_2 \ge 5$. We define the interval $J_1$ as
\[
    J_1 \coloneqq \left(2\cos \left( \pi -\frac{2\pi}{p_2}\right) ,2\cos \left( \frac{2\pi}{p_2}\right)\right) = \left(-2\cos \left( \frac{2\pi}{p_2}\right) ,2\cos \left( \frac{2\pi}{p_2}\right)\right). 
\]
\begin{claim}\label{claim intervalli}
    The interval $J_0(p_1,p_2)$ contains the connected interval $J_1$
\end{claim}

We prove the proposition assuming Claim \ref{claim intervalli}, and then we prove Claim \ref{claim intervalli}. Note that, since $p_2 \ge 5$, the interval $J_1$ is nonempty.

Let us assume Claim \ref{claim intervalli} and suppose, by contradiction, that $J_0(p_1,p_2)$ is not connected. Let $J_2$ be a connected component of $J_0(p_1,p_2)$ disjoint from $J_1$. Then
\[
\text{either } J_2 \subseteq \left[ -2 , -2 \cos \left( \frac{2 \pi }{p_2}\right) \right) \quad \text{or } J_2 \subseteq \left(2 \cos \left( \frac{2 \pi }{p_2}\right),2  \right].
\]
This means that $m(J_2)\le 2-2\cos \left( \frac{2 \pi }{p_2}\right)= 4 \sin^2 \left( \frac{\pi }{p_2}\right)$.
Let $i,j \in \mathbb{Z}$ be such that the interval $I\left(\nicefrac{i}{p_1},\nicefrac{j}{p_2}\right)$ has non-zero length.
Let $n \in \mathbb{N}$, we remark that if $k$ is an integer such that $k \notin \frac{n}{2} \mathbb{Z}$, then $\left |\sin (\nicefrac{2\pi k}{n}) \right| \ge \left |\sin (\nicefrac{\pi}{n}) \right|$.
According to the \eqref{eq: length of the interval J(a,b)}, this implies that
\begin{equation*}
    \begin{aligned}
        m \left( I \left(\frac{i}{p_1}, \frac{j}{p_2} \right) \right) = \left| \sin \left( \frac{2\pi i}{p_1}\right)\sin \left( \frac{2\pi j}{p_2}\right)\right| \ge
        4 \left |\sin \left( \frac{\pi}{p_1}\right)\sin \left( \frac{\pi}{p_2}\right) \right| \ge 4 \sin^2 \left( \frac{\pi}{p_2} \right).
    \end{aligned}
\end{equation*}
If $ p_1<p_2$, then this last inequality is strict. In this case, we get a contradiction since we supposed that the connected component $J_2$ has length smaller or equal than $4 \sin^2 \left( \frac{\pi }{p_2}\right)$.

Suppose than that $ p_1=p_2\ge 5$. The end points of the interval $I(\nicefrac{k-1}{p_1},\nicefrac{1}{p_1})$ are 
\[
    2\cos\left( 2\pi \frac{k}{p_1} \right) \quad \text{and} \quad 2\cos\left( 2\pi \frac{k-2}{p_1} \right).
\]
This implies that
\[
    J_3 \coloneqq \bigcup_{k \in \mathbb{Z}} I \left( \frac{k-1}{p_1},\frac{1}{p_1}\right) = \bigcup_{k' =1}^{p_1} \left( 2\cos\left( \frac{2\pi k'}{p_1} \right), 2 \right) \subset J_0(p_1,p_2)
\]
is a connected subinterval in $J_0(p_1,p_2)$.
If $k_1$ and $k_2$ are both not in $\frac{p_1}{2}\mathbb{Z}$, then the end points of the interval $I(\nicefrac{k_1}{p_1},\nicefrac{k_2}{p_1})$ are
\[
2\cos\left(2\pi\frac{ k_1+k_2}{p_1}\right) \quad \text{and} \quad 2\cos\left(2\pi\frac{ k_1-k_2}{p_1}\right).
\]
This implies that the end points of the interval $I(\nicefrac{k_1}{p_1},\nicefrac{k_2}{p_1})$ are both of the form $2\cos(\nicefrac{2\pi k'}{p_1})$, with $k' \in \mathbb{Z}$.
This implies that
for every $k_1,k_2\in \mathbb{Z}$, we have $I(\nicefrac{k_1}{p_1},\nicefrac{k_2}{p_2}) \subseteq J_3$. And therefore
$J_0(p_1,p_1) \subseteq J_3$ that brings us to $J_0(p_1,p_2)=J_3$. The conclusion holds by the fact that $J_3$ is connected.

\begin{proof}[Proof of Claim \ref{claim intervalli}]
    We shall prove that $J_1 \subseteq J_0(p_1,p_2)$ for $p_2 \ge 5$. Consider the union
\begin{equation}
\bigcup_{k_2 \in \mathbb{Z}} I \left(\frac{S(p_1)}{p_1},\frac{k_2}{p_2} \right) = \bigcup_{j \in \mathbb{Z}}I \left ( \frac{S(p_1)}{p_1}, \frac{S(p_2)+j}{p_2}\right).
\label{eq: la uso una volta: unione intervalli}
\end{equation}
Let $\theta_{1},\theta_{2} \in [0,\pi]$, be such that $I \left (\nicefrac{S(p_1)}{p_1}, \nicefrac{S(p_2)}{p_2}\right) = \left(2\cos\theta_{1},2\cos \theta_{2}\right)$.
The end points of $I \left ( \frac{S(p_1)}{p_1}, \frac{S(p_2)+j}{p_2}\right)$ are:
\begin{equation*}
\begin{split}
    \left\{ 2 \cos \left( 2\pi \frac{p_1 S(p_2)+p_2 S(p_1)}{p_1p_2} + \frac{2\pi j}{p_2}\right), 2 \cos \left( 2\pi \frac{p_1 S(p_2)- p_2 S(p_1)}{p_1p_2} + \frac{2\pi j}{p_2}\right) \right\} = \\[5pt]
= \left\{ 2 \cos \left( \theta_1 + \frac{2\pi j}{p_2}\right), 2 \cos \left( \theta_2 + \frac{2\pi j}{p_2}\right) \right\}.
\end{split}
\label{eq: endpoint of the interval in the proof with j}
\end{equation*}
By Lemma \ref{lemma: the interval I(S,S) has angle more than 2pi/3}, $|\theta_1 - \theta_2| \ge \frac{2\pi}{3}> \frac{2\pi}{p_2}$. In particular $\pi \ge \theta_1 > \theta_2 + \nicefrac{2\pi}{p_2} > \theta_2 \ge 0$.
This implies that
\[
    2\cos \theta_2 > 2\cos \left( \theta_2 + \frac{2\pi}{p_2}\right)> 2\cos \theta_1.
\]
Furthermore, if $\theta_1 < \pi - \nicefrac{\pi}{p_2}$, then $\pi - \nicefrac{\pi}{p_2} < \theta_1 + \nicefrac{2\pi}{p_2}< \pi + \nicefrac{\pi}{p_2}$. Thus,
\[
 2\cos \theta_2 > 2\cos \left( \theta_2 + \frac{2\pi}{p_2}\right)> 2\cos \theta_1 > 2\cos \left( \theta_1 + \frac{2\pi}{p_2}\right).
\]
Let $j \in \mathbb{N}$ be such that $\theta_1 + \nicefrac{2\pi j}{p_2} < \pi - \nicefrac{\pi}{p_2}$, then
\begin{align}
     2\cos\left(\theta_{2}+ \frac{2\pi (j+1)}{p_2}\right)> 2\cos\left(\theta_{1}+ \frac{2\pi j}{p_2}\right) > 2\cos\left(\theta_{1}+ \frac{2\pi (j+1)}{p_2}\right).
    \label{eq: una volta, i punti finali}
\end{align}
Notice that the left and right sides of identity \eqref{eq: una volta, i punti finali} are the end points of $I \left ( \frac{S(p_1)}{p_1}, \frac{S(p_2)+j+1}{p_2}\right)$. Furthermore, the central term the identity \eqref{eq: una volta, i punti finali} is an end point of $I \left ( \frac{S(p_1)}{p_1}, \frac{S(p_2)+j}{p_2}\right)$. Hence, identity \eqref{eq: una volta, i punti finali}
implies that $I\left( \frac{S(p_1)}{p_1},\frac{S(p_2)+j}{p_2} \right) \cap I\left( \frac{S(p_1)}{p_1},\frac{S(p_2)+j+1}{p_2} \right)$ is nonempty.
Thus, this implies that $I\left( \frac{S(p_1)}{p_1},\frac{S(p_2)+j}{p_2} \right) \cup I\left( \frac{S(p_1)}{p_1},\frac{S(p_2)+j+1}{p_2} \right)$ is connected.
Similarly, let $j \in \mathbb{N}$ be such that $\theta_2 + \nicefrac{2\pi j}{p_2} \ge \nicefrac{\pi}{p_2}$, then
\[
    2\cos\left(\theta_{1}+ \frac{2\pi (j-1)}{p_2}\right)>2\cos\left(\theta_{2}+ \frac{2\pi j}{p_2}\right) > 2\cos\left(\theta_{2}+ \frac{2\pi (j-1)}{p_2}\right).
\]
As before, this implies that one of the end points of $I \left ( \frac{S(p_1)}{p_1}, \frac{S(p_2)+j}{p_2}\right)$ is in $I \left ( \frac{S(p_1)}{p_1}, \frac{S(p_2)+j-1}{p_2}\right)$.
So $I\left( \frac{S(p_1)}{p_1},\frac{S(p_2)+j}{p_2} \right) \cap I\left( \frac{S(p_1)}{p_1},\frac{S(p_2)+j-1}{p_2} \right)$ is nonempty.
Furthermore, $I\left( \frac{S(p_1)}{p_1},\frac{S(p_2)+j}{p_2} \right) \cup I\left( \frac{S(p_1)}{p_1},\frac{S(p_2)+j-1}{p_2} \right)$ is connected.
Thus, the union in \eqref{eq: la uso una volta: unione intervalli} contains the connected interval $J_1$.
\end{proof} 
\end{proof}
\end{lemma}

The set $J_\pi(p_1,p_2)$ is not connected in general, for instance, $J_\pi(4,4) =(-2,0) \cup(0,2)$. The author thinks that this is the only case in which this happens.

\begin{lemma}
    Let $3 \le p_1$, the set $J_\pi(p_1,p_2)$ does not contain the element $0 \in [-2,2]$ if and only if $p_1=p_2=4$.
    \label{lemma: the interval Jpi(p1,p2) contains 0 if and only if p1,p2 is not 4,4}
    \begin{proof}
        The interval $I \left( \nicefrac{k_1}{2p_1},\nicefrac{k_2}{2p_2}\right)$ has length equal to
        \[
            m\left( I \left( \frac{k_1}{2p_1},\frac{k_2}{2p_2}\right)\right) = 4 \left| \sin \left( \frac{\pi k_1}{p_1}\right)\sin \left( \frac{\pi k_2}{p_2} \right)\right|.
        \]
        It is easy to see that if $p_1 \ge 2$ there exists an odd integer $k$ in the range $\nicefrac{p_1}{4} \le k \le \nicefrac{3p_1}{4}$. For such a $k$ we have 
        \[
            \left| \sin \left( \frac{2\pi k}{2p_1}\right) \right|=\left| \sin \left( \frac{ \pi k}{p_1}\right) \right| \ge  \frac{\sqrt{2}}{2}.
        \]
        Moreover, if $p_1 \neq 4$, we can choose $k$ to make this inequality strict. 
        Thus, if $(p_1,p_2)\neq (4,4)$, then there exists an interval 
        $
        I\left( \nicefrac{k_1}{2p_1},\nicefrac{k_2}{2p_2}\right)\subseteq J_\pi(p_1,p_2)
        $
        whose length is strictly bigger than $2$. This implies that $0 \in I\left( \nicefrac{k_1}{2p_1},\nicefrac{k_2}{2p_2}\right)$. A direct computation shows that $0 \notin J_\pi(4,4)$ and this gives the conclusion.
    \end{proof}
\end{lemma}

\begin{lemma}
    Let $3 \le p_1$ and let $g_1=\gcd(p_1,p_2)$, then
    \[
        \mathcal{S}(p_1,p_2) \coloneqq \left\{ k_1p_2 + k_2 p_1  \, \middle| \, k_i \in \mathbb{Z}, k_i \notin \frac{p_i}{2}\mathbb{Z}\right\} = 
        \begin{dcases}
            \mathbb{Z} \setminus \left( \frac{p_1}{2}\mathbb{Z} \cup \frac{p_2}{2}\mathbb{Z}\right)   & \text{if } g_1=1, \\[5pt]
            2\mathbb{Z} \setminus \left( p_1\mathbb{Z} \cup p_2\mathbb{Z}\right)   & \text{if } g_1=2, \\[5pt]
            g_1\mathbb{Z}    & \text{if } g_1 \ge 3.
        \end{dcases}
    \]
    \label{lemma: description of the set S(p1,p2)}
    \begin{proof}
        
        Let us start with $g_1 = 1$. Let $x=k_1p_2+k_2p_1 \in \mathcal{S}(p_1,p_2)$. Let us suppose by contradiction that $x \in \frac{p_1}{2}\mathbb{Z}$. This implies that $k_1 p_2 \in \frac{p_1}{2}\mathbb{Z}$. Since $\gcd(p_1,p_2)=1$, we obtain that $k_1  \in \frac{p_1}{2} \mathbb{Z}$ that is a contradiction.
        Similarly, if $x \in \frac{p_2}{2}\mathbb{Z}$, then $k_2p_1 \in \frac{p_2}{2} \mathbb{Z}$. Hence, $k_2 \in \frac{p_2}{2} \mathbb{Z}$ that is a contradiction. Conversely, let $x \in \mathbb{Z} \setminus \left( \frac{p_1}{2}\mathbb{Z} \cup \frac{p_2}{2}\mathbb{Z}\right)$. Let us write $x=k_1p_2+k_2p_1$. Let us suppose, by contradiction, that $k_1 \in \frac{p_1}{2}\mathbb{Z}$. This implies that $x \in\frac{p_1}{2}\mathbb{Z}$, that is a contradiction.

        Let us suppose that $g_1=2$. We write $p_1=2n_1$ and $p_2=2n_2$, with $\gcd(n_1,n_2)=1$. Let $x=2k_1n_2+2k_2n_1 \in \mathcal{S}(p_1,p_2)$. We note that since $p_1 \ge 3$, the set $ \mathcal{S}(p_1,p_2)$ is nonempty. Let us suppose, by contradiction, that $x \equiv_{2n_1} 0$. This implies that $2k_1n_2 \equiv_{2n_1} 0$. We obtain that $k_1 \equiv_{n_1} 0$, this is a contradiction.
        Similarly, if $x \equiv_{2n_2} 0$, then $2k_2n_1 \equiv_{2n_2}0$. Thus, $k_2 \equiv_{n_2} 0$ and this is a contradiction.
        Conversely, let $x \in 2\mathbb{Z} \setminus \left( p_1\mathbb{Z} \cup p_2\mathbb{Z}\right)$. Let us write $x$ as $k_1p_2+k_2p_1$ and suppose that $k_1 \equiv_{n_1}0$. This assumption implies that $2k_1 \equiv_{2n_1}0$ and hence $x \equiv_{p_1}0$ that is a contraction. Similarly, if $k_2 \equiv_{n_2}0$, then $x \equiv_{p_2}0$, and here is our contraction.
        
        Let us move to the case $g_1 \ge 3$. Notice that $\mathcal{S}(p_1,p_2) \subseteq g_1 \mathbb{Z}$. Let us suppose that $x \in g_1 \mathbb{Z}$ and let us write $x=p_2k_1 + p_1k_2$. More generally, we have that for every $n \in \mathbb{N}$
        \[
            x= p_1 k_{2,n} + p_2 k_{1,n} \quad \text{with} \quad k_{1,n}= k_1 + n\left(\frac{p_1}{g_1}\right) \quad \text{and} \quad k_{2,n}= k_2 - n\left(\frac{p_2}{g_1}\right).
        \]
        We need to show that there exists a $n \in \mathbb{N}$ such that
        \begin{equation}
            k_{1,n} \notin \frac{p_1}{2} \mathbb{Z} \quad \text{and} \quad k_{2,n} \notin \frac{p_2}{2} \mathbb{Z}.
            \label{eq: desired property}
        \end{equation}
        If $k_{1,0}=k_1$ and $k_{2,0}=k_2$ are not in $\frac{p_1}{2} \mathbb{Z}$ and $\frac{p_2}{2} \mathbb{Z}$ respectively, we get the conclusion. Suppose then without loss of generality that $k_{1,0} \in \frac{p_1}{2} \mathbb{Z}$. Then $k_{1,1}=k_1 + \frac{p_1}{g_1} \notin \frac{p_1}{2} \mathbb{Z}$. If the couple $(k_{1,1},k_{2,1})$ has the property in \eqref{eq: desired property}, then we get the conclusion. Let us suppose that $k_{2,1}=k_2+ \frac{p_2}{g_1} \in \frac{p_2}{2} \mathbb{Z}$. Then, since $g_1 \ge 3$, we obtain $k_{1,2} \notin \frac{p_2}{2} \mathbb{Z}$ and $k_{2,2}\notin \frac{p_2}{2} \mathbb{Z}$. This implies that $x \in \mathcal{S}(p_1,p_2)$ and this completes the proof.
        \end{proof}
\end{lemma}

Let $\mathcal{S}(p_1,p_2)$ be as in Lemma \ref{lemma: description of the set S(p1,p2)}, we define
$x_{\min},x_{\max} \in [-2,2]$ as
\begin{equation}
x_{\min} = \min_{x \in \mathcal{S}(p_1,p_2)} \left\{ 2 \cos \left( \frac{2\pi x}{p_1p_2}\right)\right\} \quad \text{and} \quad x_{\max}=\max_{x \in \mathcal{S}(p_1,p_2)} \left\{ 2 \cos \left( \frac{2\pi x}{p_1p_2}\right)\right\}.
\label{eq: explitic expression of J0(p1,p2)}
\end{equation}

\begin{lemma}\label{lemma: computation of J_0(p1,p2)}
    Let $x_{\min}, x_{\max} \in [-2,2]$ be as in \eqref{eq: explitic expression of J0(p1,p2)}, then $J_0(p_1,p_2)=\left( x_{\min}, x_{\max} \right)$.
    \begin{proof}
        Let $k_1,k_2 \in \mathbb{Z}$ be two integers such that $k_i \notin \frac{p_i}{2}\mathbb{Z}$. As we proved in Lemma \ref{lemma: J0(2,p2) è vuoto}, the interval $I(\nicefrac{k_1}{p_1},\nicefrac{k_2}{p_2})$ is nonempty and 
        $I(\nicefrac{k_1}{p_1},\nicefrac{k_2}{p_2}) \subseteq J_0(p_1,p_2)$. In particular,
        \[
            \partial I\left(\frac{k_1}{p_1},\frac{k_2}{p_2}\right) = \left\{ 2\cos \left( \frac{2\pi (k_1p_2 \pm k_2p_1)}{p_1p_2}\right) \right\}\subset \overline{J_0(p_1,p_2)},
        \]
        Thus,
        the interval $\overline{J_0(p_1,p_2)}$ must contain $2 \cos(\nicefrac{2\pi x}{p_1p_2})$ for all $x \in \mathcal{S}(p_1,p_2)$. Since $J_0(p_1,p_2)$ is connected by Lemma \ref{lemma: J(p_1,p_2) is connected}, we obtain that $\overline{J_0(p_1,p_2)}=\left[ x_{\min}, x_{\max} \right]$ and hence $J_0(p_1,p_2)=\left( x_{\min}, x_{\max} \right)$.
    \end{proof}
\end{lemma}

\begin{cor}
    Suppose $p_1 \ge 3$ and let $\alpha$ be the angle the supports the interval $J_0(p_1,p_2)$, then $\alpha \ge \pi - \frac{4\pi g_1}{p_1p_2} \ge \pi - \frac{4\pi}{p_2}$.
    \label{cor: angle supporting J0(p1,p2) is bigger that pi-4pi/p2}
\end{cor}

\begin{lemma}\label{lemma: the interval Jpi(2,p2) is nested and connected}
     If $p_2\ge 2$, then $J_\pi(2,p_2)=I\left(\frac{1}{4},\frac{k}{2p_2}\right)$, where $k$ is an odd number satisfying $1\le k\le {p_2}$ that maximizes $|\sin \left( \nicefrac{\pi k}{p_2}\right)|$. In particular $J_\pi(2,p_2)$ is connected and supported by an angle greater than or equal to $\pi - \frac{2\pi}{p_2}$.
    \begin{proof}
        Let $I\left(\frac{1}{4},\frac{k_2}{2p_2}\right)$ be an interval of $J_\pi(2,p_2)$. The end points of $I\left(\frac{1}{4},\frac{k_2}{2p_2}\right)$ are
    \begin{equation}
    \begin{aligned}
    \partial I\left(\frac{1}{4},\frac{k_2}{2p_2}\right)& = \left\{ 2 \cos \left( 2\pi \frac{ 2p_2 \pm 4k_2 }{8p_2} \right) \right\} \\
    &= \left\{ 2 \cos \left( \frac{\pi}{2} + \frac{\pi k_2}{p_2}\right) , 2 \cos \left( \frac{\pi}{2} - \frac{\pi k_2}{p_2}\right)\right\}
    = \left\{ \pm 2 \sin \left( \frac{\pi k_2}{p_2}\right) \right\}.
    \end{aligned}
    \label{eq: la uso una volta nested p1=2}
    \end{equation}
    Hence, all intervals in $J_\pi(2,p_2)$ are nested, and we conclude that $J_\pi(2,p_2)$ is connected.
    Let $\alpha \in [0,\pi]$ be the angle supporting $J_\pi(2,p_2)$. The identity \eqref{eq: la uso una volta nested p1=2} implies that
    \[
        \alpha \ge \pi - \frac{\pi}{p_2} - \frac{\pi}{p_2} = \pi - \frac{2\pi}{p_2}.
    \]
    \end{proof}
\end{lemma}

\begin{cor}\label{cor: l'angolo di Jpi(2,p2) è sempre maggione di 1/2pi}
    If $p_2\ge 2$, then $J_\pi(2,p_2)$ is supported by an angle greater than or equal to $\nicefrac{\pi}{2}$. Moreover, $J_\pi(2,4) \subseteq J_\pi(2,p_2)$ for every $p_2 \ge 2$. This last inclusion is strict if $p_2\neq 4$.
    \begin{proof}
    Let $\alpha_{p_2} \in [0,\pi]$ be the angle supporting $J_\pi(2,p_2)$.
    A direct computation shows that $\alpha_{p_2}>\nicefrac{\pi}{2}$ if $p_2=2,3$. Moreover, if $p_2 \ge 4$, then $\alpha_{p_2} \ge \nicefrac{\pi}{2}$ by Lemma \ref{lemma: the interval Jpi(2,p2) is nested and connected}.

    Let $p_2 \ge 2$. Lemma \ref{lemma: the interval Jpi(2,p2) is nested and connected} implies that there exists a real number $x \in (0,2]$ such that $J_\pi(2,p_2) = (-x,x)$. Since $\alpha(J_\pi(2,p_2)) \ge \nicefrac{\pi}{2}$ and $\alpha(J_\pi(2,4))=\nicefrac{\pi}{2}$, we obtain that $J_\pi(2,4)$ is the smallest possible interval among the ones of the form $J_\pi(2,p_2)$. This also implies that $J_\pi(2,4) \subseteq J_\pi(2,p_2)$. The identity \eqref{eq: la uso una volta nested p1=2} implies that if $p_2 \neq 4$, then $J_\pi(2,4) \neq J_\pi(2,p_2)$. Thus if $p_2 \neq 4$, then $J_\pi(2,4) \subset J_\pi(2,p_2)$.
    \end{proof}
\end{cor}

\section{The intersections between \texorpdfstring{$P_1$}{P\_1} and \texorpdfstring{$P_2$}{P\_2}}
\label{sec: P1 cap P2}

    Recall that we assumed $g_1 \le g_2$.
\begin{prop}
    The set $P_1 \cap P_2$ is empty if and only if one of the following conditions holds:
    \begin{enumerate}
        \item $g_1=1$;
        \item $g_1=2$ and $o_1=1$;
        \item $g_2=2$ and $o_2=1$;
        \item $g_1=2$, $o_2$ odd and $\Delta(h_1,\lambda_{M_2})$ is even.
    \end{enumerate}
    \label{prop: when P1 P2 contains SU(2)-abelian representations}
    \begin{proof}
        Let us suppose that either the condition $(1)$ or $(2)$ holds. By Lemma \ref{lemma: description of the repersentation P1}, the set $P_1$ is empty. Thus, $P_1 \cap P_2= \emptyset$.
        Similarly, if the condition $(3)$ holds, then Lemma \ref{lemma: description of the repersentation P1} implies that $P_2 = \emptyset$ and then, $P_1 \cap P_2 = \emptyset$.
        
        Let us suppose that the condition $(4)$ holds. Without loss of generality, we can suppose that conditions $(1)$, $(2)$, and $(3)$ do not hold. This implies that $g_1=o_1=2$. Let $\rho \in \mathcal{R}(M)$ be such that $\rho(\fund{\Sigma}) \subset \ZSU$.  
        We will show that $\rho_1(h_1)=1$. According to Lemma \ref{lemma: description of the repersentation P1}, this implies that $\rho_1$ is a central. Thus, $P_1 \cap P_2= \emptyset$.
        
        Since $\rho(\fund{\Sigma}) \subset \ZSU$, the restrictions $\rho_1=\restr{\rho}{\fund{M_1}}$ and $\rho_2=\restr{\rho}{\fund{M_2}}$ are abelian.
        Let $\xi$ be an oriented simple closed curve in $\partial M_2$ such that $\{\lambda_{M_2},\xi\}$ generates $\fund{\partial M_2}$.
        The representation $\rho$ maps both $\lambda_{M_2}$ and $\xi$ into $\pm 1$. Since $\rho_2$ is abelian, as a consequence of Corollary \ref{cor: an abelian representation extends if and only if the rational longitude is zero}, we obtain that 
    \[
        \rho(\lambda_{M_2})^{o_2}=\rho_2(\lambda_{M_2})^{o_2}=1.
    \]
    Since $o_2$ is odd by hypothesis, $\rho(\lambda_{M_2})=1$.
    There exist integers $n,m \in \mathbb{Z}$, with $|n|=\Delta(\xi,h_1)$ and $|m|=\Delta(\lambda_{M_2},h_1)$, such that
        \begin{equation}\label{eq: la uso una volta rho(h1)=1}
            h_1 = n \cdot \lambda_{M_2}+m \cdot \xi \subset \partial M_1, \quad \text{and hence} \quad \rho_1(h_1) = \cancel{\rho_2(\lambda_{M_2})^n}\rho_2(\xi)^m=\rho_2(\xi)^m.
        \end{equation}    
        Since $\Delta(h_1,\lambda_{M_2})=|m|$ is even,
        we get that $\rho_1(h_1)=1$.

        Conversely, let us suppose that $ P_1\cap P_2 = \emptyset$.
        This implies that if the representation $\rho \in \mathcal{R}(M)$ is such that $\rho(\fund{\Sigma}) \subset \ZSU$ then, either $\rho_1$ or $\rho_2$ is central.
        Let $\eta \colon \fund{\Sigma} \to \{1\}$ be the trivial representation. Then, we have $\eta(\lambda_{M_1})^{o_1}=\eta(\lambda_{M_2})^{o_2}=1$. According to Corollary \ref{cor: an abelian representation extends if and only if the rational longitude is zero}, the representation $\eta$ extends to a representation $\rho \colon \fund{M} \to SU(2)$. Since either $\rho_1$ or $\rho_2$ is $SU(2)$-central, Lemma \ref{lemma: g1 >=3, allora non-central} implies that either $g_1 \le 2$ or $g_2 \le 2$. Since we supposed that $g_1 \le g_2$, we obtain $g_1 \le 2$.
        
        We complete the proof by proving that if conditions $(1)$, $(2)$, and $(3)$ do not hold, then condition $(4)$ holds.
        If conditions $(1)$, $(2)$, and $(3)$ do not hold, then, since $g_1\le 2$, we get that $o_1=g_1=2$ and either $o_2=g_2=2$ or $g_2 \ge 3$.

         Suppose first that $o_2 \equiv_20$. Then $g_2 \equiv_20$ by Lemma \ref{lemma: description of the rational longitude}. Since $o_1$ and $o_2$ are both even, Corollary \ref{cor: an abelian representation extends if and only if the rational longitude is zero} implies that every central representation $\fund{\Sigma} \to \ZSU$ extends to both $\fund{M_1}$ and $\fund{M_2}$. Let $\eta \colon \fund{\Sigma} \to \ZSU$ be a representation such that $\eta(h_1)=\eta(h_2)=-1$.
        Lemma \ref{lemma: g1 pari e h1=-1 allora non-central} implies that there exist two non-central representations $\rho_1 \in\mathcal{R}(M_1)$ and $\rho_2 \in \mathcal{R}(M_2)$ such that $\restr{\rho_1}{\fund{\Sigma}}\equiv\eta$ and $\restr{\rho_2}{\fund{\Sigma}}\equiv\eta$. This implies that $\eta \in P_1$ and $\eta \in P_2$.
        Thus, $P_1 \cap P_2 \neq \emptyset$, and this is a contradiction. This implies that if $P_1 \cap P_2= \emptyset$, and the conditions $(1)$, $(2)$, and $(3)$ do not hold, then $o_2$ is odd.
        
        Suppose next that $o_2$ is odd. Since $2=g_1 \le g_2$ and condition $(3)$ do not hold, we obtain that $g_2 \ge 3$. Lemma \ref{lemma: g1>=3 the trivial representation is not central} implies that there exists a representation $\eta \colon \fund{\Sigma} \to \ZSU$ that extends to $\fund{M_2}$ in a non-central way.
        In particular $\eta \in P_2$.
        Since $o_1=2$, the representation $\eta$ extends to $\fund{M_1}$ as well. Since $P_1 \cap P_2 = \emptyset$ by hypothesis, if $\rho_1 \colon \fund{M_1} \to SU(2)$ is an extension of $\eta$, then $\rho_1$ is central. Notice that, since $o_1=2$, then $\eta$ extends to $\fund{M_1}$ by Corollary \ref{cor: an abelian representation extends if and only if the rational longitude is zero}. Since $g_1=o_1=2$, Lemma \ref{lemma: description of the repersentation P1} implies that $\rho_1$ is central if and only if $\rho_1(h_1)=1$.
        The identity \eqref{eq: la uso una volta rho(h1)=1} implies that this happens if and only if $\Delta(h_1,\lambda_{M_2})$ is even. Thus, $P_1 \cap P_2= \emptyset$, and the conditions $(1)$, $(2)$, and $(3)$ do not hold, then condition $(4)$ applies. This concludes the theorem.
        \end{proof}
\end{prop}

It is conjectured that $SU(2)$-abelian rational homology $3$-spheres are Heegaard Floer L-spaces (see Conjecture \ref{conj: SU(2)-abelian implies L-space}). It is known that the converse of this conjecture is false, as shown below.
By \cite[Theorem 5]{Boyer2011OnLA}, the gluing of two twisted I-bundles over the Klein bottle is an L-space if and only if the gluing map does not identify the two rational longitudes. However, we prove in Example \ref{exmp: two evil manifolds} that such a manifold is not $SU(2)$-abelian.

\begin{exmp}
    Let $M_1$ and $M_2$ be two copies of the twisted I-bundle over the Klein bottle. Hence, $M_1=M_2=\mathbb{D}^2(\nicefrac{2}{1},\nicefrac{2}{1})$. Let $\varphi:\partial M_1 \to \partial M_2$ be a diffeomorphism.
    Proposition \ref{prop: when P1 P2 contains SU(2)-abelian representations} implies that the graph manifold $M_1 \cup_{\Sigma}M_2=M_1 \cup_{\varphi} M_2$ admits an irreducible $SU(2)$-representation whose restriction to $\fund{\Sigma}$ is in $P_1\cap P_2$. 
    This implies that it does not exist a diffeomorphism $\varphi: \partial M_1 \to \partial M_2$ such that $M_1\cup_{\varphi}M_2$ is $SU(2)$-abelian. 
    \label{exmp: two evil manifolds}
\end{exmp}

\section{The intersection between \texorpdfstring{$H_1$}{H\_1} and \texorpdfstring{$A_2$}{A\_2} and between \texorpdfstring{$H_2$}{H\_2} and \texorpdfstring{$A_1$}{A\_1}}
\label{sec: A1 e H2}
We use the coordinates $(\theta_1,\psi_1)$ for the space $\R(\partial M_1)$ with respect to the basis $\{\mu_1,h_1\}$ as in \eqref{eq: coordinate for the torus mu1 and h1}. 
Choose integers $n,m \in \mathbb{Z}$ such that
\begin{align}
    \lambda_{M_2}= n h_1+ m \mu_1 \subset \partial M_1.
    \label{eq: la uso una volta, longitudine}
\end{align}
Then $|n|=\Delta(\lambda_{M_2},\mu_1)$ and $|m|=\Delta(\lambda_{M_2},h_1)$.
Define a subspace $L \subset \R(\partial M_1)$ by
\[
    L=\left\{ (\theta_1,\psi_1) \, \middle| \, o_2 n \psi_1 + o_2 m \theta_1\equiv_{2\pi}0\right\} \subset \R(\partial M_1).
\]
Let $\eta \colon \fund{\partial M_1} \to SU(2)$ be a representation such that $\eta \in L$, then
\[
    \eta(\lambda_{M_2})^{o_2}= \left(\eta(h_1)^n\eta(\mu_1)^m \right)^{o_2}=
    e^{i \left(o_2 n \psi_1 + o_2m \theta_1\right)}=1.
\]
Therefore, as a result of Corollary \ref{cor: Y SU(2) allora i due pezzi sono SU(2)-abelian}, we obtain that the representation $\eta$ is in $A_2$.
In fact, $L$ is the set $A_2$ in $(\theta_1,\psi_1)$ coordinates.
Let us suppose that $\Delta(\lambda_{M_2},h_1)\neq 0$, then there are
$o_2 \Delta(\lambda_{M_2},h_1)$ equally spaced
intersections between $A_2$ and the line $\{\psi_1=0\}$.
Similarly, if $\Delta(\lambda_{M_2},h_1)\neq 0$, then there are $o_2 \Delta(\lambda_{M_2},h_1)$ equally spaced intersection points of $A_2$ with the line $\{\psi_1=\pi\}$.
If $p_1=2$, then $\alpha(J_\pi(2,p_2))\ge \nicefrac{\pi}{2}$ by Corollary \ref{cor: l'angolo di Jpi(2,p2) è sempre maggione di 1/2pi}. Furthermore, if $p_1 \ge 3$, then $\alpha(J_{0}(p_1,p_2))\ge \nicefrac{2\pi}{3} > \nicefrac{\pi}{2}$ by Lemma \ref{lemma: the interval I(S,S) has angle more than 2pi/3}. Thus, according to Lemma \ref{lemma: surjective maps}, either $H_{1,0}$ or $H_{1,\pi}$ contains a connected subset of length at least $\nicefrac{\pi}{2}$.
Hence if $o_2 \Delta(\lambda_{M_2},h_1) \ge 5$, then $A_2$ intersects $H_{1}= H_{1,0}\cup H_{1,\pi}$ at least once. We conclude that if $H_1 \cap A_2$ is empty, then $o_2 \Delta(\lambda_{M_2},h_1) \le 4$.

\begin{prop}\label{prop: H1 and A2 sono disjoint if and only if}
    The intersection $H_1 \cap A_2$ is empty if and only if one of the following holds:
    \begin{itemize}
        \item $\Delta(\lambda_{M_2},h_1)=0$, $p_1=2$ and $o_2 \equiv_2 1$ ;
        \item $\Delta(\lambda_{M_2},h_1)=1$ and $o_2\le 2$
        \item $\Delta(\lambda_{M_2},h_1)=1$, $o_2=3$, $p_1=p_2=3$, and $\Delta(\lambda_{M_2},\lambda_{M_1})\equiv_2 0$,
        \item $\Delta(\lambda_{M_2},h_1)=2$, $p_1=p_2=4$ and $o_2=1$,
        \item $\Delta(\lambda_{M_2},h_1)=3$, $p_1=p_2=3$, $o_2=1$, and $\Delta(\lambda_{M_2},\lambda_{M_1})\equiv_2 0$,
        \item $\Delta(\lambda_{M_2},h_1)=4$, $p_1=2$, $p_2=4$, and $o_2=1$.
    \end{itemize}
    \begin{proof}
        We assume that $\fund{M_1}$ and $\fund{M_2}$ are presented as in \eqref{eq: presentation of the fundamental group of M1 and M2} with all the $q_i$ odd. We recall that Lemma \ref{lemma: surjective maps} states that the maps
        \begin{alignat*}{2}
            f_0 \colon  H_{1,0}  &\to J_0(p_1,p_2) \quad \quad \text{and} \quad \quad f_{\pi} \colon H&_{1,\pi} \to  J_\pi(p_1,p_2) \\
            \eta &\mapsto \Tr\eta(\mu_1) &\eta \mapsto \Tr\eta(\mu_1) 
        \end{alignat*}
        are surjective. If $\eta \in \R (\partial M_1)$ has coordinates equal to $(\theta_1,\psi_1)$, then $\Tr \eta (\mu_1)= 2\cos \theta_1$.
        
        As we stated before, if $o_2 \Delta(\lambda_{M_2},h_1) \ge 5$, then $H_1\cap A_2$ is nonempty. Hence, we can suppose that $o_2\Delta(\lambda_{M_2},h_1) \le 4$. Thus, $0 \le\Delta(\lambda_{M_2},h_1) \le 4$ and if this latter is positive, $o_2 \le \frac{4}{\Delta(\lambda_{M_2},h_1)}$. Notice that the quantities $\Delta(\lambda_{M_2},h_1)$ and $\Delta(\lambda_{M_2},\mu_1)$ are coprime.
        
        If $o_2 \Delta(\lambda_{M_2},h_1)=3$, then either $\Delta(\lambda_{M_2},h_1)=1$ and $o_2=3$ or $\Delta(\lambda_{M_2},h_1)=3$ and $o_2=1$. We divide the proof in five cases:
        \begin{enumerate}[label=\roman*)]
        \item $\Delta(\lambda_{M_2},h_1)=0$,
        \item $\Delta(\lambda_{M_2},h_1)=1$ and $o_2 \neq 3$,
        \item $\Delta(\lambda_{M_2},h_1)=2$,
        \item $\Delta(\lambda_{M_2},h_1)=4$,
        \item $o_2 \Delta(\lambda_{M_2},h_1)=3$.
        \end{enumerate}
        
        \textbf{Case \romannum{1})} $\boldsymbol{\Delta(\lambda_{M_2},h_1)=0}$. In this case
        \[
           A_2 = \left\{\left(\theta_1,\varphi_1\right) \middle| \theta_1 \in [0,2\pi] \text{ and } o_2 \Delta(\lambda_{M_2},\mu_1) \psi_1 \equiv_{2\pi} 0\right\} \subset \R(\partial M_1).
        \]
        In particular $A_2 \cap H_{1,0}=H_{1,0}$. Thus, $A_2 \cap H_{1,0}= \emptyset$ if and only if $H_{1,0}$ is empty. By Lemma \ref{lemma: J0(2,p2) è vuoto} and Lemma \ref{lemma: surjective maps}, this happens if and only if $p_1=2$. If $p_1=2$, then $H_1 \cap A_2$ is empty if and only if $H_{1,\pi}\cap A_2$ is empty as well.
        Since, by Lemma \ref{lemma: J0(2,p2) è vuoto} and Lemma \ref{lemma: surjective maps}, the set $H_{1,\pi}$ is not empty, 
        $H_{1,\pi}\cap A_2$ is empty if and only if $\psi_1=\pi$ is not a solution of
        \[
            o_2 \Delta(\lambda_{M_2},\mu_1) \psi_1 \equiv_{2\pi} 0.
        \]
        Since $\Delta(\lambda_{M_2},h_1)$ is zero, $\Delta(\lambda_{M_2},\mu_1)$ is odd. Therefore, $A_2\cap H_1 = \emptyset$ if and only if $p_1=2$ and $o_2\equiv_2 1$.

        \textbf{Case \romannum{2})} $\boldsymbol{\Delta(\lambda_{M_2},h_1)=1}$ \textbf{and} $\boldsymbol{o_2 \neq 3}$. Thus, $o_2 \le 4$.
        The intersection between $\{\psi_1 =0\}$ and $A_2$ are
        \[
            (\theta_1,\psi_1)= \left(\frac{2\pi k}{o_2} , 0  \right) \quad \text{with} \quad k \in \{1, \cdots, o_2\}.
        \]
        If $o_2 \le 2$, then $H_{1,0} \cap A_2 = \emptyset$ as a consequence of Corollary \ref{lemma: H1 doesn't contain -2 and 2}.
        
       If $p_1 \ge 3$, then $J_0(p_1,p_2)$ is supported by an angle $\alpha \ge \nicefrac{2\pi}{3}$ by Lemma \ref{lemma: the interval I(S,S) has angle more than 2pi/3}. Thus, $0 \in J_0(p_1,p_2)$. Lemma \ref{lemma: surjective maps} implies that the set $H_{1,0}$ contains at least one of the points $(\nicefrac{\pi}{2},0)$ and $(\nicefrac{3\pi}{2},0)$. This implies that if $o_2=4$, then $H_{1,0}\cap A_2 = \emptyset$ if and only if $H_{1,0}=\emptyset$. This happens if and only if $p_1=2$.
        
        The intersection between $\{\psi_1=\pi\}$ and $A_2$ are
        \[
            (\theta_1,\psi_1)= \left(\pi \Delta(\lambda_{M_2},\mu_1) +\frac{2\pi k}{o_2} , \pi  \right).
        \]
        If $o_2 \le 2$, then $H_{1,\pi} \cap A_2= \emptyset$ since $H_{1,\pi}$ does not contain the points $(0,\pi)$ and $(\pi,\pi)$ by of Corollary \ref{lemma: H1 doesn't contain -2 and 2}.
        We conclude that, if $o_2\le 2$, then $A_2 \cap H_1 = \emptyset$.
        
        As a result of Lemma \ref{lemma: the interval Jpi(p1,p2) contains 0 if and only if p1,p2 is not 4,4}, if $o_2=4$, then $H_{1,\pi}$ is disjointed from $A_2$ if and only if $p_1=4$ and $p_2=4$. Thus, if $o_2 =4$, then either $H_{1,0}\cap A_2 \neq \emptyset$ or $H_{1,\pi}\cap A_2 \neq \emptyset$. This implies that if $o_2=4$, then $H_1 \cap A_2 \neq \emptyset$.

        \textbf{Case \romannum{3})} $\boldsymbol{\Delta(\lambda_{M_2},h_1)=2}$. 
        In this case $\Delta(\lambda_{M_2},\mu_1)\equiv_21$ and $o_2\le 2$. The intersection points of $A_2$ with $\{\psi_1=0\}$ are
        \[
            (\theta_1,\psi_1)= \left(\frac{\pi k}{o_2} ,0  \right) \quad \text{with} \quad k \in \{1, \cdots, 2o_2\}.
        \]
        According to Corollary \ref{lemma: H1 doesn't contain -2 and 2}, if $o_2=1$ then $H_{1,0}\cap A_2$ is empty. 
        
        As we said before, if $p_1 \ge 3$, then at least one of the points $(\nicefrac{\pi}{2},0)$ and $(\nicefrac{3\pi}{2},0)$ is in $H_{1,0}$.
        Hence, if $o_2=2$, then $H_{1,0}\cap A_2$ is empty if and only if $p_1=2$.
        
        Let us focus on the line $\{\psi_1=\pi\}$. We see that the points in $A_2\cap \{\psi_1=\pi\}$ are
        \[
            (\theta_1,\psi_1)= \begin{dcases}
                \left(\frac{\pi}{2} + k \pi,\pi  \right)_{k =0,1} &\text{if }o_2=1,\\
                \left(\frac{\pi k}{2},\pi  \right)_{k=1,\cdots,4} & \text{if }o_2=2.
            \end{dcases}
        \]
        
        If $o_2=1$, then $H_{1,\pi}\cap A_2$ is empty if only if $p_1=p_2=4$, as a result of Lemma \ref{lemma: surjective maps} and Lemma \ref{lemma: the interval Jpi(p1,p2) contains 0 if and only if p1,p2 is not 4,4}. 
        If $o_2=2$, then, since we suppose that $p_1=2$, we have that $H_{1,\pi} \cap A_2$ is nonempty by Lemma \ref{lemma: the interval Jpi(p1,p2) contains 0 if and only if p1,p2 is not 4,4}.
        We conclude that $A_2\cap H_1$ is empty if and only if $p_1=p_2=4$ and $o_2=1$.
        
        \textbf{Case \romannum{4})} $\boldsymbol{\Delta(\lambda_{M_2},h_1)=4}$.
        As before, $\Delta(\lambda_{M_2},\mu_1)\equiv_2 1$ and $o_2=1$.
        The intersection points of $A_2$ with $\{\psi_1=0\}$ are
        \[
            (\theta_1,\psi_1)= \left(\frac{\pi}{2}k , 0  \right) \quad \text{with} \quad k \in \{1, \cdots, 4\}.
        \]
        Lemma \ref{lemma: the interval I(S,S) has angle more than 2pi/3} implies that if $p_1 \ge 3$, then $0 \in J_0(p_1,p_2)$. Furthermore, Lemma \ref{lemma: surjective maps} implies that either $(\nicefrac{\pi}{2},0)$ or $(3\nicefrac{\pi}{2},0)$ is in $H_{1,0}$. Thus, $H_{1,0}\cap A_2 = \emptyset$ if and only if $H_{1,0}=\emptyset$. This latter happens if and only if $p_1=2$.

        Let us suppose that $p_1=2$.
        The intersection points of $A_2$ with $\{\psi_1=\pi\}$ are
        \[
            (\theta_1,\psi_1)= \left(\frac{\pi}{4}+\frac{\pi k}{2} , \pi  \right) \quad \text{with} \quad k \in \{1, \cdots, 4\}.
        \]
        If $p_2=4$, then $J_\pi(p_1,p_2)=J_\pi(2,4)=(-\sqrt{2},\sqrt{2})$. Lemma \ref{lemma: surjective maps} implies that $H_{1,\pi} \cap A_2 = \emptyset$. Conversely, if $p_2 \neq 4$ Corollary \ref{cor: l'angolo di Jpi(2,p2) è sempre maggione di 1/2pi} implies that there exists an $\varepsilon >0$ such that 
        \[
        (-\sqrt{2}-\varepsilon,\sqrt{2}+\varepsilon) \subseteq J_\pi(2,p_2)
        \]
        Thus, by Lemma \ref{lemma: surjective maps} $H_{1,\pi}\cap A_2 \neq \emptyset$.
        This leads us to the conclusion that if $\Delta(\lambda_{M_1},h_1)=4$, then $A_2 \cap H_1$ is empty if and only if $p_1=2$, $p_2=4$, and $o_2=1$.

       \textbf{Case \romannum{5})} $\boldsymbol{o_2\Delta(\lambda_{M_2},h_1)=3}$.
        This case includes ones for which either $\Delta(\lambda_{M_2},h_1)=3$ and $o_2=1$ or $\Delta(\lambda_{M_2},h_1)=1$ and $o_2=3$. In particular, $o_2 \equiv_21$. We need to prove that $H_1 \cap A_2$ if and only if $p_1=p_2=3$ and $\alg{\lambda_{M_2}}{\lambda_{M_1}} \equiv_2 0$.
        
        The intersections between $\{\psi_1=0\}$ and $A_2$ are
        \[
            (\theta_1,\psi_1)= \left(\frac{2\pi k}{3}, 0 \right) \quad \text{with} \quad k \in \{1, 2, 3\}.
        \]
        If $p_1 =2$, then $H_{1,0} = \emptyset$ and this implies that $H_{1,0} \cap L= \emptyset$. If $3 \le p_1$, then Lemma \ref{lemma: the interval I(S,S) has angle more than 2pi/3} and Lemma \ref{lemma: surjective maps} imply that $H_{1,0}$ is disjoint from $A_2$ if and only if
        \[
            J_0(p_1,p_2)= \left( 2\cos\left( \frac{2\pi}{3}\right),2\right)= J_0(3,3).
        \]
        Lemma \ref{lemma: computation of J_0(p1,p2)} implies that $J_0(p_1,p_2)=J_0(3,3)$ if and only if $p_1=p_2=3$. Thus, $H_{1,0}\cap A_2 = \emptyset$ if and only if either $p_1=2$ or $p_1=p_2=3$.

        The intersection of $A_2$ with $\{\psi_1=\pi\}$ contains either 
        \[
        (\theta_1,\psi_1)=\left( \frac{2\pi}{3},\pi\middle) \quad \text{or} \quad (\theta_1,\psi_1)=\middle(\frac{\pi}{3},\pi\right), 
        \]
        depending if the integer $n$ is even or odd, where $n= \pm \Delta(\lambda_{M_2},h_1)$. If $p_1=2$, then Corollary \ref{cor: l'angolo di Jpi(2,p2) è sempre maggione di 1/2pi} states that for every $p_2 \ge 2$, the interval $J_\pi(2,p_2)$ contains $J_\pi(2,4) = (-\sqrt{2},\sqrt{2})$. This implies that $2\cos\left( \nicefrac{2\pi}{3}\right) \in J_\pi (2,p_2)$ and $2\cos\left(     \nicefrac{\pi}{3}\right) \in J_\pi(2,p_2)$.
        As a result of Lemma \ref{lemma: surjective maps}, if $p_1=2$, then $A_2 \cap H_{1,\pi}=A_2 \cap H_1$ is nonempty.
        
        Let us now assume that $p_1=p_2=3$. Lemma \ref{lemma: description of the rational longitude} implies that $o_1 \equiv_2 1$. Since
        $
            J_\pi(p_1,p_2)= \left( 2\cos \left(\frac{2\pi}{3}\right),2 \right),
       $
        Lemma \ref{lemma: surjective maps} implies that
        \begin{equation}
            H_{1,\pi} = \left\{ (\theta_1,\pi)  \middle| \theta_1 \in  \left( 0,\frac{2\pi}{3} \right) \cup \left(\frac{4\pi}{3},2\pi\right) \right\} \subset \R{\partial M_1}.
            \label{eq: lo uso una volta: Jpi(3,3)}
        \end{equation}

        Lemma \ref{lemma: description of the rational longitude} implies that $o_1$ divides $\frac{p_2q_1+p_1q_2}{g_1}$. Since we supposed that $p_1=p_2=3$, the quantity $\frac{q_1+q_2}{o_1}$ is an integer. Furthermore, $\frac{q_1+q_2}{o_1}$ is even as we suppose that $q_1$ and $q_2$ are odd.
        Let $\eta \colon \fund{\partial M_1}\to SU(2)$ so that $\eta(h_1)=-1$ and $\eta(\mu_1)=1$.
        Lemma \ref{lemma: description of the rational longitude} implies that
        \[
            \eta(\lambda_{M_1})=\cancel{\eta(\mu_1)^{3}}\eta(h_1)^{\frac{q_1+q_2}{o_1}}= (-1)^{\frac{q_1+q_2}{o_1}}=1.
        \]
        The representation $\eta$ has $(\theta_1,\psi_1)$-coordinates equal to $(0,\pi)$. The \eqref{eq: lo uso una volta: Jpi(3,3)} implies that $\eta$ is in $\overline{H_{1,\pi}} \setminus H_{1,\pi}$.
        We note that the set $\overline{H_{1,\pi}} \setminus H_{1,\pi}$ consists in three points, see for instance Figure \ref{Some examples of toroidal representation spaces}.C. Since the intersections between $\{\psi_1=\pi\}$ and $A_2$ are three equally spaced point, $H_{1,\pi}\cap A_2$ is empty if and only if $\{\psi_1=\pi\} \cap A_2 = \overline{H_{1,\pi}} \setminus H_{1,\pi}$.
        This latter happens if and only if $\eta \in A_2$. As a result of Proposition \ref{prop: a representation from the boundary to the whole manifold extends iff is trivial on the kernel of iota}, $\eta \in A_2$ if and only if $\eta(\lambda_{M_2})^{o_2}=1$.
        Since $o_2$ is odd by assumption and $\eta$ is a central representation, this is equivalent to say that $\eta(\lambda_{M_2})=1$.
        Let $\xi$ be a simple closed curve in $\partial M_1$ so that $\{\lambda_{M_1},\xi\}$ is a basis of $\fund{\partial M_1}$, then
        \[
            \lambda_{M_2}= a \xi+ b  \lambda_{M_1} \in \fund{\partial M_2},
        \]
        where $|a|= \Delta(\lambda_{M_2},\lambda_{M_1})$ and $|b|= \alg{\lambda_{M_2}}{\xi}$.
        This brings us to the following:
        \[
            \eta(\lambda_{M_2})^{o_2}= \eta(\lambda_{M_2})= \eta(\xi)^a \cancel{\eta(\lambda_{M_1})^b} = (\pm 1) ^{\alg{\lambda_{M_2}}{\lambda_{M_1}}}.
        \]
        Thus, $\eta(\lambda_{M_2})^{o_2}=1$ if and only if $\alg{\lambda_{M_2}}{\lambda_{M_1}} \equiv_2 0$.
        Summing up, if $o_2\Delta(\lambda_{M_2},h_1)=3$, then $H_1 \cap A_2$ is empty if and only if $p_1=p_2=3$ and $\Delta(\lambda_{M_1},\lambda_{M_2})$ is even.
    \end{proof}
\end{prop}

The following is the symmetric version of Proposition \ref{prop: H1 and A2 sono disjoint if and only if}.

\begin{prop}
    The intersection $H_2 \cap A_1$ if empty if and only if one of the following holds:
    \begin{itemize}
        \item $\Delta(\lambda_{M_1},h_2)=0$, $p_3=2$ and $o_1\equiv_21$;
        \item $\Delta(\lambda_{M_1},h_2)=1$ and $o_1\le 2$,
        \item $\Delta(\lambda_{M_1},h_2)=1$, $p_3=p_4=3$, $o_1=3$, and $\Delta(\lambda_{M_2},\lambda_{M_1})\equiv_2 0$,
        \item $\Delta(\lambda_{M_1},h_2)=2$, $p_3=p_4=4$ and $o_1=1$,
        \item $\Delta(\lambda_{M_1},h_2)=3$, $o_1=1$, $p_3=p_4=3$, and $\Delta(\lambda_{M_2},\lambda_{M_1})\equiv_2 0$,
        \item $\Delta(\lambda_{M_1},h_2)=4$, $p_3=2$, $p_4=4$ and $o_1=1$.
    \end{itemize}
    \label{prop: H2 and A1 sono disjoint if and only if}
\end{prop}

\section{The intersection between \texorpdfstring{$H_1$}{H\_1} and \texorpdfstring{$H_2$}{H\_2}}
\label{sec: H1 cap H2}
In this section we shall prove that $H_1 \cap H_2 = \emptyset$ if and only if $|\beta|\coloneqq \Delta(h_1,h_2) \le 2$ up to some exceptions. 
\begin{lemma}
    Let $M$ be the manifold $M_1 \cup_{\varphi} M_2$. Let us suppose than neither $(p_1,p_2,|\beta|)$ nor $(p_3,p_4,|\beta|)$ are in the set $\{(2,4,4),(3,3,3)\}$.
    The intersection $H_1 \cap H_2$ is empty if and only if one of the following condition holds:
    \begin{itemize}
        \item $|\beta|=1$,
        \item $|\beta|=2$ and $p_1=p_3=2$,
        \item $|\beta|=2$ and $p_1=p_2=p_3=p_4=4$.
    \end{itemize}
    \label{lemma: H1 cap H2 empty if and only if beta 1 or 2}
    \begin{proof}
        The conclusion holds by Lemma \ref{lemma: beta =1 then H1 intersection H2 is empty}, Lemma \ref{lemma: beta = 2 then H1 intersection H2} and Lemma \ref{lemma: if b >= 3 then H1 and H2 contain an intersection}.
    \end{proof} 
\end{lemma}

Let us fix the presentations for $\fund{M_1}$ and $\fund{M_2}$ as in \eqref{eq: presentation of the fundamental group of M1 and M2}.
We choose $\{\mu_1,h_1\}$ and $\{\mu_2,h_2\}$ as bases for the groups $\fund{\partial M_1}$ and $\fund{\partial M_2}$.
We recall that, with respect to the chosen bases for $\fund{\partial M_1}$ and $\fund{\partial M_2}$, the map $\varphi^\ast \colon \fund{\partial M_1}\to \fund{\partial M_2}$ is represented by the matrix
\begin{align}
    \begin{bmatrix}
     \alpha & \beta \\ \gamma & \delta
    \end{bmatrix}  \quad \text{with} \quad \alpha\delta - \beta\gamma = -1.
    \label{eq: una volta la mappa per le entrate}
\end{align}
In particular $|\beta|=\Delta(h_1,h_2)$. For $i \in \{1,2\}$, we are going to fix the coordinates for the torus $\R (\partial M_i)$ in the following way: the point $(\theta_i,\psi_i)\in [0,2\pi]^2/_\sim$ corresponds to the representation in $\R(\partial M_i)$ such that
\[
\mu_i \mapsto \begin{bmatrix}
    e^{i\theta_i} & 0 \\ 0 & e^{-i\theta_i}
\end{bmatrix} \quad \text{and} \quad h_i \mapsto \begin{bmatrix}
    e^{i\psi_i} & 0 \\ 0 & e^{-i\psi_i}
\end{bmatrix}.
\]
We recall that, since in $M$ the boundary $\partial M_1$ coincides with $\partial M_2$, the space $\R(\partial M_1)$ coincides with $\R(\partial M_2)$. Consequently, $(\theta_1,\psi_1)$ and $(\theta_2,\psi_2)$ are two ways of parameterizing the space $\R(\Sigma)$, where $\Sigma$ is the torus corresponding to $\partial M_1 =\partial M_2$.

Let us define the the lines $L_0$ and $L_\pi$ of $\R (\partial M_1)$ as
\begin{multline*}
L_0 \coloneqq \left\{ (\theta_1,\psi_1)=(\alpha \theta_2, \beta \theta_2) \; \middle| \; \theta_2 \in [0,2\pi]\right\}\subset \R (\partial M_1) \quad \text{and} \\ \quad L_\pi \coloneqq \left\{ (\theta_1,\psi_1)=(\alpha \theta_2 + \pi \gamma, \beta \theta_2+ \pi \delta) \; \middle| \; \theta_2 \in [0,2\pi]\right\} \subset \R (\partial M_1).
\end{multline*}
Here the quantities $\alpha, \beta, \gamma$ and $\delta$ are as in \eqref{eq: una volta la mappa per le entrate}.

Let $\rho \colon \fund{M} \to SU(2)$ be a representation and let $\rho_1$ and $\rho_2$ the restrictions on $\fund{M_1}$ and $\fund{M_2}$ respectively. The group presentation \eqref{fund(M) descrizione totale} of $\fund{M}$ gives that
\begin{equation}
\rho_2(a_2b_2)^\alpha\rho_2(h_2)^\gamma=\rho_1(a_1b_1) \quad \text{and} \quad \rho_2(a_2b_2)^\beta\rho_2(h_2)^\delta=\rho_1(h_1).
\label{eq: conditizione per incollare due rappresentazioni dai due lati}
\end{equation}
Therefore, $L_0$ is the line $\{\psi_2=0\}\subset \R(\partial M_2)$ in $(\theta_1,\psi_1)$ coordinates. Similarly, $L_\pi$ is the line $\{\psi_2=\pi\}\subset \R(\partial M_2)$ in $(\theta_1,\psi_1)$ coordinates. 
Hence, $L_0$ and $L_\pi$ respectively contain the sets $H_{2,0}$ and $H_{2,\pi}$ in the torus $\R (\partial M_1)$. Thus,
\[
    H_1 \cap H_2  \subseteq \left(\{\psi_1 = 0\}\cup \{\psi_1=\pi\}\right) \cap (L_0 \cup L_\pi). 
\]
We define $4P\coloneq\{0,\pi\}^2 \subset \R (\partial M_1)$.
Lemma \ref{lemma: H1 doesn't contain -2 and 2} states that $H_1 \cap 4P = \emptyset$.

\begin{lemma}
    If $|\beta|=1$, then $H_1 \cap H_2$ is empty.
    \label{lemma: beta =1 then H1 intersection H2 is empty}
    \begin{proof}
       The intersection of $L_0 \cup L_\pi$ with $\{\psi_1 =0 \} \cup \{\psi_1=\pi\}$ is contained in $4P$. Since $H_1$ is disjoint from $4P$ by Lemma \ref{lemma: H1 doesn't contain -2 and 2}, we get the conclusion. 
    \end{proof}
\end{lemma}

\begin{lemma}
    Let us suppose that $|\beta|=2$. The set $H_1 \cap H_2$ is empty if and only either $p_1=p_3=2$ or $p_1=p_2=p_3=p_4=4$.
    \label{lemma: beta = 2 then H1 intersection H2}
    \begin{proof}
        Since $\gcd(\alpha,\beta)=1$, the integer $\alpha$ is odd.
        We note that $\{\psi_1=0\} \cap L_0$ and $\{\psi_1=\pi\} \cap L_\pi$ are in $4P$. By Lemma \ref{lemma: H1 doesn't contain -2 and 2}, $H_{1,0}\cap H_{2,0}=H_{1,\pi}\cap H_{2,\pi}=\emptyset$. The intersection $ \{\psi_1=\pi\} \cap L_0$ is equal, in $(\theta_1,\psi_1)$ and $(\theta_2,\psi_2)$ coordinates, to
        \[
            \left(\theta_1,\psi_1\right) \in \left\{ \left(\frac{\pi}{2},\pi\right),\left(\frac{3\pi}{2},\pi\right)\right\} \quad \text{and} \quad \left(\theta_2,\psi_2\right) \in \left\{ \left(\frac{\pi}{2},0\right),\left(\frac{3\pi}{2},0\right)\right\}.
        \]
        According to Lemma \ref{lemma: surjective maps}, the points
        \begin{equation}
            \left\{ \left(\frac{\pi}{2},\pi\right),\left(\frac{3\pi}{2},\pi \right) \right\} \subset \R(\partial M_1)
            \label{eq: una volta beta=2 M1}
        \end{equation}
        are in $H_{1,\pi}$
        if and only if $0 \in J_\pi(p_1,p_2)$. Similarly, Lemma \ref{lemma: surjective maps} states that the points
        \begin{equation}
            \left\{ \left(\frac{\pi}{2},0\right),\left(\frac{3\pi}{2},0 \right) \right\} \subset \R (\partial M_2)
            \label{eq: una volta beta=2 M2}
        \end{equation}
        are in $H_{2,0}$ and only if $0 \in J_0(p_3,p_4)$.
        Thus, $H_{1,\pi} \cap H_{2,0} \neq \emptyset$ if and only if $0 \in J_\pi(p_1,p_2)$ and $0 \in J_0(p_3,p_4)$.
        
        As a consequence of Lemma \ref{lemma: the interval I(S,S) has angle more than 2pi/3}, if $J_0(p_1,p_2)$ is nonempty, then it contains $0$.
        We conclude, by Lemma \ref{lemma: the interval Jpi(p1,p2) contains 0 if and only if p1,p2 is not 4,4}, that $H_{1,\pi} \cap H_{2,0}$ is empty if and only if either $p_1=p_2=4$ or $p_3=2$. By a similar analysis, we can conclude that of $H_{2,\pi} \cap H_{1,0}$ is empty if and only if either $p_1=2$ or $p_3=p_4=4$.
        Hence, we have that $H_{1,\pi} \cap H_{2,0} = \emptyset$ and $H_{1,0} \cap H_{2,\pi}=\emptyset$ if and only if either $p_1=p_3=2$ or $p_1=p_2=p_3=p_4=4$.
    \end{proof}
    \end{lemma}
    
    The next lemma uses four claims that are proved in next subsection.
    
    \begin{lemma}
        If $|\beta|\ge 3$ and $(p_1,p_2,|\beta|),(p_3,p_4,|\beta|)\notin \{(2,4,4),(3,3,3)\}$, then $H_1 \cap H_2 \neq \emptyset$.
        \label{lemma: if b >= 3 then H1 and H2 contain an intersection}
        \begin{proof}
            With an abuse of notation, we will write $\beta$ for its absolute value.
            The proof is divided in $4$ cases:
            \begin{enumerate}[label=\roman*)]
                \item $p_1 \ge 3$ and $p_3 \ge 3$,
                \item $p_1 \ge 3$ and $p_3 =2$,
                \item $p_1=2$ and $p_3 \ge 3$,
                \item $p_1=2$ and $p_2=2$.
            \end{enumerate}
             \textbf{Case \romannum{1})} $\boldsymbol{p_1 \ge 3}$ \textbf{and} $\boldsymbol{p_3 \ge 3}$. We will prove that $H_{1,0} \cap H_{2,0} \neq \emptyset$. We define the sets $\mathcal{S}_1$ and $\mathcal{S}_2$ as follows:
            \begin{equation}
                \begin{split}
                    \mathcal{S}_1 \coloneqq \left\{ 2 \cos\left( \frac{2\pi k}{\beta} \right) \middle| k \in \{1,\cdots, \beta\}, 2 \cos\left( \frac{2\pi k}{\beta} \right) \in J_0(p_1,p_2) \right\} \quad \text{and} \\[3pt] \mathcal{S}_2 \coloneqq \left\{ 2 \cos\left( \frac{2\pi \alpha k }{\beta} \right) \middle| k \in \{1,\cdots, \beta\}, 2 \cos\left( \frac{2\pi k }{\beta} \right) \in J_0(p_3,p_4)  \right\}.
                \end{split}
                \label{eq: definizione S1 e S2 nella dimostrazione H1capH2}
            \end{equation}
            The two sets in \eqref{eq: nuova definizione S1 e S2} are both subsets of
            $
                \left\{ 2 \cos\left( \frac{2\pi k}{\beta} \right) \right\}_{k \in \{1,\cdots, \beta\}}.
            $
            Since the intervals $J_0(p_1,p_2)$ and $J_0(p_3,p_4)$ do not contain $\pm 2$, the sets $\mathcal{S}_1$ and $\mathcal{S}_2$ do not contain $\pm 2$ by Lemma \ref{lemma: surjective maps}.
            Hence,
            \[
                |\mathcal{S}_1|\le \frac{\beta}{2}-1 \quad \text{and} \quad |\mathcal{S}_2|\le \frac{\beta}{2}-1.
            \]
            According to Lemma \ref{lemma: the interval I(S,S) has angle more than 2pi/3}, $\mathcal{S}_1$ (resp. $\mathcal{S}_2$) is empty if and only if $J_0(p_1,p_2)=J_0(3,3)$ (resp. $J_0(p_3,p_4)=J_0(3,3)$) and $\beta =3$.
            Lemma \ref{lemma: computation of J_0(p1,p2)} implies that $\mathcal{S}_1$ (resp. $\mathcal{S}_2$) is empty if and only if $p_1=p_2=\beta=3$ (resp. $p_3=p_4=\beta=3$). Thus neither $\mathcal{S}_1$ nor $\mathcal{S}_2$ is empty.
            \begin{claim}\label{claim: S1 e S2 danno H10 e H20 empty}
                Let $\mathcal{S}_1$ and $\mathcal{S}_2$ be as in \eqref{eq: definizione S1 e S2 nella dimostrazione H1capH2}. Then
                $H_{1,0} \cap H_{2,0}$ is empty if and only if $\mathcal{S}_1 \cap \mathcal{S}_2$ is empty. 
                \end{claim}
                \begin{claim}\label{claim: gli angli sono troppo grandi allora H10 H2}
             Let $\mathcal{S}_1$ and $\mathcal{S}_2$ be as in \eqref{eq: definizione S1 e S2 nella dimostrazione H1capH2}. If either $2\beta < p_2$ or $2\beta < p_4$, then $\mathcal{S}_1 \cap \mathcal{S}_2 \neq \emptyset$.
            \end{claim}
            Let $\alpha_1$ (resp. $\alpha_2$) be the angle supporting $J_0(p_1,p_2)$ (resp. $J_0(p_3,p_4)$), Corollary \ref{cor: angle supporting J0(p1,p2) is bigger that pi-4pi/p2} implies that
            $\alpha_1 \ge \pi - \frac{4\pi}{p_2}$ and $\alpha_2\ge \pi - \frac{4\pi}{p_4}$.
            If either $2\beta < p_2$ or $2\beta < p_4$, then the conclusion holds by Claim \ref{claim: S1 e S2 danno H10 e H20 empty} and Claim \ref{claim: gli angli sono troppo grandi allora H10 H2}.
            Thus, we assume $p_2 \le 2\beta$ and $p_4 \le 2\beta$.
            Lemma \ref{lemma: the interval I(S,S) has angle more than 2pi/3} implies that $\alpha_1 \ge \nicefrac{2\pi}{3}$ and $\alpha_2 \ge \nicefrac{2\pi}{3}$. 
            Moreover, Corollary \ref{cor: angle supporting J0(p1,p2) is bigger that pi-4pi/p2} implies that if $\alpha_1$ and $\alpha_2$ are both smaller than $\nicefrac{4\pi}{5}$, then $ p_2\le 20$ and $p_4 \le 20$. Hence, up to checking by hand the conclusion for those finite many cases, we can suppose $\alpha_1 \ge \nicefrac{2\pi}{3}$ and $\alpha_2 \ge \nicefrac{4\pi}{5}$. 
            Therefore, we have that
            \[
            \left| \mathcal{S}_1\right| \ge  \left\lfloor \frac{\beta}{3} \right\rfloor-2 \quad \text{and} \quad \left| \mathcal{S}_2\right| \ge  \left\lfloor \frac{2\beta}{5} \right\rfloor-2.
            \]
            Thus, if $\mathcal{S}_1 \cap \mathcal{S}_2 = \emptyset$, then $|\mathcal{S}_1|+ |\mathcal{S}_2| \le \frac{\beta}{2}-1$.
            Then we get the following:
            \begin{equation*}
                \begin{aligned}
               \left(\frac{\beta}{3} - \frac{2}{3} \right)+\left(\frac{2\beta}{5} - \frac{4}{5} \right)-4 \le  \left\lfloor \frac{\beta}{3} \right\rfloor + \left\lfloor \frac{2\beta}{5} \right\rfloor-4 \le \left| \mathcal{S}_1\right|+ \left| \mathcal{S}_2\right| - \cancel{\left| \mathcal{S}_1 \cap \mathcal{S}_2\right|} & \le \frac{\beta}{2}-1 \\
               \left(\frac{\beta}{3} - \frac{2}{3} \right)+\left(\frac{2\beta}{5} - \frac{4}{5} \right)-4 & \le \frac{\beta}{2}-1 \\
                \beta &\le  19.
                \end{aligned}
            \end{equation*}
            Up to checking finite many cases (i.e. the ones such that $p_1\le p_2 \le 38$ and $p_3 \le p_4 \le 38$) we conclude that if $p_1\ge3$ and $p_3 \ge 3$, then $H_{1,0}\cap H_{2,0} \neq \emptyset$. The author checked these finite many cases with Algorithm \ref{alg1}.
            \begin{algorithm}
            \caption{If $3\le \beta\le 19$, then $H_{1,0}\cap H_{2,0} \neq \emptyset$}
            \label{alg1}
                \begin{algorithmic}
                \For{$\beta \le 19$}
                    \For{$\alpha \le \beta \wedge \gcd(\alpha, \beta)=1$}
                        \For{$p_i \le 2\beta$}
                            \If{$(p_1,p_2,\beta)\neq (3,3,3) \wedge (p_3,p_4,\beta)\neq (3,3,3)$}
                                \State Compute $J_0(p_1,p_2)$ and $J_0(p_3,p_4)$ using Lemma \ref{lemma: computation of J_0(p1,p2)}
                                \State Compute $\mathcal{S}_1$ and $\mathcal{S}_2$ as in \eqref{eq: definizione S1 e S2 nella dimostrazione H1capH2}
                                \State Check $\mathcal{S}_1 \cap \mathcal{S}_2 \neq \emptyset$
                            \EndIf
                        \EndFor
                    \EndFor
                \EndFor
                \end{algorithmic}
            \end{algorithm}

             \textbf{Case \romannum{2})} $\boldsymbol{p_1 \ge 3}$ \textbf{and} $\boldsymbol{p_3 =2}$. We shall prove that $H_{1,0}\cap H_{2,\pi} \neq \emptyset$. 
            We define the sets $\mathcal{S}_3$ and $\mathcal{S}_4$ as follows:
            \begin{equation}
                \begin{split}
            \mathcal{S}_3 \coloneqq \left\{ 2 \cos\left( \frac{\pi k}{\beta} -\gamma \pi \right) \middle| k \in \{1,3, \cdots , 2\beta-1\}, 2 \cos\left( \frac{\pi k}{\beta} \right) \in J_0(p_1,p_2) \right\} 
            \quad \text{and} \\ \mathcal{S}_4 \coloneqq \left\{ 2 \cos\left( \frac{2\pi \alpha k }{\beta}-\frac{\delta \alpha\pi }{\beta} \right) \middle| k \in \{1,\cdots, \beta\}, 2 \cos\left( \frac{2\pi k }{\beta} - \frac{\delta \pi}{\beta
            } \right) \in J_\pi(2,p_4)  \right\}.
            \label{eq: nuova definizione S1 e S2}
            \end{split}
            \end{equation}
            Since $(p_1,p_2,|\beta|),(p_3,p_4,|\beta|)\notin \{(2,4,4),(3,3,3)\}$, neither $\mathcal{S}_3$ nor $\mathcal{S}_4$ is empty.
            The sets $\mathcal{S}_3$ and $\mathcal{S}_4$ are both subsets of either
            \[
            \left\{2\cos \left(\frac{\pi (2k+1)}{\beta} \right)\right\}_{k \in \mathbb{Z}} \quad \text{or} \quad \left\{2\cos \left(\frac{2\pi k}{\beta}\right)\right\}_{k \in \mathbb{Z}},
            \]
            depending if $\beta \gamma \equiv_20$ or $\beta \gamma \equiv_2 1$. In particular,
            \[
                |\mathcal{S}_3|\le \frac{\beta}{2}-1 \quad \text{and} \quad |\mathcal{S}_4|\le \frac{\beta}{2}-1.
            \]
            \begin{claim}
                Let $\mathcal{S}_3$ and $\mathcal{S}_4$ as in \eqref{eq: nuova definizione S1 e S2}. Then $\mathcal{S}_3 \cap \mathcal{S}_4$ is empty if and only if  $H_{1,0}\cap H_{2,\pi}$ is empty.
                \label{claim: S3 S4 iff H10 H2pi}
            \end{claim}

            \begin{claim}
                 Let $\mathcal{S}_3$ and $\mathcal{S}_4$ be as in \eqref{eq: nuova definizione S1 e S2}. If either $p_4 > \beta$ or $p_2 > 2\beta$, then $\mathcal{S}_3 \cap \mathcal{S}_4 \neq \emptyset$.
                 \label{claim: angoli grossi per S3 e S4}
            \end{claim}
            Without loss of generality, we can assume that $p_4 \le \beta$ and $p_2 \le 2\beta$. Let $\alpha_1$ and $\alpha_2$ be the angles supporting $J_0(p_1,p_2)$ and $J_\pi(2,p_4)$ respectively, then, according to Lemma \ref{lemma: the interval I(S,S) has angle more than 2pi/3} and Corollary \ref{cor: l'angolo di Jpi(2,p2) è sempre maggione di 1/2pi},
            $
                \alpha_1 \ge \nicefrac{2\pi}{3}$ and $ \alpha_2 \ge \nicefrac{\pi}{2}.
            $
            Lemma \ref{lemma: the interval Jpi(2,p2) is nested and connected} and Corollary \ref{cor: angle supporting J0(p1,p2) is bigger that pi-4pi/p2} imply that if $\alpha_2 < \nicefrac{4\pi}{5}$ and $\alpha_1 < \nicefrac{4\pi}{5}$ then, $p_2 \le 20$ and $p_4 \le 10$. Up to checking by hand the conclusion for these finite many cases, we can suppose that $\alpha_2 \ge \nicefrac{4\pi}{5}$.
            Therefore, we have that
            \[
            \left| \mathcal{S}_3\right| \ge  \left\lfloor \frac{\beta}{3} \right\rfloor-2 \quad \text{and} \quad \left| \mathcal{S}_4\right| \ge  \left\lfloor \frac{2\beta}{5} \right\rfloor-2.
            \]
            Thus, if $\mathcal{S}_3 \cap \mathcal{S}_4 = \emptyset$, then $|\mathcal{S}_3|+ |\mathcal{S}_4| \le \frac{\beta}{2}-1$.
            As we stated before, if $\mathcal{S}_3 \cap \mathcal{S}_4 = \emptyset$, then
             \begin{equation*}
                \begin{aligned}
               \left(\frac{2\beta}{5} - \frac{4}{5} \right) + \left(\frac{\beta}{3} - \frac{2}{3} \right)-4 \le  \left\lfloor \frac{2\beta}{5} \right\rfloor + \left\lfloor \frac{\beta}{3} \right\rfloor-4 \le \left| \mathcal{S}_1\right|+ \left| \mathcal{S}_2\right| - \cancel{\left| \mathcal{S}_3 \cap \mathcal{S}_4\right|} & \le \frac{\beta}{2} -1 \\
               \left(\frac{2\beta}{5} - \frac{4}{5} \right) + \left(\frac{\beta}{3} - \frac{2}{3} \right)-4 & \le \frac{\beta}{2} -1 \\
                \beta&\le  19.
                \end{aligned}
            \end{equation*}
            Again, up to checking by hand the remaining cases, namely the ones such that $p_1 \le p_2 \le 38$ and $ p_4 \le 19$, we get the conclusion. The author checked these finite many cases with Algorithm \ref{alg2}.
            \begin{algorithm}
            \caption{If $3\le \beta\le 19$, then $H_{1,0}\cap H_{2,\pi} \neq \emptyset$}
            \label{alg2}
                \begin{algorithmic}
                \For{$\beta \le 19$}
                    \For{$\alpha \le \beta \wedge \gcd(\alpha, \beta)=1$}
                        \For{$p_1,p_2 \le 2\beta \text{ and } p_4 \le \beta$}
                            \If{$(p_1,p_2,\beta)\neq (3,3,3) \wedge (p_3,p_4,\beta)\neq (2,4,4)$}
                                \State Compute $J_0(p_1,p_2)$ and $J_\pi(2,p_4)$ with Lemma \ref{lemma: computation of J_0(p1,p2)} and Lemma \ref{lemma: the interval Jpi(2,p2) is nested and connected}
                                \State Compute $\mathcal{S}_3$ and $\mathcal{S}_4$ as in \eqref{eq: nuova definizione S1 e S2}
                                \State Check $\mathcal{S}_3 \cap \mathcal{S}_4 \neq \emptyset$
                            \EndIf
                        \EndFor
                    \EndFor
                \EndFor
                \end{algorithmic}
            \end{algorithm}
    
             \textbf{Case \romannum{3})} $\boldsymbol{p_1 = 2}$ \textbf{and} $\boldsymbol{p_3  \ge 3}$.
            This case is proven as in the previous case withe the roles of $H_1$ and $H_2$ switched.

             \textbf{Case \romannum{4})} $\boldsymbol{p_1 = 2}$ \textbf{and} $\boldsymbol{p_3  =2}$.
            We shall prove that $H_{1,\pi}\cap H_{2,\pi} \neq \emptyset$. 
            We define the sets $\mathcal{S}_5$ and $\mathcal{S}_6$ as follows:
            \begin{equation}
                \begin{split}
            \mathcal{S}_5 \coloneqq \left\{ 2 \cos\left( \frac{2 \pi k}{\beta} + \frac{\pi (\alpha+1)}{\beta} - \pi \gamma \right) \middle| k \in \{1, \cdots , \beta\}, 2 \cos\left( \frac{2 \pi k}{\beta} + \frac{\pi (\alpha+1)}{\beta} \right) \in J_\pi(2,p_2) \right\} 
            \quad \text{and} \\ \mathcal{S}_6 \coloneqq \left\{ 2 \cos\left( \frac{2 \pi \alpha k}{\beta} + \frac{\pi \alpha(\delta+1)}{\beta} \right) \middle| k \in \{1,\cdots, \beta\}, 2 \cos\left( \frac{2 \pi k}{\beta} + \frac{\pi (\delta+1)}{\beta} \right) \in J_\pi(2,p_4)  \right\}.
            \label{eq: nuova definizione S5 e S6}
            \end{split}
            \end{equation}
            Since $(p_1,p_2,|\beta|) \neq (2,4,4)$ and $(p_3,p_4,|\beta|) \neq (2,4,4)$ by hypothesis, neither $\mathcal{S}_5$ nor $\mathcal{S}_6$ is empty.
            The sets $\mathcal{S}_5$ and $\mathcal{S}_6$ are both subsets of either
            \[
            \left\{2\cos \left(\frac{\pi (2k+1)}{\beta} \right)\right\}_{k \in \mathbb{Z}} \quad \text{or} \quad \left\{2\cos \left(\frac{2\pi k}{\beta}\right)\right\}_{k \in \mathbb{Z}},
            \]
            if either $\alpha\equiv_21$ and $\delta \equiv_2 0$ or $\beta \gamma \alpha \equiv_2 0$ respectively. In particular,
            \[
                |\mathcal{S}_5|\le \frac{\beta}{2}-1 \quad \text{and} \quad |\mathcal{S}_6|\le \frac{\beta}{2}-1.
            \]
            \begin{claim}
                Let $\mathcal{S}_5$ and $\mathcal{S}_6$ as in \eqref{eq: nuova definizione S5 e S6}. Then $\mathcal{S}_5 \cap \mathcal{S}_6$ is empty if and only if  $H_{1,\pi}\cap H_{2,\pi}$ is empty.
                \label{claim: S5 S6 iff H1pi H2pi}
            \end{claim}
            \begin{claim}
                 Let $\mathcal{S}_5$ and $\mathcal{S}_6$ be as in \eqref{eq: nuova definizione S5 e S6}. If either $p_4 > \beta$ or $p_2 >\beta$, then $\mathcal{S}_3 \cap \mathcal{S}_4 \neq \emptyset$.
                 \label{claim: angoli grossi per S5 e S6}
            \end{claim}
            Without loss of generality, we can assume that $p_4 \le \beta$ and $p_2 \le 2\beta$. Let $\alpha_1$ and $\alpha_2$ be the angles supporting $J_\pi(2,p_2)$ and $J_\pi(2,p_4)$ respectively. According to Corollary \ref{cor: l'angolo di Jpi(2,p2) è sempre maggione di 1/2pi},
            $
                \alpha_1 \ge \frac{\pi}{2}$ and $\alpha_2 \ge \frac{\pi}{2}.
            $
            Corollary \ref{cor: angle supporting J0(p1,p2) is bigger that pi-4pi/p2} imply that if $\alpha_2 < \nicefrac{4\pi}{5}$ and $\alpha_1 < \nicefrac{4\pi}{5}$ then, $p_2 \le 10$ and $p_4 \le 10$. Up to checking by hand the 
            for these finite many cases, we can suppose that $\alpha_2 \ge \nicefrac{4\pi}{5}$.
            Thus,
            \[
            \left| \mathcal{S}_5\right| \ge  \left\lfloor \frac{\beta}{4} \right\rfloor-2 \quad \text{and} \quad \left| \mathcal{S}_6\right| \ge  \left\lfloor \frac{2\beta}{5} \right\rfloor-2.
            \]
            Therefore, if $\mathcal{S}_5 \cap \mathcal{S}_6 = \emptyset$, then $|\mathcal{S}_1|+ |\mathcal{S}_2| \le \frac{\beta}{2}-1$.
            If $\mathcal{S}_5 \cap \mathcal{S}_6 = \emptyset$, then
             \begin{equation*}
                \begin{aligned}
               \left(\frac{\beta}{4} - \frac{3}{4} \right)+ \left(\frac{2\beta}{5} - \frac{4}{5} \right)-4 \le  \left\lfloor \frac{\beta}{4} \right\rfloor + \left\lfloor \frac{2\beta}{5} \right\rfloor -4 \le \left| \mathcal{S}_1\right|+ \left| \mathcal{S}_2\right| - \cancel{\left| \mathcal{S}_3 \cap \mathcal{S}_4\right|} & \le \frac{\beta}{2} -1 \\
               \left(\frac{\beta}{4} - \frac{3}{4} \right)+ \left(\frac{2\beta}{5} - \frac{4}{5} \right)-4 & \le \frac{\beta}{2} -1 \\
                \beta&\le  30.
                \end{aligned}
            \end{equation*}
            Again, up to checking by hand the remaining cases, namely the ones such that $ p_2 \le 30$ and $ p_4 \le 30$, we get the conclusion. The author checked these finite many cases with Algorithm \ref{alg3}.
            \begin{algorithm}
            \caption{If $3\le \beta\le 30$, then $H_{1,\pi}\cap H_{2,\pi} \neq \emptyset$}
            \label{alg3}
                \begin{algorithmic}
                \For{$\beta \le 30$}
                    \For{$\alpha \le \beta \wedge \gcd(\alpha, \beta)=1$}
                        \For{$p_2 \le \beta \text{ and } p_4 \le \beta$}
                            \If{$(p_1,p_2,\beta)\neq (2,4,4) \wedge (p_3,p_4,\beta)\neq (2,4,4)$}
                                \State Compute $J_\pi(2,p_2)$ and $J_\pi(2,p_4)$ Lemma \ref{lemma: the interval Jpi(2,p2) is nested and connected}
                                \State Compute $\mathcal{S}_5$ and $\mathcal{S}_6$ as in \eqref{eq: nuova definizione S1 e S2}
                                \State Check $\mathcal{S}_5 \cap \mathcal{S}_6 \neq \emptyset$
                            \EndIf
                        \EndFor
                    \EndFor
                \EndFor
                \end{algorithmic}
            \end{algorithm}
        \end{proof}
    \end{lemma}
\subsection{Proof of the claims}

In this section we prove the claims we used in Lemma \ref{lemma: if b >= 3 then H1 and H2 contain an intersection}.

            \begin{proof}[Proof of Claim \ref{claim: S1 e S2 danno H10 e H20 empty}]
            Let us suppose that $\mathcal{S}_1 \cap \mathcal{S}_2$ is nonempty. Then, there exist two integers $k_1, k_2 \in \{1,\cdots, \beta\}$ such that
           \begin{equation}
                \begin{split}
                2 \cos\left( \frac{2\pi k_1}{\beta} \right) \in J_0(p_1,p_2), \quad 2 \cos\left( \frac{2\pi k_2 }{\beta} \right) \in J_0(p_3,p_4), \quad \text{and} \quad 2 \cos\left( \frac{2\pi k_1}{\beta} \right)=2 \cos\left( \frac{2\pi \alpha k_2}{\beta} \right).
                \end{split}
                \label{eq: la uso una volta: coseni}
            \end{equation}
            By Lemma \ref{lemma: surjective maps}, there exist two irreducible representations $\rho_1 \in \mathcal{R}(M_1)$ and $\rho_2 \in \mathcal{R}(M_2)$ such that $\rho_1(h_1)=\rho_2(h_2)=1$,
            \[
                \Tr \rho_1(a_1b_1)= 2 \cos\left( \frac{2\pi k_1}{\beta} \right),  \quad \text{and} \quad \Tr \rho_2(a_2b_2)=  2 \cos\left( \frac{2\pi k_2 }{\beta} \right).
            \]
            The third condition in \eqref{eq: la uso una volta: coseni} implies that $\rho_1(a_1b_1)=\rho_2(a_2b_2)^\alpha$ up to conjugation. Thus, the representations $\rho_1$ and $\rho_2$ satisfy the conditions in \eqref{eq: conditizione per incollare due rappresentazioni dai due lati}. This guarantees the existence of an irreducible representation $\rho \colon \fund{M} \to SU(2)$ such that $\restr{\rho}{\fund{M_1}}=\rho_1$ and $\restr{\rho}{\fund{M_2}}=\rho_2$. This implies that the intersection $H_{1,0}\cap H_{2,0}$ contains a representation of $\fund{\partial M_1}=\fund{\partial M_2}$ that extends to $\fund{M}$. This implies that $H_{1,0}\cap H_{2,0}$ is nonempty.
            
            Conversely, if $H_{1,0}\cap H_{2,0}$ is nonempty, then there exists and irreducible representation $\rho \colon \fund{M}\to SU(2)$ such that $\restr{\rho}{\fund{M_1}}$ (resp. $\restr{\rho}{\fund{M_2}}$) is irreducible and $\rho(h_1)=\rho(h_2)=1$.
            The \eqref{eq: la uso una volta: coseni} implies that $\rho_2(a_2b_2)^\beta=1$. Moreover, the \eqref{eq: la uso una volta: coseni} implies that $\Tr \rho(a_1b_1)= \Tr \rho(a_2b_2)^\alpha$. In particular, we have that $\rho(a_1b_1)^\beta=1$. Hence, there exists two integers $k_1,k_2 \in \{1,\cdots \beta\}$ such that
           \begin{equation*}
                \begin{split}
                \Tr \rho(a_1b_1)=2 \cos\left( \frac{2\pi k_1}{\beta} \right) \in J_0(p_1,p_2), \quad \Tr \rho(a_2b_2) = 2\cos\left( \frac{2\pi k_2 }{\beta} \right) \in J_0(p_3,p_4), \\ \text{and} \quad 2 \cos\left( \frac{2\pi k_1}{\beta} \right)=2 \cos\left( \frac{2\pi \alpha k_2}{\beta} \right).
            \end{split}
            \end{equation*}
            We conclude that $\mathcal{S}_1 \cap \mathcal{S}_2$ is not empty.
            \end{proof}

            \begin{proof}[Proof of Claim \ref{claim: gli angli sono troppo grandi allora H10 H2}]
            If $2\beta < p_2$, meaning that $\pi - \nicefrac{4\pi}{p_2}> \pi - \nicefrac{2\pi}{\beta}$, then $2\cos(\nicefrac{2\pi k}{\beta}) \in \mathcal{S}_1$ if and only if $k \notin \frac{\beta}{2}\mathbb{Z}$.
            Let $k' \in \{1,\cdots,\beta\}$ be such that $2\cos\left(\nicefrac{2\pi k' \alpha}{\beta}\right) \in \mathcal{S}_2$. Such $k'$ exists as $\mathcal{S}_2$ is nonempty, in particular $k' \notin \frac{\beta}{2}\mathbb{Z}$. Since $\alpha$ and $\beta$ are coprime, we have that 
            $k' \alpha \notin \frac{\beta}{2}\mathbb{Z}$. This means that $2\cos\left(\nicefrac{2\pi k' \alpha}{\beta}\right) \in \mathcal{S}_1$. Hence, $\mathcal{S}_1 \cap \mathcal{S}_2 \neq \emptyset$. If $2 \beta < p_4$, then $\mathcal{S}_1 \cap \mathcal{S}_2 \neq \emptyset$ by an identical strategy.
                \end{proof}

                \begin{proof}[Proof of Claim \ref{claim: S3 S4 iff H10 H2pi}]
                    Let $\Sigma$ be the torus $\partial M_1 = \partial M_2 \subset M$. Let us suppose that there exists a representation $\rho \in \mathcal{R}(M)$ so that $\restr{\rho}{\fund{\Sigma}} \in H_{1,0}\cap H_{2,\pi}$, then $\rho(h_1)=1$ and $\rho(h_2)=-1$. The relations in \eqref{eq: conditizione per incollare due rappresentazioni dai due lati} imply that
            \[
                \rho_2(a_2b_2)^\beta (-1)^{\delta} = 1 , \quad \text{and} \quad \rho_1(a_1b_1)^\beta = \rho(a_2b_2)^{\alpha \beta}(-1)^{\gamma \beta} =(-1)^{-\alpha\delta+\gamma \beta} =-1.
            \]
            In particular, $\rho_2(a_2b_2)^{\beta}=(-1)^{\delta}$.
            This implies that there exists two integers $1\le k_1\le 2\beta -1$ and $1\le k_2\le \beta$ with $k_1$ odd, such that
            \[
            \Tr \rho_1(a_1b_1)= 2\cos \left( \frac{\pi k_1}{\beta} \right) \in J_0(p_1,p_2), \quad \text{and} \quad \Tr \rho_2(a_2b_2)= 2\cos \left( \frac{2\pi k_2}{\beta} - \frac{\pi \delta}{\beta} \right) \in J_{\pi}(2,p_2). 
            \]
            Furthermore, the \eqref{eq: conditizione per incollare due rappresentazioni dai due lati} implies that $\Tr\left( (- 1)^{\gamma}\rho_1(a_1b_1)\right)= \Tr \rho_2(a_2b_2)^{\alpha}$ and this implies that
            \[
                2\cos \left( \frac{\pi k_1}{\beta} - \gamma \pi \right) = 2 \cos \left( \frac{2\pi \alpha k_2}{\beta} - \frac{\delta \alpha}{\beta}\right).
            \]

            Conversely, let us suppose that $\mathcal{S}_3 \cap \mathcal{S}_4 \neq \emptyset$. Following a strategy similar to Claim \ref{claim: S1 e S2 danno H10 e H20 empty}, we construct the irreducible representations $\rho_1$ and $\rho_2$ of $\fund{M_1}$ and $\fund{M_2}$ respectively, such that
            \[
                \rho_1(h_1)=1, \quad \rho_2(h_2)=-1, \quad \text{and} \quad \restr{\rho_1}{\fund{\Sigma}}\equiv \restr{\rho_2}{\fund{\Sigma}}.
            \]
            This implies that $H_{1,0} \cap H_{2,\pi} \neq \emptyset$.
            \end{proof}

            \begin{proof}[Proof of Claim \ref{claim: angoli grossi per S3 e S4}]
                The proof follows from the same strategy of Claim \ref{claim: gli angli sono troppo grandi allora H10 H2}.
            \end{proof}

            \begin{proof}[Proof of Claim \ref{claim: S5 S6 iff H1pi H2pi}]
                Let us suppose that $\eta \in H_{1,\pi}\cap H_{2,\pi}$. Then, there exists a representation $\rho \colon \fund{M} \to SU(2)$ such that $\restr{\rho}{\fund{M_1}}$ (resp. $\restr{\rho}{\fund{M_1}}$) is irreducible and $\rho(h_1)=\rho(h_2)=-1$. The relations \eqref{eq: conditizione per incollare due rappresentazioni dai due lati} become
                \[
                    \rho_2(a_2b_2)^\beta=(-1)^{\delta +1} \quad \text{and} \quad \rho_1(a_1b_1)=\rho_2(a_2b_2)^\alpha(-1)^\gamma.
                \]
                Hence, $\rho_1(a_1b_1)^\beta = (-1)^{\alpha+1}$. According to Remark \ref{rmk: A n =-+1 se e solo se traccia} we obtain that there exist two integers $k_1$ and $k_2$ so that
                \[
                    2\cos \left( \frac{2\pi k_1}{\beta}+ \frac{\pi (\alpha+1)}{\beta}\right) \in J_\pi(2,p_2) \quad \text{and} \quad 2\cos \left( \frac{2\pi k_2}{\beta}+ \frac{\pi (\delta+1)}{\beta}\right) \in J_\pi(2,p_4).
                \]
                The condition $\rho_1(a_1b_1)=\rho_2(a_2b_2)^\alpha(-1)^\gamma$ implies that
                \[
                    \Tr \rho_1(a_1b_1) = 2\cos \left( \frac{2\pi k_1}{\beta}+ \frac{\pi (\alpha+1)}{\beta}- \pi\gamma\right) = 2\cos \left( \frac{2\pi \alpha k_2}{\beta}+ \frac{\pi \alpha(\delta+1)}{\beta}\right) = \Tr\rho(a_2b_2)^\alpha.
                \]
                Therefore $\mathcal{S}_5 \cap \mathcal{S}_6 \neq \emptyset$.

                Conversely, let us suppose that $\mathcal{S}_5 \cap \mathcal{S}_6 \neq \emptyset$. Let $\Sigma$ be the torus corresponding to $\partial M_1=\partial M_2 \subset M$ As in Claim \ref{claim: S1 e S2 danno H10 e H20 empty}, this implies the existence of the irreducible representations $\rho_1$ and $\rho_2$ of $\fund{M_1}$ and $\fund{M_2}$ respectively, such that
            \[
                \rho_1(h_1)=\rho_2(h_2)=-1 \quad \text{and} \quad \restr{\rho_1}{\fund{\Sigma}}\equiv \restr{\rho_2}{\fund{\Sigma}}.
            \]
            This implies that $H_{1,\pi} \cap H_{2,\pi} \neq \emptyset$.
                \end{proof}

            \begin{proof}[Proof of Claim \ref{claim: angoli grossi per S5 e S6}]
                The proof follows from the same strategy of Claim \ref{claim: gli angli sono troppo grandi allora H10 H2}.
            \end{proof}
\section{The main theorem}
\label{section: the main theorem}
In this section we will prove the two main results of this work.

\begin{lemma}\label{lemma: beta cannot be 2}
    If $\Delta(h_1,h_2)=2$ and $p_1=p_3=2$, then $M=M_1 \cup_{\Sigma} M_2$ is not $SU(2)$-abelian.
    \begin{proof}
        Let us choose two presentations so that $M_1=\mathbb{D}^2(\nicefrac{2}{1},\nicefrac{p_2}{q_2})$ and $M_2=\mathbb{D}^2(\nicefrac{2}{1},\nicefrac{p_4}{q_4})$.
        Notice that $g_1$ and $g_2$ are both either $1$ or $2$. We use the same notation as in Section \ref{sec: H1 cap H2}. 
        Lemma \ref{lemma: description of the rational longitude} implies that
        \[
            \Delta(\lambda_{M_2},h_1) = \frac{2}{o_2g_2} \left|p_4 \delta - 2 q_4 - p_4\right| \quad \text{and} \quad \Delta(\lambda_{M_1},h_2)= \frac{2}{o_1g_1} \left|p_2 \alpha + 2 q_2 + p_2\right|.
        \]

       Since $g_1 \le g_2 \le 2$, we divide the proof in three cases:
        \begin{enumerate}[label=\roman*)]
            \item $g_1= g_2=1$,
            \item $g_1=1$ and $g_2=2$,
            \item $g_1=g_2=2$.
        \end{enumerate}
        
        \textbf{Case \romannum{1})} $\boldsymbol{g_1= g_2=1}$. According to Lemma \ref{lemma: description of the rational longitude}, $o_1=o_2=1$. In particular $\Delta(\lambda_{M_1},h_2)$ is even (or zero). Proposition \ref{prop: H2 and A1 sono disjoint if and only if} implies that $H_2 \cap A_1$ is empty if and only if $\Delta(\lambda_{M_1},h_2)=0$. However, since $g_1=1$, $\Delta(\lambda_{M_1},h_2)=0$ if and only if $p_2$ divides $q_2$. This cannot happen since $\gcd(p_2,q_2)=1$ and $p_2 \ge 2$. We conclude that $H_1 \cap A_2$ is nonempty.

         \textbf{Case \romannum{2}) $\boldsymbol{g_1=1}$ and $\boldsymbol{g_2=2}$}. Again, $\Delta(\lambda_{M_1},h_2)$ is even and cannot be zero, Proposition \ref{prop: H2 and A1 sono disjoint if and only if} states that $A_1$ and $H_2$ are disjointed if and only if $\Delta(\lambda_{M_1},h_2)=4$, $p_3=2$ and $p_4=4$. This implies that $o_2=1$. Thus, $\Delta(\lambda_{M_2},h_1) = 2|2\delta -q_4-2|$ is even. According to Proposition \ref{prop: H1 and A2 sono disjoint if and only if}, $A_2 \cap H_1= \emptyset$ if and only if $\Delta(\lambda_{M_2},h_1)=0$. Since $|2\delta -q_4 -2|$ is an odd number, we conclude that $\Delta(\lambda_{M_2},h_1)\neq0$ and that $A_2 \cap H_1 \neq \emptyset$.

        \textbf{Case \romannum{3}) $\boldsymbol{g_1=g_2=2}$}. We study the cases $o_1=1$ and $o_1=2$ separately.
        
        If $o_1=1$, then $\Delta(\lambda_{M_1},h_2)=|p_2 \alpha  + 2q_2 + p_2|$ is even. Proposition \ref{prop: H2 and A1 sono disjoint if and only if} implies that $H_2 \cap A_1$ is empty if and only if $\Delta(\lambda_{M_1},h_2)$ is either $0$ or $4$.
        
        If $\Delta(\lambda_{M_1},h_2)=0$, then $p_2$ divides $2q_2$ and then $p_2=2$. This is a contradiction because if $p_1=p_2=2$ then the order of its rational longitude is $2$ according to Lemma \ref{lemma: description of the rational longitude}.
        
        If $\Delta(\lambda_{M_1},h_2)=4$ then $p_3=2$ and $p_4=4$. Moreover, $\Delta(\lambda_{M_2},h_1) = 2|2\delta -q_4-2|$ is even. Proposition \ref{prop: H1 and A2 sono disjoint if and only if} implies that $\Delta(\lambda_{M_2},h_1)$ has to be either $0$ or $4$. Both cases are are not possible since $|2\delta -q_4 -2|$ is an odd number. We conclude that if $g_1=g_2=2$ and $o_1 \equiv_2 1$, then $M$ is not $SU(2)$-abelian.
        
       If $o_1=2$, Lemma \ref{lemma: description of the rational longitude} implies that $p_2=2$. Proposition \ref{prop: H1 and A2 sono disjoint if and only if} implies that $\Delta(\lambda_{M_2},h_1)$ is either $0$ or $1$. If $\Delta(\lambda_{M_2},h_1)=0$, then $p_4$ divides $2$, hence $p_4=2$.
       If $\Delta(\lambda_{M_2},h_1)=1$, since $(p_4 \delta - 2q_4 -p_4)$ is an even number, we have that $o_2=2$ and again $p_4=2$.
        Hence, if $g_1=g_2=o_1=2$, then $p_1=p_2=p_3=p_4=2$.
        The conclusion holds by Example \ref{exmp: two evil manifolds}.
        \end{proof}
\end{lemma}

\begin{lemma}
    If either $(p_1,p_2,|\beta|)$ or $(p_3,p_4,|\beta|)$ is in the set $\{(2,4,4),(3,3,3)\}$, then $M=M_1 \cup_{\Sigma} M_2$ is not $SU(2)$-abelian.
    \label{lemma: some edge cases}
    \begin{proof}
        In this proof use the same notation we used in Lemma \ref{lemma: H1 cap H2 empty if and only if beta 1 or 2}. We divide the proof in four cases:
        \begin{enumerate}[label=\roman*)]
            \item $(p_1,p_2,|\beta|)=(3,3,3)$,
            \item $(p_1,p_2,|\beta|)=(2,4,4)$,
            \item $(p_3,p_4,|\beta|)=(2,4,4)$,
            \item $(p_3,p_4,|\beta|)=(3,3,3)$.
        \end{enumerate}
    
       \textbf{Case \romannum{1})} $\boldsymbol{(p_1,p_2,|\beta|)=(3,3,3)}$. Since $p_1=p_2=3$, $g_1=3$.
       Proposition \ref{prop: when P1 P2 contains SU(2)-abelian representations} and Theorem \ref{teo: M SU(2)-abeliano se e solo se i pezzi sono empty} implies that $M$ is not $SU(2)$-abelian.

        \textbf{Case \romannum{2})} $\boldsymbol{(p_1,p_2,|\beta|)=(2,4,4)}$. Since $\beta$ is even, $\gamma$ and $\delta$ are odd. The intersection points between $\{\psi_1=\pi\}$ and $L_\pi$ have coordinates
        \[
            (\theta_1,\psi_1) = \left( \frac{\pi}{2}k\alpha + \gamma\pi, \pi\right) \quad \text{and} \quad (\theta_2,\psi_2) = \left( \frac{\pi}{2}k , \pi\right), 
        \]
        with $k \in \{0,1,2,3\}$.
        Lemma \ref{lemma: surjective maps} implies that the points $\left( \frac{\pi}{2}k\alpha + \gamma\pi, \pi\right)$ are contained in $H_{1,\pi}$. 
        Thus, the set $H_{1,\pi}\cap H_{2,\pi}$ is empty if and only if $p_3=p_4=4$ as a consequence of Lemma \ref{lemma: surjective maps} and Lemma \ref{lemma: the interval Jpi(p1,p2) contains 0 if and only if p1,p2 is not 4,4}. If $p_3=p_4=4$, then Proposition \ref{prop: H2 and A1 sono disjoint if and only if} implies that $A_1 \cap H_2$ is empty if and only if $\Delta(\lambda_{M_1},h_2)=2$ and $o_1=1$. Since $p_1=2$ and $p_2=4$, we obtain that
        \[
            \Delta(\lambda_{M_1},h_2)= 4| \alpha +q_2+2q_1| \neq 2.
        \]
        Thus, $A_1 \cap H_2 \neq \emptyset$. The conclusion follow from Theorem \ref{teo: M SU(2)-abeliano se e solo se i pezzi sono empty}.

        \textbf{Case \romannum{3})} $\boldsymbol{(p_3,p_4,|\beta|)=(2,4,4)}$. This case is proven as in the previous one by switching the roles of $(p_1,p_2)$ and $(p_3,p_4)$.

        \textbf{Case \romannum{4})} $\boldsymbol{(p_3,p_4,|\beta|)=(3,3,3)}$. The intersection points in $\{\psi_1=\pi\}\cap L_0$ have coordinates
        \[
            (\theta_1,\psi_1) = \left( \frac{\pi \alpha}{3} +\frac{2\pi \alpha k}{3}, \pi\right) \quad \text{and} \quad (\theta_2,\psi_2) = \left(  \frac{\pi}{3} +\frac{2\pi k}{3}, 0\right),
        \]
        with $k \in \{0,1,2\}$. Since $J_0(p_3,p_4)=J_0(3,3)=(-1,2)$, Lemma \ref{lemma: surjective maps} implies that the point $(\theta_2,\psi_2)=(\nicefrac{\pi}{3},0)$ is in $H_{2,0}$.
        Lemma \ref{lemma: description of the set S(p1,p2)} implies that, $H_{1,\pi}\cap H_{2,0}$ is empty if and only if $p_1=p_2=3$ and $\alpha=2$. As we stated before, this would imply that $3=g_1 \le g_2$ and, as a consequence of Proposition \ref{prop: when P1 P2 contains SU(2)-abelian representations}, the existence of an irreducible $SU(2)$-representation of $\fund{M}$.
    \end{proof}
\end{lemma}

\begin{cor}\label{cor: M SU(2) ab then beta =1}
    If $M=M_1 \cup_{\Sigma} M_2$ is $SU(2)$-abelian, then $\Delta(h_1,h_2)=1$.
    \begin{proof}
        The conclusion is a consequence of Lemma \ref{lemma: if b >= 3 then H1 and H2 contain an intersection}, Lemma \ref{lemma: beta cannot be 2} and Lemma \ref{lemma: some edge cases}.
    \end{proof}
\end{cor}

\begin{rmk} If $\Delta(h_1,h_2)=1$ and $\Delta(\lambda_{M_1},h_2)=0$, then $\alpha p_1p_2 \pm (p_1q_2+p_2q_1)=0$. This implies that $p_1=p_2$. Similarly, if $\Delta(h_1,h_2)=1$ and $\Delta(\lambda_{M_2},h_1)=0$, then $p_3=p_4$.
    \label{Please Note: Delta =0 implica p1=p2}
\end{rmk}

We remind the reader that the condition $A$, $B$, and $C$ are reported in the Introduction.
Condition $A$ holds if and only the hypoteses of Proposition \ref{prop: when P1 P2 contains SU(2)-abelian representations} are satisfied. Condition $B$ holds if and only if hypoteses of Proposition \ref{prop: H1 and A2 sono disjoint if and only if} are satisfied with the additional requirement that $g_1=\gcd(p_1,p_2)\le 2$. Remark \ref{Please Note: Delta =0 implica p1=p2} states that if $\Delta(h_1,h_2)=1$ and $\Delta(\lambda_{M_1},h_2)=0$, then $p_1=p_2$. In particular, if $g_1 \le 2$, then $p_1=p_2=2$.
The order of the rational longitude of the manifold $\mathbb{D}^2(\nicefrac{2}{q_1},\nicefrac{2}{q_2})$ is $2$. Hence, condition $C$ holds if and only if hypoteses of Proposition \ref{prop: H2 and A1 sono disjoint if and only if} are satisfies with the additional requirements that $g_1 \le 2$ and $\Delta(h_1,h_2)=1$.

\begin{repteo}{MAIN THEOREM}
    Let manifold $M=M_1 \cup_{\Sigma} M_2$ be the graph manifold as above and let us suppose that $g_1 \le g_2$. The manifold $M$ is $SU(2)$-abelian if and only if $\Delta(h_1,h_2)=1$ and the conditions $A$, $B$, and $C$ hold.
    \begin{proof}
        If $\Delta(h_1,h_2)=1$ and conditions $A$, $B$ and $C$ hold, then we have that the sets $H_1\cap H_2$, $A_1 \cap H_2$, $H_1 \cap A_2$, and $P_1 \cap P_2$ are empty as a consequence of Lemma \ref{lemma: beta =1 then H1 intersection H2 is empty}, Proposition \ref{prop: H2 and A1 sono disjoint if and only if}, Proposition \ref{prop: H1 and A2 sono disjoint if and only if}, and Proposition \ref{prop: when P1 P2 contains SU(2)-abelian representations}. The conclusion holds by Theorem \ref{teo: M SU(2)-abeliano se e solo se i pezzi sono empty}.

        Conversely, if $M$ is $SU(2)$-abelian, then Corollary \ref{cor: M SU(2) ab then beta =1} implies that $\Delta(h_1,h_2)=1$.
        Theorem \ref{teo: M SU(2)-abeliano se e solo se i pezzi sono empty} implies the sets $P_1 \cap P_2$, $H_1\cap A_2$, and $A_1 \cap H_2$ are empty. By Proposition \ref{prop: when P1 P2 contains SU(2)-abelian representations} condition $A$ holds. In particular, we have that $g_1\le2$. Since $H_1 \cap A_2$ and $A_1 \cap H_2$ are empty and $g_1 \le 2$, Condition $B$ and Condition $C$ hold as a consequence of Proposition \ref{prop: H1 and A2 sono disjoint if and only if}, Proposition \ref{prop: H2 and A1 sono disjoint if and only if}, and Remark \ref{Please Note: Delta =0 implica p1=p2}. 
    \end{proof}
\end{repteo}

Let us suppose $M=M_1\cup_{\Sigma}M_2$ is an $SU(2)$-abelian manifold, we want to give a description of $M_1$ and $M_2$ in terms of the (non-unique) Seifert coefficients.

\begin{lemma}\label{lemma: congruence conditions}
    Let $M_1 = \mathbb{D}^2 (\nicefrac{p_1}{q_1},\nicefrac{p_2}{q_2})$ and $M_2 = \mathbb{D}^2 (\nicefrac{p_3}{q_3},\nicefrac{p_4}{q_4})$ and let $\varphi: \partial M_1 \to \partial M_2$ be an orientation reversing diffeomorphism. If $M_1 \cup_{\varphi} M_2$ is $SU(2)$-abelian, then
    \begin{align*}
          \pm o_1g_1\Delta(\lambda_{M_1},h_2) \equiv_{p_1p_2} p_1q_2 + p_2q_1 \quad \text{and} \quad \pm o_2g_2\Delta(\lambda_{M_2},h_1) \equiv_{p_3p_4} p_3q_4 + p_4q_3.
    \end{align*}
    \begin{proof}
Lemma \ref{lemma: description of the rational longitude} implies that
\[
    o_1g_1 \Delta(\lambda_{M_1},h_2)= \Delta \left( \begin{bmatrix}
        \alpha & \beta \\ \gamma & \delta
    \end{bmatrix} \begin{pmatrix}
        p_1p_2 \\ p_1q_2+p_2q_1
    \end{pmatrix}, \begin{bmatrix}
        0 \\ 1
    \end{bmatrix} \right)= \left| \alpha p_1p_2 + \beta (p_1q_2+p_2q_1) \right|.
\]
    Corollary \ref{cor: M SU(2) ab then beta =1} implies that $\beta= \pm 1$. Hence, $\pm o_1g_1 \Delta(\lambda_{M_1},h_2) \equiv_{p_1p_2} p_1q_2+p_2q_1$. The second half of the conclusion holds similarly.
    \end{proof}
\end{lemma}

\begin{prop}\label{prop: l'incollamento esiste sempre se soddisma le CC}
Let $M_1 = \mathbb{D}^2 (\nicefrac{p_1}{q_1},\nicefrac{p_2}{q_2})$ and $M_2 = \mathbb{D}^2 (\nicefrac{p_3}{q_3},\nicefrac{p_4}{q_4})$. If there exist two integers $n$ and $m$ such that
    $
         o_1g_1 n \equiv_{p_1p_2} p_1q_2 + p_2q_1$ and $o_2g_2 m \equiv_{p_3p_4} p_3q_4 + p_4q_3,
    $
then there exists an orientation reversing diffeomorphism $\varphi \colon \partial M_1 \to \partial M_2$ such that 
\[
    \Delta(\varphi(h_1),h_2)=1, \quad \Delta(\varphi(\lambda_{M_1}),h_2)=|n|, \quad \text{and} \quad \Delta(\varphi(h_1),\lambda_{M_1})=|m|.
\]
\begin{proof}
    By hypothesis there exists two integers $k_1$ and $k_2$ such that
    \[
         o_1g_1 n = k_1p_1p_2 + p_1q_2 + p_2q_1 \quad \text{and} \quad o_2g_2 m = k_2 p_3p_4 + p_3q_4 + p_4q_3.
    \]
    Let us define the matrix $A$ as
    \[
        A=\begin{bmatrix}
            -k_1 & -1 \\ k_1k_2-1 & k_2
        \end{bmatrix}.
    \]
    The matrix $A$ has determinant equal to $-1$. Hence, there exists an orientation reversing diffeomorphism $\varphi \colon \partial M_1 \to \partial M_2$ such that $\varphi_\ast = A$ with respect to the ordered bases $\{\mu_1,h_1\}$ and $\{\mu_2,h_2\}$.
    In particular $\Delta(\varphi(h_1),h_2)=1$.
    By Lemma \ref{lemma: description of the rational longitude} we get that
    \[
    \Delta\left(\varphi(\lambda_{M_1}),h_2\right) = \Delta\left(\lambda_{M_1},\varphi^{-1}(h_2)\right) = \frac{1}{o_1g_1} \left|k_1p_1p_2 + (p_1q_2+p_2q_1) \right| = |n|
    \]
    Similarly, we have that
    \[
        \Delta\left(\varphi(h_1),\lambda_{M_2}\right) = \frac{1}{o_2g_2} \left|k_2p_3p_4 + (p_3q_4+p_4q_3) \right| = |m|
    \]
\end{proof}
\end{prop}

If $M$ is $SU(2)$-abelian, then $\Delta(h_1,h_2)=1$ by Corollary \ref{cor: M SU(2) ab then beta =1}. Thus, the set of the homotopy classes of the regular fibers $h_1$ and $h_2$ forms a basis for $\fund{\partial M_1}=\fund{\partial M_2}$. Therefore, there exist integers $n_1,m_1,n_2,m_2$ so that
\[
    \lambda_{M_1}=n_1h_1 + m_1h_2 \subset \partial M_1 \quad \text{and} \quad \lambda_{M_2}=n_2h_1 + m_2h_2 \subset \partial M_2.
\]
Then $|n_i|=\alg{\lambda_{M_i}}{h_2}$ and $|m_i|=\alg{\lambda_{M_i}}{h_1}$ for $i \in \{1,2\}$.\\Thus,
$
\Delta(\lambda_{M_1},\lambda_{M_2})= \big|n_1m_2 - n_2m_1 \big|.
$
If $M_1$ and $M_2$ are presented as in \eqref{eq: presentation of the fundamental group of M1 and M2}, then, according to Lemma \ref{lemma: description of the rational longitude}, $\alg{\lambda_{M_1}}{h_1} = \frac{p_1p_2}{o_1g_1}$ and $\alg{\lambda_{M_2}}{h_2} = \frac{p_3p_4}{o_2g_2}$. Therefore,
\begin{equation}
\Delta(\lambda_{M_1},\lambda_{M_2}) \equiv_2 \alg{\lambda_{M_1}}{h_2}\alg{\lambda_{M_2}}{h_1} +\frac{p_1p_2p_3p_4}{o_1g_1o_2g_2} .
\label{eq: parity of rational logitudes intersection}
\end{equation}

\begin{repteo}{teo: classification}
     Let $M=M_1\cup_{\Sigma}M_2$ be as in Theorem \ref{MAIN THEOREM}, and further suppose that $M_1$ and $M_2$ are presented in such a way that $0<q_i<p_i$ and $g_1 \le g_2$. The manifold $M$ is $SU(2)$-abelian if and only if $\Delta(h_1,h_2)=1$ and it lays in one of the seven classes in Table \ref{table: classification}.
    \begin{proof}
    
        If a manifold $M$ is homeomorphic to a manifold in one of the $7$ classes in Table \ref{table: classification} and $\Delta(h_1,h_2)=1$, 
        then conditions $A$, $B$ and $C$ holds. Hence, $M$ is $SU(2)$-abelian by Theorem \ref{MAIN THEOREM}.

        Conversely, let us suppose that $M$ is $SU(2)$-abelian.
        Corollary \ref{cor: M SU(2) ab then beta =1} implies that $\Delta(h_1,h_2)=1$. The conditions $A$, $B$, and $C$ follow by Theorem \ref{MAIN THEOREM}. In particular,
        Condition $B$ implies that $\Delta(\lambda_{M_2},h_1)$ is in $\{0,1,4\}$. We divide the proof in three cases:
        \begin{enumerate}[label=\roman*)]
            \item $\Delta(\lambda_{M_2},h_1)=0$,
            \item $\Delta(\lambda_{M_2},h_1)=4$,
            \item $\Delta(\lambda_{M_2},h_1)=1$.
        \end{enumerate}
        
        \textbf{Case \romannum{1}) $\boldsymbol{\Delta(\lambda_{M_2},h_1)=0}$}. According to Condition $B$, $p_1=2$ and $o_2 \equiv_2 1$. Remark \ref{Please Note: Delta =0 implica p1=p2} implies that $p_3=p_4$ and in particular $o_2=\gcd(p_3,q_3+q_4)$. 
        Lemma \ref{lemma: congruence conditions} implies $q_3+q_4 \equiv_{p_3} 0$ and hence that $o_2=p_3$. Since $o_2$ is odd by assumption, $p_3$ is odd as well.
        Since $p_3=p_4$, Condition $C$ states that $\Delta(\lambda_{M_1},h_2)$ is either $1$, $2$ or $3$.

        If $\Delta(\lambda_{M_1},h_2)=1$, then the manifold $M$ is in the class $(1)$.
        
        If $\Delta(\lambda_{M_1},h_2)=2$, then $p_3=p_4=4$ by Condition $C$. Thus, according to Lemma \ref{lemma: description of the rational longitude}, $o_2\equiv_2 0$. This is a contradiction since we supposed that $o_2\equiv_2 1$.
        
        If $\Delta(\lambda_{M_1},h_2)=3$, then Condition $C$ implies that $p_3=p_4=3$.
        Lemma \ref{lemma: congruence conditions} implies that
        $
            3(q_3+q_4)\equiv_9 0,
        $
        and hence $q_3\equiv_31$ and $q_4\equiv_32$. Therefore, the manifold $M$ is in the class $(2)$.
        
        \textbf{Case \romannum{2}) $\boldsymbol{\Delta(\lambda_{M_2},h_1)=4}$}. As a consequence of Condition $B$, $p_1=2$, $p_2=4$, and $o_2=1$. Lemma \ref{lemma: description of the rational longitude} implies that $o_1=1$ and
        \[
        \Delta(\lambda_{M_1},h_2)= \frac{1}{o_1g_1} \left| \alpha p_1p_2 + \beta \left(p_1q_2+p_2q_1 \right) \right| =  \left| 4\alpha + \pm \left(q_2+2q_1 \right) \right|\equiv_2 1. 
        \]
        Hence, by Condition $C$, $\Delta(\lambda_{M_1},h_2)$ is either $1$ or $3$.
        
        If $\Delta(\lambda_{M_1},h_2)=1$ then $M$ lies in the class $(3)$. If $\Delta(\lambda_{M_1},h_2)=3$, then Condition $C$ implies that $p_3=p_4=3$. Lemma \ref{lemma: congruence conditions} implies that
        \[
        3(q_3+q_4) \equiv_9 \pm 12 \equiv_9 \mp 3.
        \]
        Therefore, $q_3 +q_4 \equiv_3 \pm 1$. This implies that $q_3\equiv_3 q_4$ and that $M$ is in the class $(4)$.

        \textbf{Case \romannum{3}) $\boldsymbol{\Delta(\lambda_{M_2},h_1)=1}$}. Condition $B$ implies that $o_2 \le 2$. Condition $C$ states that $\Delta(\lambda_{M_1},h_2)\in\{1,2,3,4\}$.
        
        If $\Delta(\lambda_{M_1},h_2)=4$, then, up to switching $M_1$ and $M_2$, the manifold $M$ is in the class $(3)$.
        
        If $\Delta(\lambda_{M_1},h_2)=3$, then as a consequence of Condition $C$, $p_3=p_4=3$, $o_1=1$ and $\Delta(\lambda_{M_1},\lambda_{M_2})\equiv_20$. Since $o_2$ divides $g_2=3$ by Lemma \ref{lemma: description of the rational longitude}, $o_2=1$ and hence $\gcd(3,q_3+q_4)=1$. Thus, $q_3 \equiv_3 q_4$. The identity \eqref{eq: parity of rational logitudes intersection} shows that
        \[
            \Delta(\lambda_{M_1},\lambda_{M_2}) \equiv_2 \frac{3p_1p_2}{g_1} + 3.
        \]
        Thus, $\Delta(\lambda_{M_1},\lambda_{M_2})$ is even if and only if $p_1p_2/g_1$ is odd and hence if and only if $p_1p_2 \equiv_2 1$ and $g_1=1$. This implies that the manifold lies in the class $(5)$.
        
        Let us suppose that $\Delta(\lambda_{M_1},h_2)=2$. As before, Condition $C$ imposes that $p_3=p_4=4$ and $o_1=1$. In this case the manifold $M$ is in class $(6)$.
        Finally, let $\Delta(\lambda_{M_1},h_2)=1$. According to Lemma \ref{lemma: description of the repersentation P1} $g_1 \le 2$. Since $o_1$ divides $g_1$ by Lemma \ref{lemma: description of the rational longitude}, $o_1 \le 2$ and the manifold $M$ is in class $(7)$.
    \end{proof}
\end{repteo}

We conclude this section we an explicit description of the class $4$ in Table \ref{table: classification}. A combinatorial calculation gives that if $M$ is of class $4$ in Table \ref{table: classification}, then it is one of the following manifolds:

\begin{equation}
    \begin{alignedat}{2}
    &\mathbb{D}^2\left( \frac{2}{1},\frac{4}{1}\right) \bigcup_{\pm\left[\begin{smallmatrix}
            0 & 1 \\
            1 & -2
        \end{smallmatrix}\right]} \mathbb{D}^2\left(\frac{3}{1},\frac{3}{1} \right), \qquad \qquad  &&\mathbb{D}^2\left( \frac{2}{1},\frac{4}{3}\right) \bigcup_{\pm\left[\begin{smallmatrix}
            -2 & 1 \\
            5 & -2
        \end{smallmatrix}\right]} \mathbb{D}^2\left(\frac{3}{1},\frac{3}{1} \right),\\[5pt]
        &\mathbb{D}^2\left( \frac{2}{1},\frac{4}{3}\right) \bigcup_{\pm \left[\begin{smallmatrix}
            -2 & 1 \\
            1 & 0
        \end{smallmatrix}\right]} \mathbb{D}^2\left(\frac{3}{2},\frac{3}{2} \right), &&\mathbb{D}^2\left( \frac{2}{1},\frac{4}{1}\right) \bigcup_{\pm \left[\begin{smallmatrix}
            0 & 1 \\
            1 & 0
        \end{smallmatrix}\right]} \mathbb{D}^2\left(\frac{3}{2},\frac{3}{2} \right). 
        \end{alignedat}
        \label{eq: the 8 manifolds in class 4}
\end{equation}
It can be proven the manifold in the right hand corner of \eqref{eq: the 8 manifolds in class 4} has positive Betti number. Furthermore, the remaining manifolds in \eqref{eq: the 8 manifolds in class 4} are rational homology spheres.

\section{Applications}
\label{sec: applications} 

Let $p,q,r,s \in \mathbb{Z}$ be such that $\gcd(p,q)=\gcd(r,s)=1$. Let $E(T_{p,q})$ be the exterior of a open tubular neighborhood of the torus knot $T(p,q) \subset S^3$ with knot meridian $\mu_{p,q}$. It is known that $E(T_{p,q})$ is a Seifert fibred space, let us denote its regular fiber by $h_{p,q}$. We denote by $Y_{(p,q),(r,s)}$ the graph manifold obtained by gluing $E(T_{p,q})$ and $E_{r,s}$ along the diffeomorphism $\varphi : \partial E(T_{p,q}) \to \partial E_{r,s}$ such that $\varphi(\mu_{p,q})= h_{r,s}$ and $\varphi(h_{p,q})=\mu_{r,s}$. In \cite{Motegi_haken_manifolds} Motegi proved that the manifold $Y_{(p,q)(r,s)}$ is $SU(2)$-abelian.

\begin{teo}\label{teo: Y(p,q,r,s) is kinda of unique}
    Let $M_1$ and $M_2$ be as in Theorem \ref{MAIN THEOREM} and let us suppose that $g_1=g_2=1$. If $M=M_1 \cup_{\varphi} M_2$ is $SU(2)$-abelian, then $M=M_1 \cup_{\Sigma} M_2$ is diffeomorphic to $Y_{(p,q),(r,s)}$, for suitable $p,q,r,s \in \mathbb{Z}$ with $\gcd(p,q)=\gcd(r,s)=1$.
    \label{teo: esiste un'unico incollamento per i nodi torici}
    \begin{proof}
        Theorem \ref{teo: classification} implies that $\Delta(\lambda_{M_1},h_2)=1$ and $\Delta(\lambda_{M_2},h_1)=1$. Lemma \ref{lemma: congruence conditions} states that
        \[
          p_1q_2+p_2q_1 \equiv_{p_1p_2} \pm 1 \quad \text{and} \quad p_3q_4+p_4q_3 \equiv_{p_3p_4} \pm 1.
        \]
        Hence, there exist two torus knots $T_1$ and $T_2$ embedded in $S^3$ such that $M_1 = S^3 \setminus \nu(T_1)$ and $M_2 = S^3 \setminus \nu(T_2)$, where $\nu(K)$ is small open neighborhood of the knot $K \subset S^3$. Let us denote as $\mu_{T_1}$ (resp. $\mu_{T_2}$) the meridian of $T_1\subset S^3$ (resp. $T_2 \subset S^3$). 
        Up to changing the presentations of $M_1$ and $M_2$, we can suppose that
        \begin{equation}
            p_1q_2+p_2q_1 = \pm 1 \quad  \text{and} \quad p_3q_4+p_4q_3= \pm 1.
            \label{lo uso sono una volta: presetnazioni dei due nodi torici}
        \end{equation}
        Let $\mu_1$ and $\mu_2$ be the corresponding fibration meridians of the chosen presentations. Notice that $\mu_1=\mu_{T_1}$ and $\mu_2=\mu_{T_2}$ in $\fund{M_1}$ and $\fund{M_2}$ respectively.
        Lemma \ref{lemma: description of the rational longitude} implies that
        \[
        \Delta(\lambda_{M_1},h_2 )= \left| \alpha p_1p_2 \pm 1 \right|=1 \quad \text{and} \quad \Delta(\lambda_{M_2},h_1 )=\left| \delta p_3p_4 \pm 1 \right|=1.
        \]
        Since $|p_1p_2| \ge 6$ and $|p_3p_4|\ge 6$, we obtain that $|\alpha|=\Delta(\mu_{T_1},h_2)=0$, $|\delta|=\Delta(\mu_{T_2},h_1)=0$. This implies that $\varphi(h_1)=\mu_{T_2}$ and $\varphi(h_2)=\mu_{T_1}$. In other words, if $T_1=T(p,q)$ and $T_2=T(r,s)$, then $\varphi(\mu_{p,q})= h_{r,s}$ and $\varphi(h_{p,q})=\mu_{r,s}$. Thus, $M_1 \cup_{\Sigma} M_2$ is diffeomorphic to $Y_{(p,q),(r,s)}$
    \end{proof}
\end{teo}

\begin{repcor}{cor: torus knots complements are unique}
    For $i \in \{1,2 \}$, let $E(T_i)$ be the exterior of a open tubular neighborhood of the torus knot $T_i \subset S^3$. We denote by $\lambda_i$ and $\mu_i$ the null-homologous longitude and the meridian of $T_i$. The manifold $E(T_1)\cup_{\Sigma} E(T_2)$ is $SU(2)$-abelian if and only if $\Delta(\lambda_1,\mu_2)=0$ and $\Delta(\lambda_2,\mu_1)=0$.
    \begin{proof}
        If $E(T_1)\cup_{\Sigma} E(T_2)$ is $SU(2)$-abelian, then the conclusion holds by Theorem \ref{teo: esiste un'unico incollamento per i nodi torici}. Conversely, if $\Delta(\lambda_1,\mu_2)=0$ and $\Delta(\lambda_2,\mu_1)$ then the conclusion holds by \cite{Motegi_haken_manifolds}.
    \end{proof}
\end{repcor}

Given two Seifert pieces $M_1$ and $M_2$ with $g_1=g_2=1$, Theorem \ref{teo: Y(p,q,r,s) is kinda of unique} implies that there exists a unique $SU(2)$-abelian $3$-manifold obtained by gluing $M_1$ and $M_2$. In general, it is not true that for a given couple of Seifert fibred manifolds there exists a unique gluing that produces an $SU(2)$-abelian manifold, as is shown in the following example.

\begin{exmp}
    Let $M_1 =\mathbb{D}^2(\nicefrac{4}{1},\nicefrac{5}{4})$ and $M_2 = \mathbb{D}^2(\nicefrac{2}{1},\nicefrac{2}{1})$. Let us choose the usual basis $\{\mu_1,h_1\}$ and $\{\mu_1,h_1\}$ for $\fund{\partial M_1}$ and $\fund{\partial M_2}$. Let $\varphi_{1}$ and $\varphi_2$ be the two orientation reversing diffeomorphisms $\partial M_1 \to \partial M_2$ such that
    \[
        \varphi_1^\ast = \begin{bmatrix}
            -1 & 1 \\ 1 & 0
        \end{bmatrix} \quad \text{and} \quad
        \varphi_2^\ast = \begin{bmatrix}
            -1 & 1 \\ -1 & 2
        \end{bmatrix}.
    \]
    Theorem \ref{MAIN THEOREM} implies that $M_1 \cup_{\varphi_{1}} M_2$ and $M_1 \cup_{\varphi_2} M_2$ are $SU(2)$-abelian. Theorem \ref{teo: classification} implies that $M_1 \cup_{\varphi_{1}} M_2$ and $M_1 \cup_{\varphi_2} M_2$ are both in class $7$ of Table \ref{table: classification}. It can be proven that $\Delta(\varphi_1(\lambda_{M_1}),\lambda_{M_2})=19$ and that $\Delta(\varphi_2(\lambda_{M_1}),\lambda_{M_2})=21$. Thus, $M_1 \cup_{\varphi_{1}}M_2$ is not diffeomorphic to $M_1 \cup_{\varphi_{2}}M_2$.
    \label{exmp: non unique diffeomorphism}
    \end{exmp}
\section{L-spaces}
\label{sec: L-spaces}

Let $p$,$q$,$r$, and $s$ be integers such that $\gcd(p,q)=\gcd(r,s)=1$. In \cite{classofknots} Zentner proved that if $pqrs$ is an even number, then the graph manifold $Y_{(p,q)(r,s)}$ is the double branched cover of an alternating knot. This implies that $Y_{(p,q)(r,s)}$ is a strong L-space. 
Moreover, in \cite{zhang2019remarks} the Zhang proved that the manifold $Y_{(p,q)(r,s)}$ is an L-space. In general, every known $SU(2)$-abelian rational homology sphere is an Heegaard Floer L-space.

Let $M=M_1 \cup_{\varphi} M_2$ be a manifold as in Theorem \ref{MAIN THEOREM}, we are going to give a further evidence of Conjecture \ref{conj: SU(2)-abelian implies L-space} by proving that if $M$ is an $SU(2)$-abelian rational homology sphere, then $M$ is an L-space.

\begin{defn}
    Let $Y$ be a compact oriented $3$-manifold with torus boundary. We denote by $\calS(\partial Y)$ the set of slopes in $\partial Y$. We define the
    L-space interval of $Y$ to be
    \[
    \mathcal{L}(Y) = \left\{ \gamma \in \calS(\partial Y) \,\middle|\, Y(\gamma) \text{ is an L-space} \right\}.
    \]
\end{defn}
For a more detailed description of $\mathcal{L}(Y)$ we give \cite{LspaceIntervalGraphManifold} as a reference.
Let $Y$ be a $3$-manifold with torus boundary. The set $\calS(\partial Y)$ can be identified as the projectivization of the first homology  $\mathbb{P}\left(H_1(\partial Y ; \mathbb{Z})\right)$. Thus, there exists a canonical embedding 
\[
\calS(\partial Y) \hookrightarrow \mathbb{P}\left(H_1(\partial Y ; \mathbb{R})\right) \cong S^1 = \mathbb{R} \cup \{\infty\}.
\]

\begin{teo}[{\cite[Proposition 1.5]{LspaceIntervalGraphManifold}}]
    If $Y_1$ and $Y_2$ are non-solid-torus graph manifolds with torus boundary, then the union $Y_1 \cup_{\varphi} Y_2$, with gluing map $\varphi \colon \partial Y_1 \to - \partial Y_2$, is an L-space
    if and only if
    \[
        \varphi_\ast^{\mathbb{P}}\left( \mathcal{L}^\circ(Y_1)\right) \cup \mathcal{L}^\circ(Y_2) = \mathbb{P} \left( H_1 \left( \partial Y_2;\mathbb{Q}\right)\right).
    \]
    \label{teo: L space gluing graph manifold}
\end{teo}

By Corollary \ref{cor: M SU(2) ab then beta =1}, if $M=M_1 \cup_{\varphi} M_2$ is $SU(2)$-abelian, then $\Delta(h_1,h_2)=1$. Let us call $\Sigma$ the embedded torus in $M=M_1 \cup_{\varphi} M_2$ corresponding to $\partial M_1 = \partial M_2$.
Let us fix the convention according to which the fraction $\nicefrac{p}{q} \in \mathbb{Q}\cup \{\infty\}$ corresponds to the slope $ph_1+qh_2 \subset \Sigma$.
If we compute $\mathcal{L}(M_1)$ and $\mathcal{L}(M_2)$ using this convention, then Theorem \ref{teo: L space gluing graph manifold} is equivalent to say that the manifold $M=M_1 \cup_{\varphi} M_2$ is an L-space if and only if
\[
    \mathcal{L}(M_1)^\circ \cup \mathcal{L}(M_2)^\circ = \calS(\Sigma).
\]

The L-space intervals $\mathcal{L}(M_1)$ and $\mathcal{L}(M_2)$ can be found explicitly via \cite[Proposition 3.9]{LspaceIntervalGraphManifold}. By means of calculations and using the convention above, it can be proven that
\[
    \mathcal{L}(M_1)^\circ \cup \mathcal{L}(M_2)^\circ = \calS(\Sigma).
\]
Since this is done via the combinatorial calculation of the L-space interval $\mathcal{L}(M_i)$, the details are not given in this work but will appear in the author's PhD thesis.

\begin{repteo}{teo: M SU(2)-abelian then Lspace}
    Let $M=M_1 \cup_{\Sigma} M_2$ be a $3$-manifold as in Theorem \ref{teo: classification}. If $M$ is an $SU(2)$-abelian rational homology sphere, then $M$ is an L-space.
\end{repteo}
\printbibliography
\end{document}